\documentclass[11pt]{article}

\usepackage[margin=1in]{geometry}
\usepackage{amsmath,amssymb,amsthm,mathtools, mathrsfs, amsfonts}
\usepackage[thmtools-compat]{keytheorems}

\usepackage{bbm, bm}
\usepackage{enumitem}
\usepackage[pdfusetitle,colorlinks=true,linkcolor=blue,citecolor=blue,urlcolor=blue,pagebackref]{hyperref}
\hypersetup{bookmarksnumbered}
\usepackage[capitalise,nameinlink]{cleveref}
\usepackage[dvipsnames]{xcolor}
\usepackage{comment}
\usepackage{microtype}
\usepackage[small]{caption}
\usepackage[normalem]{ulem}
\usepackage{flafter}
\usepackage{placeins}

\DeclareMathOperator{\Core}{Core}

\newcommand{\RedCore}{\mathrm{RedCore}}

\theoremstyle{plain}
\newtheorem{theorem}{Theorem}[section]
\newtheorem{lemma}[theorem]{Lemma}
\newtheorem{proposition}[theorem]{Proposition}
\newtheorem{corollary}[theorem]{Corollary}

\theoremstyle{definition}
\newtheorem{definition}[theorem]{Definition}
\newtheorem{remark}[theorem]{Remark}

\newcommand\numberthis{\addtocounter{equation}{1}\tag{\theequation}}

\newcommand{\E}{\mathbb E}
\newcommand{\R}{\mathbb R}
\newcommand{\N}{\mathbb N}
\newcommand{\1}{\mathbbm 1}
\newcommand{\AMP}{\mathrm{AMP}}
\newcommand{\MMSE}{\mathrm{MMSE}}

\newcommand{\ex}{\mathrm{ex}}

\newcommand{\Var}{\mathrm{Var}}

\newcommand{\ip}[2]{\left\langle #1,#2\right\rangle}

\newcommand{\dd}{\mathrm d}
\newcommand{\eps}{\varepsilon}

\newcommand{\T}{\top}
\newcommand{\RR}{\mathbb{R}}
\newcommand{\EE}{\mathbb{E}}

\newcommand{\MSE}{\mathsf{MSE}}

\newcommand{\In}{\mathbf{1}}

\def\multiset#1#2{\ensuremath{\left(\kern-.3em\left(\genfrac{}{}{0pt}{}{#1}{#2}\right)\kern-.3em\right)}}

\newcommand{\mch}[2]{
\left.\mathchoice
  {\left(\kern-0.48em\binom{#1}{#2}\kern-0.48em\right)}
  {\big(\kern-0.30em\binom{\smash{#1}}{\smash{#2}}\kern-0.30em\big)}
  {\left(\kern-0.30em\binom{\smash{#1}}{\smash{#2}}\kern-0.30em\right)}
  {\left(\kern-0.30em\binom{\smash{#1}}{\smash{#2}}\kern-0.30em\right)}
\right.}
\newcommand{\mcn}{\mch{[n]}{2}}

\newcommand{\calC}{{\cal C}}

\newcommand{\calE}{{\cal E}}
\newcommand{\calF}{{\cal F}}

\newcommand{\calI}{{\cal I}}

\newcommand{\calN}{{\cal N}}

\newcommand{\calS}{{\cal S}}
\newcommand{\calT}{{\cal T}}

\let\cC\calC

\let\cE\calE
\let\cF\calF

\let\cI\calI

\let\cN\calN

\let\cP\calP

\let\cR\calR

\let\cT\calT

\let\cV\calV

\let\bG\boldG

\let\bR\boldR
\let\bS\boldS

\let\bX\boldX
\let\bY\boldY
\let\bZ\boldZ

\let\bh\boldh

\let\bp\boldp

\let\bv\boldv

\let\bx\boldx
\let\by\boldy

\newcommand{\bme}{{\bm e}}
\newcommand{\bmf}{{\bm f}}
\newcommand{\bmg}{{\bm g}}
\newcommand{\bmh}{{\bm h}}

\newcommand{\bmp}{{\bm p}}

\newcommand{\btheta}{{\bm \theta}}
\newcommand{\bmzero}{{\bm 0}}

\newcommand{\fH}{\mathsf{H}}

\newcommand{\sC}{\mathscr{C}}

\newcommand{\sF}{\mathscr{F}}
\newcommand{\sG}{\mathscr{G}}

\newcommand{\sT}{\mathscr{T}}

\newcommand{\lam}{\lambda}

\newcommand{\al}{\alpha}

\newcommand{\gam}{\gamma}

\newcommand{\sig}{\sigma}

\DeclareMathOperator{\He}{He}

\newcommand{\Pb}{\mathbb{P}}
\newcommand{\Qb}{\mathbb{Q}}

\newcommand{\dx}{\mathrm{d}}

\newcommand{\IT}{{\mathrm{IT}}}
\newcommand{\Bayes}{{\mathrm{Bayes}}}

\usepackage{amsmath, amssymb, graphicx}
\usepackage{tikz}
\usetikzlibrary{arrows.meta, calc, decorations.pathreplacing}

\definecolor{ldink}
{HTML}{24323B}
\definecolor{ldteal}
{HTML}{087F8C}
\definecolor{ldorange}
{HTML}{B65C29}
\definecolor{ldpurple}
{HTML}{7853A2}
\definecolor{ldgray}
{HTML}{9AA2A8}
\definecolor{ldpale}
{HTML}{F2F5F6}
\definecolor{ldpaleteal}
{HTML}{E7F4F4}
\definecolor{ldpaleorange}
{HTML}{FBF0E7}

\tikzset{
  ldroot/.style={circle, draw=ldink, fill=ldink, text=white, line width=.7pt, 
    minimum size=5.4mm, inner sep=0pt, font=\small\bfseries}, 
  ldvertex/.style={circle, draw=ldink, fill=white, line width=.8pt, 
    minimum size=4.7mm, inner sep=0pt, font=\small}, 
  ldretained/.style={circle, draw=ldink, fill=white, double=white, double distance=.6pt, 
    line width=.7pt, minimum size=5.4mm, inner sep=0pt, font=\small}, 
  ldsmall/.style={circle, draw=ldink, fill=white, line width=.8pt, 
    minimum size=2.8mm, inner sep=0pt}, 
  ldghost/.style={ldvertex, draw=ldgray, text=ldgray}, 
  ldtreevertex/.style={ldsmall, draw=ldorange, fill=white}, 
  ldedge/.style={draw=ldink, line width=1pt, line cap=round}, 
  ldcore/.style={ldedge}, ldparent/.style={draw=ldteal, line width=1.15pt, -{Stealth[length=2mm, width=1.4mm]}, line cap=round}, 
  ldselected/.style={draw=ldpurple, line width=1.05pt, -{Stealth[length=2mm, width=1.4mm]}}, 
  ldboth/.style={draw=ldpurple, line width=1.05pt, {Stealth[length=2mm, width=1.4mm]}-{Stealth[length=2mm, width=1.4mm]}}, 
  ldunused/.style={draw=ldgray, line width=.8pt, dash pattern=on 3pt off 2.5pt}, 
  ldbranch/.style={draw=ldorange, line width=1pt, -{Stealth[length=1.7mm, width=1.25mm]}}, 
  ldjoin/.style={ldbranch, dash pattern=on 2.5pt off 2pt}, 
  ldarrow/.style={draw=ldink, line width=.8pt, -{Stealth[length=2mm]}}, 
  ldlabel/.style={font=\small, text=ldink}, 
  ldhead/.style={font=\small\bfseries, text=ldink, anchor=west}, 
  ldnote/.style={font=\footnotesize, text=ldink, align=left}, 
}

\title{Degree-$o(n)$ optimality of Approximate Message Passing}
\author{Daniel Fu\thanks{Division of Applied Mathematics, Brown University. \texttt{daniel\_fu1@brown.edu}}\and Chris Jones\thanks{Department of Mathematics, UC Davis. \texttt{chijones@ucdavis.edu}}\and Youngtak Sohn\thanks{Division of Applied Mathematics, Brown University. \texttt{youngtak\_sohn@brown.edu}}\and Alex Wein\thanks{Department of Mathematics, UC Davis. \texttt{aswein@ucdavis.edu}}}
\date{}

\begin{document}
\maketitle

\begin{abstract}   
We study estimation in the spiked Wigner model $Y=\frac{\theta\theta^\top}{\sqrt n}+Z$ where the entries of $\theta\in\mathbb R^n$ are drawn independently from a known prior and $Z$ is symmetric Gaussian noise. For both nonzero-mean and zero-mean priors, we prove that, under a stability condition on the Approximate Message Passing (AMP) fixed point, no polynomial estimator of degree $D=o(n)$ achieves asymptotically smaller mean-squared error than the Bayes AMP. On the other hand, a degree-$O(n)$ estimator can approach the information-theoretically optimal error. The proof builds on the duality framework of Sohn and Wein (2025) by constructing an explicit function-valued dual certificate through a recursive graph expansion inspired by message passing.

We also construct explicit polynomial estimators of degree $O_\varepsilon(\log n)$ whose mean-squared error is within $\varepsilon+o(1)$ of the Bayes AMP error, for every fixed $\varepsilon>0$. Together, these results characterize the optimal low-degree estimation error and provide evidence that the Bayes AMP achieves the best accuracy possible among computationally efficient algorithms.

\end{abstract}
\setcounter{tocdepth}{2}
\tableofcontents

\newpage
\section{Introduction}

A fundamental task in high-dimensional statistics is to determine how accurately a hidden signal can be estimated from noisy observations. Alongside the \textit{information-theoretic} question of identifying the best possible accuracy, achieved by Bayes estimation, the high dimensionality of these problems introduces a \textit{computational challenge}: can the accuracy of the Bayes estimator be attained by an efficient algorithm with practical runtime? For many models of interest, including sparse regression~\cite{celentano2022fundamental, gamarnik2022sparse}, sparse PCA~\cite{BR-reduction,krauthgamer2015semidefinite,lesieur2017constrained, brennan19b}, and Gaussian mixture clustering~\cite{lesieur2016phase, BMVVX, even2024computation, even2025computational}, there are parameter regimes where the best known efficient algorithms fall short of achieving the accuracy of Bayes estimation. Understanding these \emph{statistical--computational gaps}, and determining the best estimation accuracy achievable by efficient algorithms, are central goals in high-dimensional~statistics.

Approximate message passing (AMP)~\cite{donoho2009message} is a family of iterative algorithms for efficient estimation in many high-dimensional statistics models. Building on the ideas from the TAP equations in statistical physics~\cite{thouless1977solution, bolthausen2014iterative}, AMP has an iterative structure (see~\cref{eq:amp}): multiplication by the data matrix alternates with coordinatewise nonlinear denoising, together with an \emph{Onsager correction} that compensates for dependencies introduced by re-using the data matrix. Under suitable assumptions, the estimation accuracy of AMP in the high-dimensional limit is characterized by a deterministic recursion known as \emph{state evolution}~\cite{bayati2011dynamics}. When the signal prior is known, the denoising functions can be chosen as posterior means, giving rise to \emph{Bayes AMP}.

The widespread success of Bayes AMP has led to conjectures that it achieves the best asymptotic estimation accuracy among all polynomial-time algorithms for several regression and low-rank estimation problems
(see e.g.~\cite[Conjecture~1.1]{celentano2022fundamental} and~\cite{lesieur2017constrained}).  However, establishing such conjectures remains out of reach, in part because classical worst-case notions of complexity, such as NP-hardness, do not apply directly to these average-case problems, where estimation accuracy is measured over random inputs drawn from a specified distribution. Instead, a popular and versatile approach is to rigorously study estimation performance within restricted classes of algorithms. In particular, a line of work originating from~\cite{SW-22} studies degree-$O(\log n)$ polynomial estimators as a proxy for polynomial-time algorithms, inspired by the low-degree framework for hypothesis testing problems~\cite{sos-clique,HS17,HKP+17,hopkins-thesis,ld-notes} (see also~\cite{wein2025computational} for a survey).
We will explain this framework in more detail later.

In this work, we use the low-degree framework to investigate the optimality of Bayes AMP. We focus on the rank-one spiked Wigner model
as a basic setting in which to develop our theory.
\begin{definition}[Spiked Wigner model]
\label{def:wig}
For $n\geq 1$ and a prior distribution $\pi$ on $\R$,
let $\btheta\in\R^n$ have i.i.d.\ coordinates $\theta_i\sim\pi$. We observe the symmetric matrix
\[
    \bY=\frac{1}{\sqrt n}\btheta\btheta^\T+\bZ,
\]
where $\bZ=(Z_{ij})_{1\le i\le j\le n}$ is a Gaussian noise independent of $\btheta$ with independent upper triangle entries $(Z_{ij})_{i\leq j}\stackrel{iid}{\sim}\cN(0,1)$. The statistical task is to estimate the signal matrix $\btheta\btheta^\T$ from $\bY$.
\end{definition}

\begin{figure}[!t]
\centering
\includegraphics[width=0.8\linewidth]{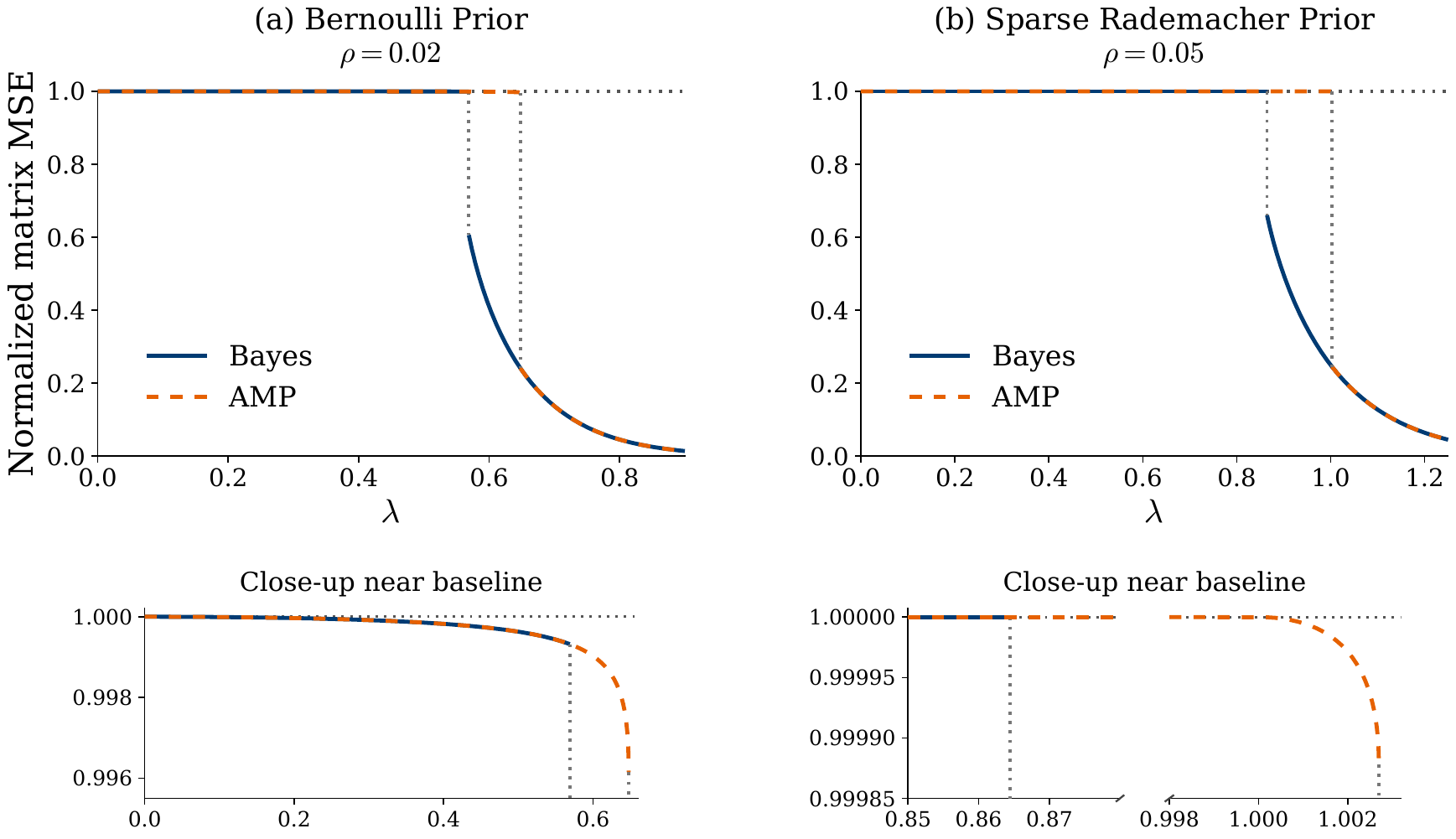}
\caption{$\MSE_{\mathrm{Bayes}}$ (blue solid line) and $\MSE_{\mathrm{AMP}}$ (orange dotted line) for the prior $\pi=\sqrt{\lambda/\rho}B$: (a) $B\sim\mathrm{Bernoulli}(\rho)$ with $\rho=0.02$;
(b) $\mathbb P(B=\pm1)=\rho/2$ and $\mathbb P(B=0)=1-\rho$
with $\rho=0.05$. The horizontal axis is
$\lambda=\mathbb{E}[\Theta^2]$, the signal-to-noise ratio. Both errors are normalized by $\lambda^2-\mu^4$, where
$\mu=\E[\Theta]$, so the trivial estimator
$\E[\btheta\btheta^\T]$ has normalized MSE $1$
(horizontal dotted line). Each curve has a single point of discontinuity, and the two curves match except on the interval between these two ``phase transitions.'' The close-ups illustrate qualitative differences between $\mu \ne 0$ and $\mu = 0$: in (a), the normalized AMP error is strictly below 1 for all $\lambda > 0$; in (b), it is exactly 1 for all $\lambda \le 1$, then improves continuously before jumping to the lower curve at
$\lambda\approx1.0027$. 
}
\label{fig:stat-comp-gap}
\end{figure}

\noindent 
Our regime of interest is $n \to \infty$ with $\pi$ held fixed. We take $\btheta\btheta^\T$ instead of $\btheta$ as the estimand, since it is impossible to distinguish $\btheta$ from $-\btheta$ when $\pi$ is symmetric about zero.
The accuracy of an estimator
$\bmf=(f_{ij})_{1 \leq i,j\leq n}:\R^{n\times n}\to\R^{n\times n}$
can be measured by its \emph{matrix mean-squared error (MSE)}
\begin{align}
    \MSE(\bmf;\btheta\btheta^\T)
    :=\;\frac{1}{n^2} \E \|\bmf(\bY)-\btheta\btheta^\T\|_F^2.
      \label{eq:mat-mse}
\end{align}
We will see in \cref{eq:matrix:scalar:MMSE} that the diagonal entries of $\btheta\btheta^\T$ have a negligible contribution to the error, and it suffices to consider estimation of the scalar $\theta_1 \theta_2$. The \emph{Bayes error} is the MSE achieved by the information-theoretically optimal estimator, the posterior mean
$\bmf_{\mathrm{Bayes}}(\bY)
:=\E[\btheta\btheta^\T\mid\bY]$. We recall the formula for the asymptotic Bayes error
from~\cite{dia2016mutual,lelarge2019fundamental} in order to explain when Bayes AMP falls
short of its accuracy.

For $\Theta\sim\pi$, let $\mu := \E[\Theta], \lambda:=\E[\Theta^2]$ and assume that
$\lambda<\infty$. In our normalization, $\lambda$ measures the signal-to-noise ratio (SNR). Define the
\emph{replica-symmetric potential}~\cite{el2018estimation}
\begin{equation}\label{eq:replica:symmetric:potential}
    \Psi(q;\pi)
    :=\frac{q^2}{4}-\frac{\lambda q}{2}
      +I(\Theta;\sqrt q\,\Theta+G),
    \qquad q\ge0,
\end{equation}
where $G\sim\calN(0,1)$ is independent of $\Theta$ and
$I(\cdot;\cdot)$ denotes mutual information.
If $\lambda=\E[\Theta^2]<\infty$ and $\Psi(\cdot;\pi)$ admits a unique
global minimizer, then it follows from~\cite{dia2016mutual,lelarge2019fundamental} that as $n\to\infty$ with $\pi$ held fixed, 
\[
\MSE(\bmf_{\mathrm{Bayes}};\btheta\btheta^\T)
    \longrightarrow
    \MSE_{\mathrm{Bayes}}
    :=\lambda^2-q_{\mathrm{IT}}^2,
    \qquad
    q_{\mathrm{IT}}
    :=\operatorname*{arg\,min}_{q\ge0}\Psi(q;\pi).
\]
Here, $q_{\mathrm{IT}}$ corresponds to the information-theoretic value, and $\lambda^2$ is the asymptotic error of the zero estimator $\bmf=\bmzero$ (assuming $\E[\Theta^4]<\infty$), so
$q_{\mathrm{IT}}^2$ measures the improvement achieved by Bayes
estimation over this baseline.

In contrast, it follows from~\cite{bayati2011dynamics,MV21} (see \cref{prop:amp-mse}) that Bayes AMP, with an appropriate initialization described in \cref{sec:amp-intro}, achieves asymptotic error
\begin{equation}\label{eq:amp:fixed:point}
    \MSE_\AMP:=\lambda^2-q_*^2,\quad\textnormal{where}\quad q_*:=\min\bigl\{
        q\ge0:\Psi'(q;\pi)=0,\ \Psi''(q;\pi)\ge0
    \bigr\}.
\end{equation}
Here derivatives are taken with respect to $q$, with right derivatives
at $q=0$, and the AMP error is measured in the double limit where $n\to\infty$ and
then the number of iterations goes to
infinity.

For priors such as sparse Bernoulli and sparse Rademacher distributions, there exist parameter regimes in which $q_*<q_{\IT}$~\cite{dia2016mutual, lelarge2019fundamental}. See \cref{fig:stat-comp-gap} and \cref{fig:fe} for an illustration. We study this gap through the low-degree framework mentioned above, which may be viewed as a proxy for efficient computation in high-dimensional statistics.

Let $\RR[\bY]_{\le D}$ be the space of multivariate polynomials with degree at most $D$ in $\bY=(Y_{ij})_{1\le i\le j\le n}$. Consider the
\emph{degree-$D$ minimum mean-squared error}~\cite{SW-22} defined by
\[
    \MMSE_{\le D}(\btheta\btheta^\T) := \inf_{\bmf \in (\RR[\bY]_{\le D})^{n \times n}}\MSE(\bmf;\btheta\btheta^\T)\,.
\]
\noindent Here the infimum ranges over all matrix-valued estimators
$\bmf=(f_{ij})$ whose entries are polynomials of degree at most
$D$, i.e., $f_{ij}\in \R[\bY]_{\leq D}$. Thus, $\MMSE_{\le D}$ is the smallest estimation error
achievable within this class, with the polynomial coefficients
chosen optimally for the model.

Our main result is the following.
\begin{theorem}\label{thm:lowdeg}
Consider the spiked Wigner model in \cref{def:wig} with a
subgaussian prior $\pi$, i.e., $\E_{\Theta\sim\pi}[e^{c\Theta^2}]<\infty$ for some $c>0$.
As $n\to\infty$ with $\pi$ fixed, the following hold.
\begin{enumerate}[label=(\alph*)]
    \item \label{item:lower}
    If $\Psi''(q_*;\pi)>0$, then for every sequence $(D_n)_{n\geq 1}$ with $D_n=o(n)$,
    \[
        \MMSE_{\le D_n}(\btheta\btheta^\T)
        \ge \MSE_\AMP-o(1),
    \]
    where we recall that $\MSE_\AMP=\lambda^2-q_*^2$.
    In the complementary case $\Psi''(q_*;\pi)=0$, the same conclusion holds
    for every $D_n=o(n^{1/5})$. 

    \item \label{item:upper}
    For every $\eps>0$, there exists a constant
    $C=C_{\eps,\pi}>0$ such that
    \[
        \MMSE_{\le D_{n,\eps}}(\btheta\btheta^\T)
        \le \MSE_\AMP+\eps+o(1),
    \]
    where $D_{n,\eps}=C\log n$ if
    $\E_{\Theta\sim\pi}[\Theta]=0$, and
    $D_{n,\eps}=C$ otherwise.
\end{enumerate}
\end{theorem}

\cref{thm:lowdeg} has two components: \ref{item:lower}~a lower bound against all polynomials of degree $o(n)$ under the stability assumption $\Psi''(q_*;\pi)>0$, and \ref{item:upper}~construction of a degree-$O(\log n)$ polynomial that achieves the accuracy of Bayes AMP.
Together, these show that AMP is optimal within the class of low-degree polynomial estimators. The degree $o(n)$ in part \ref{item:lower} is essentially optimal, as we show in \cref{app:bayes} that (for bounded priors at non-exceptional signal strengths) a degree-$O_\eps(n)$ polynomial can achieve the Bayes error up to an additive $\eps > 0$.
The degree in part \ref{item:upper} is also essentially minimal: for mean-zero priors, degree-$o(\log n)$ polynomials fail to achieve any nontrivial estimation~\cite[Theorem~2.6(a)]{sohn2025sharp}.

Based on the ``low-degree heuristic''~\cite[Hypothesis~2.1.5]{hopkins-thesis} (see also~\cite{wein2025computational})---degree-$D$ polynomials are as powerful as time-$n^{\widetilde\Theta(D)}$ algorithms---our lower bound suggests that AMP cannot be surpassed without expending fully exponential runtime $\exp(n^{1-o(1)})$, under the stability condition $\Psi''(q_*; \pi)>0$.
Indeed, the runtime of naively evaluating the Bayes estimator is $\exp(O(n))$,
since it requires a sum or integral over all possible values for $\btheta$.
Conversely, and on a finer scale, the low-degree heuristic posits that polynomial-time algorithms can be captured by degree-$O(\log n)$ polynomials, and our upper bound confirms this for AMP in the spiked Wigner model.
The low-degree heuristic does not hold for all models~\cite{ld-false-1,ld-false-2} but the spiked Wigner model is one of the ``core'' problems where it is most expected to hold, and our work helps corroborate this.

The critical case $\Psi''(q_*;\pi)=0$ in \ref{item:lower} is pathological, encompassing
priors with $\mu = 0$ and $\lambda = 1$ which are exactly on the ``BBP transition''~\cite{BBP,FP-wigner,capitaine2009largest,benaychgeorges2011eigenvalues},
or the discontinuity points of the AMP curve in \cref{fig:stat-comp-gap}.
While we have shown a low-degree lower bound here,
we remain ambivalent about the true computational complexity in the critical case;
the jump discontinuity of $\MSE_\AMP$ as a function of $\lambda$ means the model is not ``noise-robust'' in the way that is often required for low-degree conjectures (see e.g.~\cite[Conjecture~2.2.4]{hopkins-thesis} or~\cite[Section~6.4]{wein2025computational}).
We show in \cref{lem:critical-parameters} that the stability assumption $\Psi''(q_*;\pi)>0$ holds for almost all priors, in the sense that it holds for almost every (in Lebesgue measure) re-scaling of $\lambda$ for any given subgaussian prior.

Both low-degree polynomials and AMP are commonly used to probe the computational complexity of problems in high-dimensional statistics, so it is significant that both frameworks predict the same computational barrier for the spiked Wigner estimation problem, in a sharp sense. On the other hand, the two frameworks are not always equivalent, with \emph{tensor PCA} being a prime example of a model where low-degree polynomials surpass AMP (see~\cite{wein2019kikuchi}). Our work is a starting point for the broader question of understanding when AMP is optimal and when it is not. We refer to the surveys~\cite{wu-xu-survey,ld-notes,disordered-insights,kz-notes,wein2025computational,zadik-notes} for a broader treatment of statistical--computational gaps from various perspectives.

\paragraph{Comparison with previous work.}
Previously, Montanari and Wein~\cite{montanari2022equivalence} 
studied low-degree polynomial estimation for the spiked Wigner model with nonzero-mean priors, using the vector estimand $\btheta$ instead of the matrix $\btheta\btheta^\T$. Define the vector MSE by
\begin{equation}\label{eq:vec-mse}
    \MSE(\bmf;\btheta)
    :=\frac1n\E\|\bmf(\bY)-\btheta\|_2^2
    =\frac1n\sum_{i=1}^n
      \E[(f_i(\bY)-\theta_i)^2].
\end{equation}
For $\E_{\Theta \sim \pi}[\Theta]\ne0$, they showed that, for any $D=O(1)$,
\begin{equation}\label{eq:mw-lower}
\inf_{\bmf\in(\RR[\bY]_{\le D})^n}
    \MSE(\bmf;\btheta)
    \ge \lambda-q_*-o(1),
\end{equation}
where the right-hand side is the MSE achieved by Bayes AMP. They also showed the matching upper bound to within an additive $\eps > 0$ for $D = O_\eps(1)$.
Our results address two open questions posed by~\cite{montanari2022equivalence}:
$D=o(n)$ hardness and mean-zero priors.
First, we extend the inequality~\eqref{eq:mw-lower}
to every $D=o(n)$ under the stability condition
$\Psi''(q_*;\pi)>0$, and to every $D=o(n^{1/5})$ without this condition.
While the results in Theorem~\ref{thm:lowdeg} are stated for matrix MSE, the techniques readily extend to vector MSE (see \cref{thm:lowdeg:m:roots}).
Second, our results cover mean-zero priors,
where now the estimand is $\btheta\btheta^\T$ and the upper bound requires degree $O_\eps(\log n)$ rather than $O_\eps(1)$.

When $\E_{\Theta\sim\pi}[\Theta]=0$, Sohn and
Wein~\cite{sohn2025sharp} showed that, at or below the BBP
threshold $\lambda\le1$, polynomials of degree up to $n^c$
for some constant $c>0$ cannot asymptotically improve on
the trivial MSE.
In this regime, $q_*=0$ and
$\Psi''(0;\pi)=(1-\lambda^2)/2$.
Thus, the stability condition is equivalent to $\lambda<1$.
As an immediate corollary, our theorem extends the hardness result in~\cite{sohn2025sharp} to every $D=o(n)$, when strictly below the
BBP threshold. More substantially, our theorem addresses the $\lambda > 1$ regime in
which AMP achieves nontrivial accuracy (see the close-ups in \cref{fig:stat-comp-gap}), where the challenge is not merely to rule out recovery, but to identify the precise estimation error reachable by low-degree polynomials.

It is worth noting that our proof uses a different approach than \cite{montanari2022equivalence}, instead developing a functional extension of the low-degree duality framework appearing in~\cite{sohn2025sharp,alg-contig,CMSW25,fu2026low,SGSW26,functional-ineq}, which may be of independent
interest. The lower bound constructs a dual certificate whose
coefficients are functions of $\btheta$, and expresses this certificate through a recursive expansion over graph exploration histories inspired by
message passing. We believe this approach may extend much more broadly to
other high-dimensional statistics models in which AMP is conjectured to
achieve optimal estimation accuracy among poly-time algorithms. For the upper bound, we construct an explicit polynomial
estimator and compute its first and second moments directly.
Both our lower and upper bounds use
graph expansions in the spirit of the diagrammatic descriptions of AMP
in~\cite{bayati2015universality,jones2025fourier}.
In contrast to~\cite{montanari2022equivalence}, both our lower and upper bounds are entirely self-contained and do not rely on any prior AMP theory. We give an overview of our proof methods in \cref{sec:overview}.

Recent work by Li~\cite{li2026almost} proves that polynomials of degree $o(n^{1/60})$ cannot surpass AMP in the spiked Wigner model, assuming the prior is subgaussian with nonzero mean $\E_{\Theta\sim \pi}[\Theta]\neq 0$. Again, our results extend the degree range to the optimal $o(n)$ degree under the stability condition and also cover mean-zero priors. Like us,~\cite{li2026almost} extends the duality approach of~\cite{sohn2025sharp}. Our lower bound, which was discovered independently, uses a different dual certificate than the one constructed by~\cite{li2026almost} for the ``cyclic terms,'' as explained in \cref{app:cyclic-terms}.

\subsection{Background on AMP}
\label{sec:amp-intro}
We now describe Bayes AMP for the spiked Wigner model,
including its initialization and the state evolution recursion
for the asymptotic error
$\MSE_\AMP=\lambda^2-q_*^2$~\cite{FR18,DM14,lesieur2017constrained,dia2016mutual,MV21}.
For a broader introduction to AMP and its applications,
see the tutorial~\cite{amp-tutorial}.

\subsubsection{The algorithm}
\label{sec:amp-alg}

Given the observation matrix $\bY \in \R^{n \times n}$ from the spiked Wigner model, the AMP algorithm computes a sequence of estimates $\hat \btheta^{(t)}$ through the recursion:
\begin{align*}\numberthis\label{eq:amp}
    \bx^{(t+1)} &= \frac{1}{\sqrt{n}}\bY \hat\btheta^{(t)} - \langle \eta_t'(\bx^{(t)})\rangle \cdot \hat\btheta^{(t-1)},\\
    \hat \btheta^{(t+1)} &= \eta_{t+1}(\bx^{(t+1)}).
\end{align*}
Here, $\eta_t: \R \to \R$ is a  ``denoising function'' applied coordinatewise, and $\langle{\eta_t'(\bx^{(t)})}\rangle := \frac 1n \sum_{i=1}^n \eta_t'(x^{(t)}_i)$.
Estimate $\btheta\btheta^\T$ using $\hat \btheta^{(T)} (\hat \btheta^{(T)})^\T$ for a large constant number of iterations $T$.

The initialization of AMP can be a bit subtle. A standard practice is to take a small random initialization,
but convergence of AMP from random initialization in $O(\log n)$ iterations has only been proven in special cases~\cite{LFW23}. Rigorous applications of AMP frequently use a warm start depending on the prior.
Recall $\mu:=\E_{\Theta\sim\pi}[\Theta]$ and $\lambda:=\E_{\Theta\sim\pi}[\Theta^2]$.
The AMP algorithm can be initialized as follows, where one case uses a \emph{spectral initialization}~\cite{MV21} involving the leading eigenvector $\bv_1$ of $\bY$:
\begin{equation}\label{eq:amp-init}
    (\hat\btheta^{(0)}, \bx^{(1)}) :=
    \begin{cases}
        (\mu\bm{1}, \frac{\mu}{\sqrt n}\bY\bm{1}), & \mu\ne0,\\
        (\sqrt{n(\lambda-\lambda^{-1})}\bm v_1, \sqrt{n(\lam^3-\lam)}\bm v_1), & \mu=0,\ \lambda>1,\\
        (\bmzero,\bmzero), & \mu=0,\ \lambda\le1.
    \end{cases}
\end{equation}
When $\mu\ne0$, a multiple of the all-ones vector already gives a warm start for $\btheta$.
In the case $\mu = 0$ and $\lambda > 1$, spectral initialization starts from $\bm v_1$ which is a unit-norm top eigenvector of $\bY$, with sign chosen so that it correlates positively\footnote{If $\pi$ is not symmetric about zero, the algorithm can choose the correct sign based on the empirical distribution of $\bm v_1$~\cite[Theorem~1 at $t=0$]{MV21}. If $\pi$ is symmetric about zero, the algorithm can choose either sign and then the sign of $\btheta$ can be chosen accordingly for purposes of analysis.} with $\btheta$.
The spectral initialization is informative about $\btheta$ exactly when $\lambda > 1$,
marking the ``BBP transition''.
When $\mu = 0$ and $\lam \leq 1$, AMP does not achieve nontrivial estimation,
so in this case we use a trivial all-zero initialization.

\subsubsection{State evolution}
\label{sec:se}
A key feature of AMP is that its high-dimensional behavior
is described by \emph{state evolution}~\cite{bayati2011dynamics}. Although our proofs do not rely on state evolution, it motivates several parts of our approach, so we briefly recall it here.

Under suitable regularity conditions on the functions $\eta_t$, the AMP iterates are approximately ``signal plus noise'': the empirical distribution of $(\theta_i,x_i^{(t)})$
converges weakly in probability to the law of $(\Theta, m_t\Theta+ \sqrt{q_t}\, G)$ for certain constants $m_t, q_t$ and $(\Theta, G) \sim \pi \otimes \calN(0,1)$. That is, for every fixed $t$ and every bounded Lipschitz function $\varphi:\R^2\to\R$, as $n\to\infty$,
\begin{equation}\label{eq:signal-noise}
    \frac1n\sum_{i=1}^n\varphi(\theta_i,x_i^{(t)})
  \stackrel{p}{\to}
    \E\left[
        \varphi(\Theta,m_t\Theta+\sqrt{q_t}\,G)
    \right].
\end{equation}
The constants $m_t, q_t$ follow a key recursion, known as the state evolution formula: 
\begin{equation}\label{eq:state-evolution}
    m_{t+1} = \E[\Theta \eta_{t}(m_t \Theta + \sqrt{q_t}\, G)], \qquad \qquad q_{t+1} = \E[\eta_{t}(m_t \Theta + \sqrt{q_t}\, G)^2]\,.
\end{equation}

In the \emph{Bayes AMP} algorithm, the denoiser is chosen to be the posterior estimate of the signal from $\bx^{(t)}$ treated as an observation in this idealized model.
With this denoiser and the initialization from \eqref{eq:amp-init}, it can be shown that the two parameters $m_t$ and $q_t$ are equal for every $t$ (the ``Nishimori property'', see \cref{lem:calibration}).
The \textit{Bayes denoiser} is then $\eta_t(x):=\eta(x;q_t)$ defined by
\begin{equation}\label{eq:denoiser}
    \eta(x;q)
    :=\E[\Theta\mid q\Theta+\sqrt q G=x],\qquad (\Theta,G)\sim \pi\otimes \cN(0,1).
\end{equation}
In the edge case $q=0$, define $\eta(x;0):=\mu=\E[\Theta]$ for all
$x\in\R$. The state evolution sequence $(q_t)_{t\ge1}$ is defined recursively by
\begin{equation}\label{eq:phi}
    q_{t+1}:=\Phi(q_t),
    \qquad
    \Phi(q):=\E[\Theta\eta(q\Theta+\sqrt q\,G;q)]=\E[\eta(q\Theta+\sqrt q\,G;q)^2],
\end{equation}
with initialization
\begin{equation}\label{eq:q-init}
    q_1:=
    \begin{cases}
        \mu^2, & \mu\ne0,\\
        \lambda-\lambda^{-1}, & \mu=0,\ \lambda>1,\\
        0, & \mu=0,\ \lambda\le1.
    \end{cases}
\end{equation}
The sequence $(q_t)_{t\ge1}$ converges to a fixed point of $\Phi$, namely the value $q_*$ defined in \eqref{eq:amp:fixed:point}, and we call $q_*$ the \emph{AMP fixed point}.
It now follows from \cite{bayati2011dynamics,MV21} (see \cref{prop:amp-mse}) that the AMP iteration~\eqref{eq:amp} with Bayes denoising~\eqref{eq:denoiser} and initialization~\eqref{eq:amp-init} achieves the MSE value claimed earlier:
\begin{equation}\label{eq:amp-mse-limit}
    \lim_{n\to\infty}
    \MSE\bigl(\hat\btheta^{(T)}(\hat\btheta^{(T)})^\T;
        \btheta\btheta^\T\bigr)
    =\lambda^2-q_{T+1}^2
    \xrightarrow[T\to\infty]{}\lambda^2-q_*^2=\MSE_\AMP.
\end{equation}

The state evolution recursion also has a useful interpretation in terms of the replica symmetric potential $\Psi(q) = \Psi(q;\pi)$ defined in Eq.~\eqref{eq:replica:symmetric:potential}. Indeed, the I-MMSE relation (see e.g.~\cite{GWSV11}) implies
\begin{equation}\label{eq:Psi-prime}
    \Psi'(q) = \frac{1}{2}(q - \Phi(q))\,.
\end{equation}
Thus, fixed points of $\Phi$ are exactly stationary points of $\Psi$,
and the recursion $q_{t+1}=\Phi(q_t)$ is the gradient descent update on $\Psi$ with step size $2$: $q_{t+1}=q_{t}-2\Psi'(q_t)$.

\begin{figure}[t]
\centering
\includegraphics[width=0.9\linewidth]{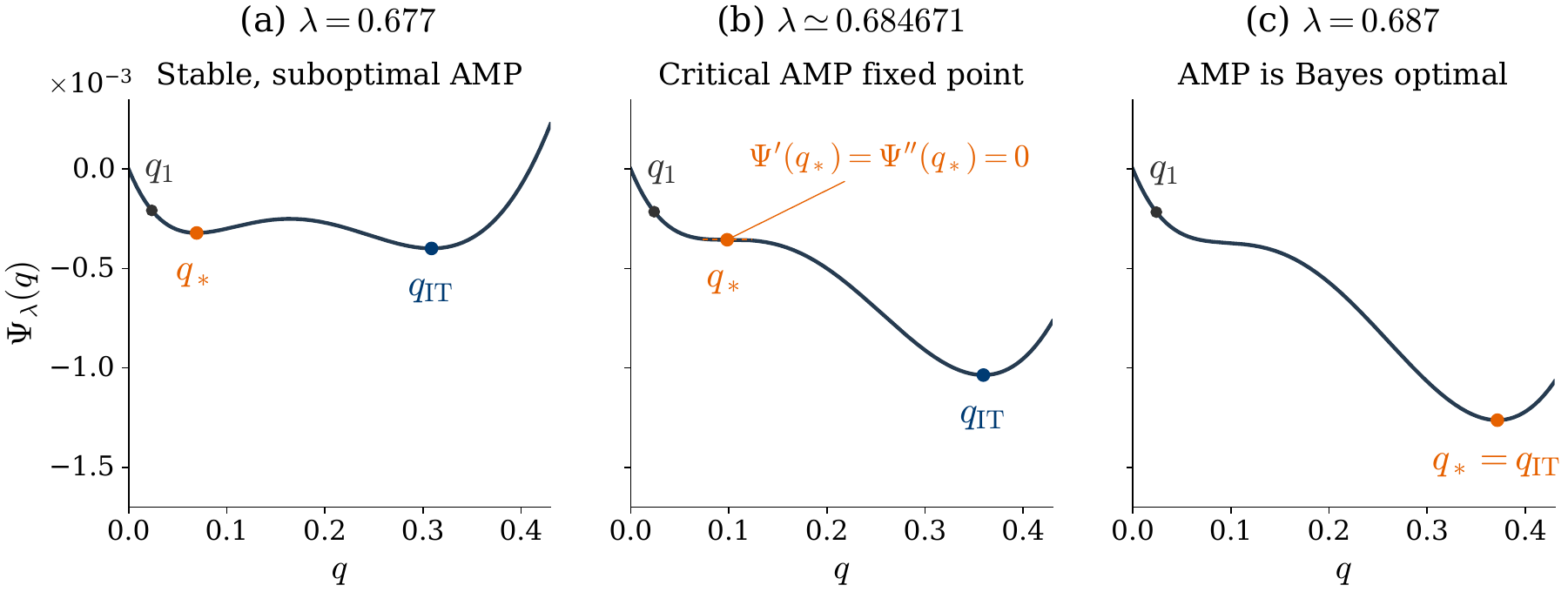}
\caption{Replica symmetric potential $\Psi(q;\pi)$ for the spiked Wigner model with prior
$\sqrt{\lambda/\rho} \cdot \mathrm{Bernoulli}(\rho)$,
where $\rho=0.035$. The three panels show increasing signal strength $\lambda$:
(a) AMP converges to a stable but suboptimal fixed point, with $\Psi'(q_*)=0$ and $\Psi''(q_*)>0$;
(b) the AMP fixed point is critical, meaning
$\Psi'(q_*)=\Psi''(q_*)=0$;
(c) AMP is Bayes optimal.
The state evolution
initialization is $q_1=\lambda\rho$, $q_*$ is its limiting fixed point obtained by moving ``downhill'' starting at $q_1$, and
$q_{\mathrm{IT}}$ is the global minimizer of $\Psi$.}
\label{fig:fe}
\end{figure}

We conclude with some standard facts about AMP's trajectory, which are included in \cref{lem:zero-fixed-point} for completeness. The sequence $q_t$ is nondecreasing,
$\Psi(q_t)$ is nonincreasing, $q_t$ stays in
$[q_1,q_*]$ and it indeed converges to $q_*$
when initialized by~\eqref{eq:q-init}.
\Cref{fig:fe} illustrates this trajectory: starting
at $q_1$, it moves downhill on $\Psi$ and converges to the
first stationary point in $[q_1,\infty)$. More precisely:
\begin{itemize}
    \item If $\mu\ne0$, then $q_*$ is the smallest nonnegative
    fixed point of $\Phi$, and $0<q_1\le q_*$.
    \item If $\mu=0$ and $\lambda>1$, then
    0 is a fixed point of $\Phi$,
    $q_*$ is the smallest positive fixed point
    of $\Phi$, and $0<q_1\le q_*$.
    \item If $\mu=0$ and $\lambda\le1$, then 0 is a fixed point of $\Phi$ and $q_t=q_*=0$
    for every $t$.
\end{itemize}
In all cases, $q_* \leq \lambda < \infty$.

\subsection{Notation}
\label{subsec:notation}

Polynomials in the entries of $\btheta, \bY, \bZ$ will be represented using graphs. We use $\N:=\{0, 1, 2, \ldots\}$ and write $[n]:=\{1, \ldots, n\}$. The space of degree-at-most-$D$ polynomials in $\bY$ is
\[
    \R[\bY]_{\leq D} := \{f \in \R[\bY]: \deg(f) \leq D\}\,.
\]
Let $\mcn = \{(i,j) : 1 \leq i \leq j \leq n\}$ and $\calI = \N^{\mcn}$.
We view $\al \in \calI$ as a multigraph with vertices in $[n]$, where $\alpha_{ij}$ represents the number of copies of edge $\{i,j\}$.

For $\alpha, \beta \in \cI$, define
\[
    |\alpha|:=\sum_{i\le j}\alpha_{ij},\qquad \alpha! := \prod_{i\le j}\alpha_{ij}!,\qquad \binom{\alpha}{\beta} := \prod_{i \leq j} \binom{\alpha_{ij}}{\beta_{ij}}.
\]
We write $\beta \leq \alpha$ whenever $\beta_{ij} \leq \alpha_{ij}$ for all $i \leq j$. For a symmetric matrix $\bX=(X_{ij})_{i,j \in [n]}$, write $\bX^\alpha:=\prod_{i\le j}X_{ij}^{\alpha_{ij}}$; in particular, $(\btheta\btheta^\T)^\alpha=\prod_{i\le j}(\theta_i\theta_j)^{\alpha_{ij}}$.
The degree of a vertex $v \in [n]$ in $\alpha$ is
\[
    \deg_\alpha(v) := \sum_{i \leq v} \alpha_{iv} + \sum_{j \geq v} \alpha_{vj}\,.
\]
Thus, self-loops count twice towards the degree.

Let $\He_k(x; \sigma^2)$ be the probabilists' Hermite polynomials at variance $\sigma^2$ and $h_k(x;\sigma^2) := \frac{\He_k(x; \sigma^2)}{\sqrt{k!}\sig^{k}}$ be their normalization. These polynomials are useful due to their orthogonality with respect to the Gaussian distribution: $\E_{G \sim \cN(0,\sigma^2)}[h_k(G;\sigma^2)h_\ell(G;\sigma^2)]=\In_{\{k=\ell\}}$.
Write $\He_k(x) := \He_k(x; 1)$ and $h_k(x) := h_k(x; 1)$.
For $\alpha\in\calI$, define the multivariate Hermite polynomials by
\[
    \He_\alpha(\bX;\sigma^2)
    :=\prod_{i\le j}\He_{\alpha_{ij}}(X_{ij};\sigma^2),
    \qquad
    \He_\alpha(\bX):=\He_\alpha(\bX;1),
\]
and their normalization at unit variance by
\[
    h_\alpha(\bX)
    :=\frac{\He_\alpha(\bX)}{\sqrt{\alpha!}}
    =\prod_{i\le j}h_{\alpha_{ij}}(X_{ij}).
\]
Since the upper-triangular entries of $\bZ$ are independent standard Gaussians,
$\E[h_\alpha(\bZ)h_\beta(\bZ)]=\In_{\{\alpha=\beta\}}$.

For a random variable $X$ taking values in $\R^d$, define the space
\[
    L^2(X) := \{f: \R^d \rightarrow \R: \text{$f$ is Borel measurable},\, \E[f(X)^2]<\infty\}\,.
\]
Identify $f$ and $g$ if $f(X)=g(X)$ almost surely, and equip $L^2(X)$ with the usual inner product $\langle f, g \rangle=\E[f(X)g(X)]$ and norm $\|f\|=\sqrt{\E[f(X)^2]}$. With these conventions, $L^2(X)$ is a Hilbert space.

\section{Proof overview}
\label{sec:overview}

This section gives a high-level overview of the proof ideas behind our main result, Theorem~\ref{thm:lowdeg}. Starting with the lower bound in part~\ref{item:lower}, in \cref{subsec:duality} we first describe the duality framework and state the requirements for the dual certificate $u$ we need to construct. Then in \cref{subsec:certificate:construction} we describe the construction of $u$, with emphasis on its motivation and derivation from first principles. Finally, \cref{sec:overview-upper} gives an overview of the upper bound in part~\ref{item:upper}, which constructs an explicit polynomial whose MSE matches that of AMP.

\subsection{Duality framework}
\label{subsec:duality}

\subsubsection{Lower bounds via duality}

Consider a general scalar estimand $x \in L^2(\btheta)$,
and recall $\R[\bY]_{\leq D}$ is the space of polynomials of degree at most $D$ in the upper-triangular entries of $\bY$. The (scalar) low-degree MMSE is defined by
\begin{equation}\label{eq:scalar-mmse}
    \MMSE_{\leq D}(x) := \inf_{f \in \R[\bY]_{\leq D}} \E[(f(\bY)-x)^2]\,.
\end{equation}
This can be directly related to the matrix MSE in \cref{eq:mat-mse}: since the joint law of $(\btheta, \bY)$ in the spiked Wigner model is invariant under permutations of $[n]$, 
\begin{equation}\label{eq:matrix:scalar:MMSE}
\begin{aligned}
    \MMSE_{\leq D}(\btheta\btheta^\top)
    &=\frac{1}{n^2} \sum_{i \neq j} \MMSE_{\leq D}(\theta_i\theta_j)+\frac{1}{n^2} \sum_{i=1}^n \MMSE_{\leq D}(\theta_i^2) \\
    &=(1-1/n) \MMSE_{\leq D}(\theta_1\theta_2)+ (1/n)\MMSE_{\leq D}(\theta_1^2) \\
    &=\MMSE_{\leq D}(\theta_1\theta_2)+O(1/n)\,,
\end{aligned}
\end{equation}
where the last step assumes $\E[\theta_1^4]<\infty$.
Therefore, to study the matrix MSE, the relevant scalar estimand is $x = \theta_1 \theta_2$. We let $x$ be a general estimand for now.

The optimal choice of polynomial estimator in \cref{eq:scalar-mmse} is the orthogonal projection in $L^2(\btheta, \bZ)$ of $x$ onto $\R[\bY]_{\leq D}$. To lower bound the MMSE, it therefore suffices to upper bound the squared norm of this projection. However, doing so directly is hard because we lack a convenient orthogonal basis under the law of $\bY$. Instead, we use the duality approach of~\cite{sohn2025sharp}.

\begin{figure}[t]
    \centering
    \begin{tikzpicture}[line width=0.8pt,font=\Large,scale=0.7]
      \coordinate (O) at (2.84,1.14);
      \coordinate (P) at (6.46,1.14);
      \coordinate (X) at (6.46,5.16);
    
\draw (0,0) -- (5.66,0) -- (8.62,2.08) -- (2.98,2.08) -- cycle;
      \node at (3.50,0.50) {$\mathbb{R}[\boldsymbol{Y}]_{\le D}$};
    
\draw (O) -- (P);
      \draw[dash pattern=on 2.4pt off 2.4pt] (X) -- (6.46,-0.96);
      \draw (6.12,1.14) -- (6.12,1.46) -- (6.46,1.46);
    
\draw[line width=1.6pt,line cap=round,line join=round,
            -{Straight Barb[length=6.5pt,width=12pt]}] (O) -- (X);
      \node[above right,inner sep=3pt,font=\LARGE] at (X) {$x$};
    \end{tikzpicture}
    \caption{The subspace of low-degree polynomials in $\bY$ inside $L^2(\btheta,\bZ)$. All points along the dotted line are feasible dual solutions to \cref{eq:dual:constraints}. The minimum-norm solution is the orthogonal projection of $x$ onto $\R[\bY]_{\le D}$, which is also the optimal estimator within $\R[\bY]_{\le D}$.}
    \label{fig:low-degree-projection}
\end{figure}
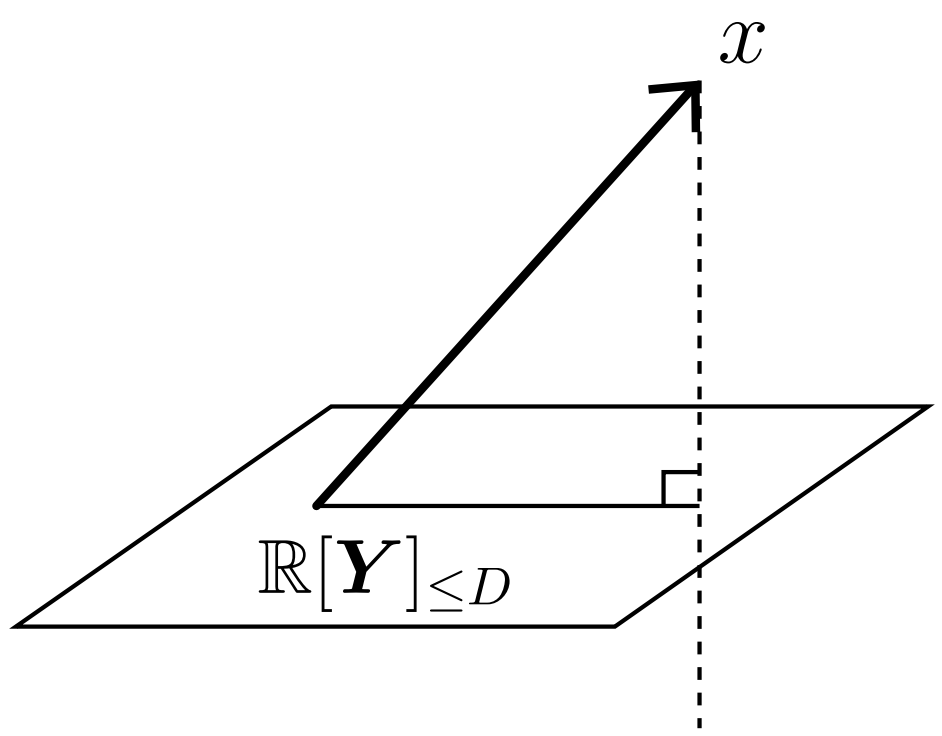

Our strategy is to construct a random variable $u \in L^2(\btheta, \bZ)$ whose projection onto the low-degree space is the same as $x$. Equivalently, we seek $u$ such that
\begin{equation}\label{eq:dual:constraints}
    \E[uf(\bY)]=\E[xf(\bY)]\,, \quad \forall f \in \R[\bY]_{\leq D}\,.
\end{equation}
We call the family of linear equations in Eq.~\eqref{eq:dual:constraints} the \textit{dual constraints}. Any $u$ which satisfies Eq.~\eqref{eq:dual:constraints} is called a \textit{dual certificate} for $x$. A crucial observation is that a dual certificate $u$ need not be a function of $\bY$, which gives us considerable freedom in choosing $u$ as long as it satisfies the dual constraints.

\begin{proposition}\label{prop:MMSE:duality}
If $u \in L^2(\btheta, \bZ)$ satisfies Eq.~\eqref{eq:dual:constraints}, then $\MMSE_{\leq D}(x) \geq \E[x^2]-\E[u^2]$.
\end{proposition}

\begin{proof}
For any $f \in \R[\bY]_{\leq D}$, the dual constraints Eq.~\eqref{eq:dual:constraints} give
\[
\begin{aligned}
    \E[(x-f(\bY))^2] &= \E[x^2]-2\E[xf(\bY)]+\E[f(\bY)^2]   \\
                     &= \E[x^2]-2\E[uf(\bY)]+\E[f(\bY)^2]   \\
                     &= \E[x^2]-\E[u^2]+\E[(u-f(\bY))^2]    \\
                     &\geq \E[x^2]-\E[u^2]\,.
\end{aligned}
\]
Taking the infimum over all $f$ yields the claim.
\end{proof}

To bound the low-degree MMSE, we will construct a dual certificate $u$ and use Proposition~\ref{prop:MMSE:duality}. The challenge is to construct $u$ whose squared norm $\E[u^2]$ is small enough; indeed, $x$ itself is always a valid certificate, but taking $u=x$ only gives the trivial bound $\MMSE_{\leq D}(x) \geq 0$. See \cref{fig:low-degree-projection} for a visualization of the space of valid certificates. The minimum-norm certificate is precisely the optimal degree-$D$ polynomial estimator, but since this is difficult to evaluate directly, we instead seek an explicit, simpler dual certificate whose norm is still asymptotically optimal.

\subsubsection{Dual constraints in the Hermite basis}
\label{subsec:dual:constraints:hermite}
We next specialize \cref{prop:MMSE:duality} to the  Gaussian additive model, where $\bZ$ is symmetric with independent standard normal entries in its upper triangle and is independent of $\btheta$. To this end, we expand the dual certificate $u$ in the Hermite basis with respect to $\bZ$, with coefficients that are functions of $\btheta$. Recall from Section~\ref{subsec:notation} the Hermite polynomials $(h_\alpha)$ orthonormal in $L^2(\bZ)$: $\E[h_\alpha(\bZ)h_\beta(\bZ)] = \In_{\{\alpha=\beta\}}$. We seek a certificate of the form
\begin{equation}\label{eq:u-expansion}
u(\btheta, \bZ)=\sum_{|\al|\leq D} \frac{c_\alpha(\btheta)h_\alpha(\bZ)}{\sqrt{\alpha!}n^{|\alpha|/2}}\,,
\end{equation}
where the coefficients $c_\alpha(\btheta)$ remain to be chosen, and $\sqrt{\alpha!}n^{|\alpha|/2}$ is a normalizing factor that will turn out to
simplify the dual constraints~\eqref{eq:dual:constraints:hermite:functional}.
The advantage of this orthogonal expansion is that the norm of $u$ can be expressed easily:
\begin{equation}\label{eq:u-norm}
\E[u^2]=\sum_{|\al|\leq D} \frac{\E[c_\alpha(\btheta)^2]}{\alpha!n^{|\alpha|}}\,.
\end{equation}
In contrast to Sohn and Wein~\cite{sohn2025sharp}, we will not further expand $c_\alpha(\btheta)$ using an orthogonal basis for $L^2(\btheta)$.

With this reformulation, the new dual variables to be chosen are the functions $\{c_\alpha(\btheta)\}_{|\alpha|\leq D}$. We now express the dual constraints~\eqref{eq:dual:constraints} in these variables. Observe that it is enough to check the dual constraints~\eqref{eq:dual:constraints} for $f$ in a basis of $\R[\bY]_{\leq D}$. We again use the Hermite polynomials; even though these are not orthogonal in $L^2(\bY)$, they will still be convenient. The dual constraints read:
\[
    \E[uh_\alpha(\bY)]=\E[xh_\alpha(\bY)]\,, \quad \forall \alpha \in \calI: |\alpha| \leq D\,.
\]
Writing $\bY = \bX + \bZ$ with $\bX := \btheta \btheta^\T/\sqrt{n}$, we apply the well-known Hermite translation formula (e.g.~\cite[Proposition 3.1]{SW-22}) using the multi-index notation from \cref{subsec:notation},
\[
    h_\alpha(\bX+\bZ)=\sum_{\beta \leq \alpha} \binom{\alpha}{\beta} \frac{\sqrt{\beta!}}{\sqrt{\alpha!}} \bX^{\alpha-\beta}h_\beta(\bZ)\,.
\]
Combining with \cref{eq:u-expansion} and using orthonormality with respect to $\bZ$, the dual constraints \eqref{eq:dual:constraints} can then be equivalently stated as:
\begin{equation}\label{eq:dual:constraints:hermite:functional}
    \sum_{\beta \leq \alpha} \binom{\alpha}{\beta} \E[(\btheta\btheta^\top)^{\alpha-\beta}c_\beta(\btheta)]=\E[x(\btheta\btheta^\top)^\alpha]\,, \quad \forall \alpha: |\alpha| \leq D\,.
\end{equation}
Thus, we have the following corollary of Proposition~\ref{prop:MMSE:duality}.
\begin{corollary}\label{cor:MMSE:lower:bound}
    For any collection $(c_{\alpha})_{|\alpha|\leq D}$ such that $c_{\alpha}\in L^2(\btheta)$ satisfies \eqref{eq:dual:constraints:hermite:functional}, we have
    \[
    \MMSE_{\leq D}(x)\geq \E[x^2]-\sum_{|\al|\leq D} \frac{\E[c_\alpha(\btheta)^2]}{\alpha!n^{|\alpha|}}\,.
    \]
\end{corollary}
More generally, the analogous result holds for a general additive Gaussian model $\bY = \bS/\sqrt{n} + \bZ$ where $(x,\bS)$ has any joint distribution independent of $\bZ$, and with $\bS$ in place of $\btheta\btheta^\T$ in the constraints~\eqref{eq:dual:constraints:hermite:functional} for $c_\alpha(x,\bS)$.

\subsection{Construction of the history certificate}
\label{subsec:certificate:construction}

Specialize to the estimand $x=\btheta_R := \prod_{i=1}^m \theta_i$ in the spiked Wigner model, where $R=\{1, \hdots, m\}$ is a set of ``roots'' and $m \ge 1$ is a fixed integer. The case $m=2$ covers the matrix MSE stated in the introduction, as shown by~\eqref{eq:matrix:scalar:MMSE}. We now motivate our construction of the dual certificate, restricting for simplicity to the case where $\pi$ has nonzero mean and $m=1$ (estimating $\theta_1$). \cref{sec:dual-analysis} will handle the full proof, including the extension to general $m \geq 1$.

\begin{itemize}
    \item Recall that the optimal choice of $u$ is the same as the optimal degree-$D$ estimator. Since we expect Bayes AMP to be a near-optimal low-degree estimator up to $o(1)$ accuracy, this suggests that $u$ should be (approximately) the estimate produced by Bayes AMP, suggesting to take the coefficients $c_\alpha(\btheta)$ to be the Hermite coefficients of Bayes AMP.
    \item The dominant terms in AMP's expansion in the $h_\alpha(\bZ)$ basis are those for which $\alpha$ is a tree rooted at vertex 1. By expanding Bayes AMP, we can systematically derive a formula for $c_\alpha(\btheta)$ when $\alpha$ is a (rooted) tree. This explicit formula is given in \cref{sec:stein-via-amp} with a heuristic derivation in \cref{app:expand-amp-Z}.
    \item These tree coefficients exactly satisfy the dual constraints for trees (i.e., \cref{eq:dual:constraints:hermite:functional} where $\alpha$ is a tree rooted at 1), plus their contribution to the norm in \cref{eq:u-norm} exactly aligns with AMP's MSE. This will be calculated in \cref{subsec:lower:excess-zero}.
    This is already enough to prove the lower bound for $D = O(1)$, as explained at the end of \cref{subsec:lower:excess-zero}.
    \item Pushing the degree to $D=o(n)$ is the major obstacle, and this requires choosing the cyclic coefficients---$c_\alpha(\btheta)$ where $\alpha$ contains a cycle---carefully. A natural attempt is to again inspect the Hermite expansion of AMP, now looking at sub-leading terms. In \cref{app:cyclic-terms} we observe that this fails. For one, these coefficients are not canonical, since there are several variants of message passing that match AMP's asymptotic performance but have different lower-order terms (such as the unapproximated message passing algorithm).
    In all cases we examined, none of the coefficient guesses obtained by this approach satisfy the dual constraints; see \cref{lem:distinct-triangle-coefficients}-\ref{item:dual-fail} for the case of AMP.
    \item We instead take a different approach. Examining why the tree coefficients satisfy the dual constraints leads to a key identity, which we call the \emph{Stein--cumulant identity} (\cref{sec:stein}). A classical form of this identity appears in Barbour's work on Stein's method~\cite{barbour86asymptotic}. Using this identity, we can ``expand'' the quantity $\E[x(\btheta\btheta^\top)^\alpha]$ in \eqref{eq:dual:constraints:hermite:functional} and collect the resulting terms to obtain coefficients $c_\beta(\btheta)$ for cyclic graphs as well. This construction is illustrated via a small example in
\cref{subsec:outline:history:certificate},
and we call the resulting choice of coefficients the \emph{history certificate}.

\end{itemize}

\subsubsection{Tree coefficients suggested by AMP}\label{sec:stein-via-amp}

Let $\cT$ be the set of trees with vertices in $[n]$ rooted at vertex~$1$, including the empty tree $\emptyset$. For a tree $\tau \in \cT$, write $V(\tau)$ for its vertex set, with the convention $V(\emptyset)=\{1\}$. We will now state the formula for tree coefficients $(c_\tau(\btheta))_{\tau \in \cT}$.

Let $(\Theta, G) \sim \pi \otimes \cN(0, 1)$ and recall the Bayes denoiser~\eqref{eq:denoiser}:
\[
    \eta(y; q)=\E[\Theta \mid q\Theta+\sqrt{q}G=y]\,.
\]
The tree coefficients suggested by Bayes AMP are built from the Hermite coefficients of the Bayes denoiser. With $y$ fixed, define the rescaled Hermite coefficients of $\eta(qy+\sqrt{q}G;q)$ by
\[
    A_r(y; q) := \E\left[\eta(qy+\sqrt{q}G; q) \frac{\He_r(G)}{q^{r/2}}\right]\,,\qquad q>0\,,\quad r\geq 0\,.
\]
To motivate the tree coefficients, consider the first coordinate $\hat\theta_1^{(t)}=\eta(x_1^{(t)};q_t)$ of the Bayes AMP iterate~\eqref{eq:amp}. In \cref{app:expand-amp-Z} we calculate that for each fixed tree $\tau$, in the limit $t \to \infty$, the Hermite coefficient of $\hat \btheta^{(t)}_1$ conditional on $\btheta$ is
\begin{equation}\label{eq:amp-hermite}
    \E[\hat\theta_1^{(t)}h_\tau(\bZ)\mid\btheta] \approx n^{-|\tau|/2}\prod_{v\in V(\tau)}A_{k_v}(\theta_v;q_*),
\end{equation}
where $k_v$ is the number of children of $v$ in $\tau$. While we expect \eqref{eq:amp-hermite} to hold in a rigorous sense, our calculation is heuristic; it motivates our choice of coefficients, but it is not used to establish the dual constraints.

While the calculation for \eqref{eq:amp-hermite} extends to general AMP iterations, the Bayes denoiser leads to an additional interpretation of these coefficients in terms of posterior cumulants.
We prove in \cref{lem:Gaussian:field} that $A_r(y;q)=\E[\kappa_{r+1}(\Theta\mid q\Theta+\sqrt q\,G)\mid\Theta=y]$, where $\kappa_{r+1}(\Theta\mid q\Theta+\sqrt q\,G)$ is the $(r+1)$st cumulant of the posterior distribution of $\Theta$ given the Gaussian observation $q\Theta+\sqrt q\,G$. That is, $A_r(y;q)$ averages this posterior cumulant over the observation noise when the signal equals $y$. We extend the definition continuously to $q=0$ by setting $A_r(y;0):=\kappa_{r+1}(\Theta)$.

Although the AMP heuristic suggests to work at the AMP fixed point $q_*$, the dual certificate constructed below will satisfy the dual constraints for any $q \geq 0$ (in fact, for any choice of a ``conditioning field'' defining the conditional cumulants; see Proposition~\ref{prop:history:certifies:dual:constraints}). We will allow $q$ to vary and choose it later to control the norm of the dual certificate. Our choice of $q$ will be near $q_*$, but sometimes slightly larger to buy some ``slack'' in the argument. With this modification, matching \eqref{eq:amp-hermite} to \eqref{eq:u-expansion}, the tree coefficients are
\begin{equation}\label{eq:tree:coefficient}
    c_\tau(\btheta; q) := \prod_{v \in V(\tau)} A_{k_v}(\theta_v; q)\,.
\end{equation}
Thus, the tree coefficients factor as a product over vertices of the tree, a property which will play an important role in controlling the norm of the resulting dual certificate; see~\cref{lem:tree:summation}.

\subsubsection{Dual constraints for stars and the Stein--cumulant identity}
\label{sec:stein}

Consider the tree coefficients stated in~\eqref{eq:tree:coefficient} above. Furthermore, if a graph $\alpha$ has a non-empty component which is isolated from vertex~$1$, we set $c_\alpha \equiv 0$. It turns out that this definition exactly satisfies the constraints~\eqref{eq:dual:constraints:hermite:functional} for which $\alpha$ is a tree rooted at vertex~$1$. We now illustrate this with some small examples, the empty graph and star graphs, which illustrate the main ideas.

\paragraph{Empty graph.}
For $\tau=\emptyset$, we have $c_\emptyset(\btheta; q)=A_0(\theta_1; q)$. The definition of the Bayes denoiser gives
\[
    \E[A_0(\Theta; q)]=\E[\eta(q\Theta+\sqrt{q}G; q)]=\E\Big[ \E[\Theta \mid q\Theta+\sqrt{q}G] \Big]=\E[\Theta]\,.
\]
This implies the dual constraint for the empty graph since $(\Theta, G) \stackrel{d}{=} (\theta_1, G)$.

\paragraph{Star graphs.}
Consider the star graph $\tau=\{e_i\}_{i=1}^d$ with edges $e_i=\{1, v_i\}$ and non-root vertices $v_1, \hdots, v_d \in [n] \setminus \{1\}$. The corresponding dual constraint reads
\[
    \sum_{\tau' \leq \tau} \E[c_{\tau'}(\btheta; q)(\btheta\btheta^\top)^{\tau-\tau'}]=\E[\theta_1(\btheta\btheta^\top)^\tau]\,,
\]
Using independence of the prior, $\E[\theta_1(\btheta\btheta^\top)^\tau]=\E[\theta_1^{d+1}] \prod_{i=1}^d \E[\theta_{v_i}]=\E[\Theta^{d+1}]\,\E[\Theta]^d$. To compute the left-hand side, note that any
subtree has the form $\tau'=\{e_i: i \in S\}$, where $S \subseteq [d]$. The definition of the tree coefficients~\eqref{eq:tree:coefficient} and independence of the prior give
\[
    \E[c_{\tau'}(\btheta; q)(\btheta\btheta^\top)^{\tau-\tau'}]=\E\left[A_{|S|}(\theta_1; q) \prod_{i \in S} A_{0}(\theta_{v_i}; q) \prod_{i \not \in S} \theta_1\theta_{v_i} \right] =\E\left[A_{|S|}(\theta_1; q)\theta_1^{d-|S|}\right] \E[\Theta]^d.
\]
In the last equality, we used that $\E[A_0(\theta_{v_i}; q)]=\E[\theta_{v_i}]=\E[\Theta]$. Note that the final quantity in the above display depends only on $|S|$, the size of the subset. Since there are $\binom{d}{r}$ possible subsets of $[d]$ of size $r$, taking summation over $S \subseteq [d]$ implies that the dual constraint is equivalent to
\[
    \sum_{r=0}^d \binom{d}{r} \E\left[A_{r}(\theta_1; q)\theta_1^{d-r}\right]\E[\Theta]^d=\E[\Theta^{d+1}]\,\E[\Theta]^d\,.
\]
It therefore suffices to establish the following identity:
\begin{equation}\label{eq:Stein:cumulant:star:q}
    \sum_{r=0}^d \binom{d}{r} \E\left[A_{r}(\Theta; q)\Theta^{d-r}\right]=\E[\Theta^{d+1}]\,, \quad \Theta \sim \pi\,.
\end{equation}
In summary, in order for the heuristic choice of tree coefficients in~\eqref{eq:tree:coefficient} to produce a valid dual certificate, we must verify that Eq.~\eqref{eq:Stein:cumulant:star:q} holds for all integers $d \geq 0$. Perhaps surprisingly, this identity holds generally with $A_r(\Theta;q)$ replaced by $\E[\kappa_{r+1}(\Theta\mid\sF)\mid\Theta]$ for any conditioning $\sigma$-algebra $\sF$, where $\kappa_{r+1}(\Theta\mid\sF)$ denotes the $(r+1)$-st cumulant of the conditional distribution of $\Theta$ given $\sF$ (see the formal definition in~\eqref{eq:def:cumulant}). We prove this general identity in \cref{lem:scalar:Stein:cumulant:star}.

The coefficients in~\eqref{eq:Stein:cumulant:star:q} correspond to the Gaussian observation $\sF=\sigma(q\Theta+\sqrt q\,G)$: as shown in \cref{lem:Gaussian:field}, $A_r(y;q)=\E[\kappa_{r+1}(\Theta\mid\sF)\mid\Theta=y]$. A classical, unconditional form of~\eqref{eq:Stein:cumulant:star:q} appears in Barbour's work on Stein's method~\cite{barbour86asymptotic}. For this reason, we call the identity~\eqref{eq:Stein:cumulant:star:q} the \textit{Stein--cumulant identity}. This identity is the basis of the remaining certificate construction.

\subsubsection{History certificate}
\label{subsec:outline:history:certificate}

The Stein--cumulant identity~\eqref{eq:Stein:cumulant:star:q} can be used inductively to check that the proposed tree coefficients~\eqref{eq:tree:coefficient} satisfy the dual constraints for trees. Furthermore, we will show how it also suggests a construction of $c_\alpha(\btheta)$ for general multigraphs $\alpha$, whose tree coefficients agree with~\eqref{eq:tree:coefficient}. The idea is to exploit the independence structure of the prior and use~\eqref{eq:Stein:cumulant:star:q} at each coordinate separately.

Let $\alpha$ be a multigraph and write $E(\alpha)$ for its multiset of edges. For the moment, assume that $\alpha$ is simple, i.e., has no parallel edges or self-loops. Let
\[
    \alpha_1=\{e\in E(\alpha):1\in e\}
\]
be the star graph consisting of all edges incident to vertex $1$. Conditioning on $\theta_2,\ldots,\theta_n$ and using the identity~\eqref{eq:Stein:cumulant:star:q} with $\Theta=\theta_1$ and $d=\deg_\alpha(1)$ gives
\begin{equation}
\label{eq:history:expansion:A}
\begin{aligned}
    \E[\theta_1(\btheta\btheta^\top)^\alpha]
    &=
    \E\left[
        \theta_1^{\deg_\alpha(1)+1}
        \prod_{v \neq 1}
            \theta_v^{\deg_\alpha(v)}
    \right] \\
    &=
    \sum_{r=0}^{\deg_\alpha(1)}
        \binom{\deg_\alpha(1)}{r}
        \E\left[
            A_r(\theta_1; q)\theta_1^{\deg_\alpha(1)-r}
            \prod_{v \neq 1}
                \theta_v^{\deg_\alpha(v)}
        \right]\,.
\end{aligned}
\end{equation}
We can interpret the sum on the RHS of this identity as choosing subgraphs of the star $\alpha_1$. There are $\binom{\deg_\alpha(1)}{r}$ choices of $\beta_1\leq\alpha_1$ with $|\beta_1|=r$. For each such choice,
\[
    \theta_1^{\deg_\alpha(1)-r}
    \prod_{v \neq 1}
        \theta_v^{\deg_\alpha(v)}
    =
    \left(\prod_{\{1,v\}\in\beta_1}\theta_v\right)
    (\btheta\btheta^\top)^{\alpha-\beta_1}\,.
\]
Thus~\eqref{eq:history:expansion:A} becomes
\begin{equation}
\label{eq:history:expansion:B}
    \E[\theta_1(\btheta\btheta^\top)^\alpha]
    =
    \sum_{\beta_1\leq\alpha_1}
    \E\left[
        A_{|\beta_1|}(\theta_1; q)
        \left(\prod_{\{1,v\}\in\beta_1}\theta_v\right)
        (\btheta\btheta^\top)^{\alpha-\beta_1}
    \right]\,.
\end{equation}
For each choice of $\beta_1$, we can repeatedly expand in the same fashion using the Stein--cumulant identity~\eqref{eq:Stein:cumulant:star:q}. To further illustrate this expansion, we work through the example of a triangle. This will also motivate the ``search procedure'' used in the definition of history certificate in~\cref{sec:dual-analysis}.

\paragraph{Expanding a triangle.}
Let $\triangle=\{12, 13, 23\}$, where we abbreviate $ij$ for the edge $\{i, j\}$. As shorthand, write $A_{r, i}:=A_r(\theta_i; q)$. We illustrate the expansion of
\[
    \E[\theta_1(\btheta\btheta^\top)^\triangle]
\]
by repeated applications of the Stein--cumulant identity~\eqref{eq:Stein:cumulant:star:q}, where the binomial coefficient $\binom{d}{r}$ is interpreted as enumerating the possible ways to select $r$ of the $d$ edges incident to a vertex.

We record these choices through a graph exploration process starting from the root $1$. Selecting an edge $ij$ at vertex $i$ is represented by an arrow $i\to j$, which discovers vertex $j$ if it has not already been discovered. At each step, ``explore'' the smallest discovered vertex that has not yet been explored. Each discovered vertex is explored once. The exploration process stops when every discovered vertex has been explored. The resulting arrows form a \emph{history} on the selected subgraph. The precise rules, including when an edge may be bidirected, are given below.
Although the definition of history relies on the order in which the vertices are explored, the resulting certificate is in fact invariant of this ordering.
Since the choice of ordering does not affect our analysis, we omit the proof of this fact.

After completing the exploration, we group the resulting terms by their selected edge set $\beta\leq\triangle$. Different histories with the same $\beta$ contribute individual summands to the coefficient $c_\beta(\btheta;q)$.

\begin{itemize}
    \item Starting from vertex~$1$, the selected incident edges can be $\emptyset$, $\{12\}$, $\{13\}$, or $\{12, 13\}$. If no edge is selected, then
no other vertex is reached, and the history stops here. This case corresponds to the term
    \[
        \E[A_{0, 1}(\btheta\btheta^\top)^\triangle]\,.
    \]
    Here, no edge is oriented.
    \item \textbf{One edge selected at the root.} Suppose that $\{12\}$ is selected at vertex~$1$. Orient $12$ as $1 \to 2$ so that vertex~$2$ is \emph{discovered} and $12$ is distinguished as a \emph{parent edge}. Among the non-parent edges incident to vertex~$2$, we can select either $\emptyset$ or $\{23\}$. The former choice terminates the history and the latter orients $23$ as $2 \to 3$, discovering vertex~$3$ with parent edge $23$. Using the same rules for vertex~$3$, we get the total contribution
    \[
        \E[A_{1, 1}A_{0, 2}(\btheta\btheta^\top)^{\triangle-\{12\}}]+\E[A_{1, 1}A_{1, 2}A_{0, 3}(\btheta\btheta^\top)^{\triangle-\{12, 23\}}]+\E[A_{1, 1}A_{1, 2}A_{1, 3}\theta_1]\,.
    \]
    The three terms correspond to the orientations $1 \to 2$, $1 \to 2 \to 3$, and $1 \to 2 \to 3 \to 1$, in that order.

    If $\{13\}$ is selected at the root, then the same reasoning applies as the case $\{12\}$, just with the vertices~$2$ and~$3$ exchanged. This case gives
    \[
        \E[A_{1, 1}A_{0, 3}(\btheta\btheta^\top)^{\triangle-\{13\}}]+\E[A_{1, 1}A_{0, 2}A_{1, 3}(\btheta\btheta^\top)^{\triangle-\{13, 23\}}]+\E[A_{1, 1}A_{1, 2}A_{1, 3}\theta_1]\,.
    \]
    The three terms correspond to $1 \to 3$, $1 \to 3 \to 2$, and $1 \to 3 \to 2 \to 1$, in that order.
\item \textbf{Two edges selected at the root.} If $\{12, 13\}$ is selected, orient $12$ as $1 \to 2$ and $13$ as $1 \to 3$, discovering vertices~$2$ and~$3$, with parent edges $12$ and $13$, respectively.
At each non-root vertex, we may select either $\emptyset$ or $\{23\}$, giving the total contribution
    \[
        \E[A_{2, 1}A_{0, 2}A_{0, 3}(\btheta\btheta^\top)^{\triangle-\{12, 13\}}]+\E[A_{2, 1}A_{1, 2}A_{0, 3}\theta_3]+\E[A_{2, 1}A_{0, 2}A_{1, 3}\theta_2]+\E[A_{2, 1}A_{1, 2}A_{1, 3}]\,.
    \]
    These four terms correspond to the orientations with $1 \to 2$, $1 \to 3$, and where the edge $23$ is unoriented, or oriented as $2 \rightarrow 3$, $2 \leftarrow 3$, or $2 \leftrightarrow 3$, in that order.
\end{itemize}
We now group histories by their shape, i.e., the multigraph consisting of its selected edges. Write $c_\beta \equiv c_\beta(\btheta; q)$.
\begin{itemize}
    \item \textbf{Trees.} For the empty graph, we have $c_\emptyset=A_{0, 1}$. For the two single-edge trees $\{12\}$, $\{13\}$, we have $c_{\{12\}}=A_{1, 1}A_{0, 2}$ and $c_{\{13\}}=A_{1, 1}A_{0, 3}$.

    Next, there are three trees with two edges: $\{12, 23\}, \{13, 23\}$, and $\{12, 13\}$. For these, we have
    \[
        c_{\{12, 23\}}=A_{1, 1}A_{1, 2}A_{0, 3}\,, \quad c_{\{13, 23\}}=A_{1, 1}A_{0, 2}A_{1, 3}\,, \quad c_{\{12, 13\}}=A_{2, 1}A_{0, 2}A_{0, 3}\,.
    \]
    \item \textbf{Triangle.} Finally, the remaining five orientations correspond to the triangle graph. Collecting these together gives
    \[
        c_\triangle=2A_{1, 1}A_{1, 2}A_{1, 3}\theta_1+A_{2, 1}A_{1, 2}A_{0, 3}\theta_3+A_{2, 1}A_{0, 2}A_{1, 3}\theta_2+A_{2, 1}A_{1, 2}A_{1, 3}\,.
    \]
\end{itemize}
By our construction, it holds that
\[
    \E[\theta_1(\btheta\btheta^\top)^\triangle]=\sum_{\beta \leq \triangle} \E[c_\beta (\btheta\btheta^\top)^{\triangle-\beta}]\,.
\]
Notice that the tree coefficients agree with the heuristic proposed in~\eqref{eq:tree:coefficient}.

\begin{remark}[Consistency of the history certificate]
A crucial observation is that the coefficients obtained from the triangle calculation are not specific to that example. Suppose we replace $\triangle$ by another graph $\alpha$ such that $\beta\leq \alpha$ and collect the terms whose selected edges form $\beta$. Then, we obtain the same coefficient $c_\beta(\btheta;q)$ as if we had started from $\beta$ itself. Indeed, the selected edges produce the same histories, and selecting $r$ edges at vertex $v$ contributes $A_r(\theta_v;q)$ in either expansion. The unselected edges contribute only to $(\btheta\btheta^\top)^{\alpha-\beta}$. In this sense, the coefficients are \emph{intrinsic} to the selected graph, allowing us to use the same family in every dual constraint. Together with the Stein–cumulant identity at each step of the expansion, this consistency ensures that the resulting coefficients satisfy all the dual constraints.

\end{remark}

\subsection{Upper bound}
\label{sec:overview-upper}

We now give an overview of the upper bound, \cref{thm:lowdeg}-\ref{item:upper}. One potential approach would be to directly approximate the AMP iteration by replacing the denoisers by polynomials, and then appeal to the known analysis of AMP state evolution. Indeed, this has already been carried out by~\cite{montanari2022equivalence} for the vector MSE~\eqref{eq:vec-mse} when $\E_{\Theta\sim \pi}[\Theta]\neq 0$. A similar approach should be possible for zero-mean priors, although the spectral initialization makes this more complicated. We instead give a direct and self-contained proof where we construct an explicit polynomial and analyze its first two moments, without relying on any prior work on state evolution.

\begin{figure}[ht]
    \centering
    \includegraphics[width=\textwidth]{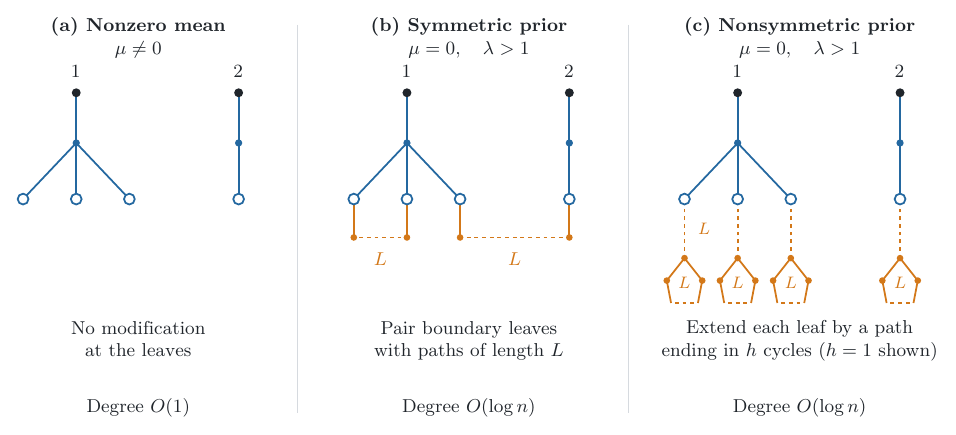}
    \caption{
        The three constructions used in the upper bound, where $\mu:=\E_{\Theta\sim\pi}[\Theta]$ and $\lambda:=\E_{\Theta\sim\pi}[\Theta^2]$.
        Blue edges form the rooted trees from AMP's expansion, and
        circles mark their boundary leaves.
        (a) For nonzero mean, no additional modification is needed (aside from a scalar factor).
        (b) For a prior that is symmetric about zero, pair the boundary
        vertices arbitrarily and connect each pair by a path of length $L=O(\log n)$.
        (c) For a non-symmetric prior with zero mean, attach to each
        boundary vertex a path of length $L$ ending at a vertex
        shared by $h$ cycles of length $L$, where
        $\mathbb{E}[\Theta^{2h+1}]\ne0$; the figure shows $h=1$.
        Dashed segments indicate omitted intermediate vertices.
    }
    \label{fig:upper-bound-constructions}
\end{figure}

To construct the polynomial, first consider the case where the estimand is $\theta_1$, and expand AMP's output as a polynomial in $\bY$. In \cref{sec:Y-basis} we give a heuristic derivation of the dominant terms in this expansion, following~\cite{jones2025fourier}. In particular, the dominant terms are $\bY^\alpha$ where $\alpha$ is a tree rooted at vertex 1, and we give a formula for their coefficients.
To instead estimate $\theta_1 \theta_2$, the expansion is similar, now dominated by $\bY^\alpha$ where $\alpha$ is a forest rooted at vertices 1 and 2.

There are three separate cases for the construction, depending on the prior, and these are depicted in \Cref{fig:upper-bound-constructions}. In all cases, the trees have constant depth and a coefficient on each internal vertex of the tree suggested by AMP. The leaf vertices correspond to the AMP initialization, where a suitable polynomial must be used as a replacement for the initialization:
\begin{enumerate}[label=(\alph*)]
    \item When $\mu \neq 0$, put the constant $\mu$ at each leaf to simulate the initialization $\mu\mathbf{1}$.
    \item When $\mu = 0$ and the prior is symmetric about zero, simulate the spectral initialization by attaching a path of length $O(\log n)$ to each leaf, akin to power iteration.
    The endpoints of the paths are furthermore paired together, which is similar to using a random initialization for the power iteration, then averaging over the random initialization with the Wick lemma.
    \item When $\mu = 0$ and the prior is not symmetric about zero, the spectral initialization furthermore should have its sign correlated with $\btheta$.
    To achieve this, a collection of cycles is placed at the endpoint of each path, which can be interpreted as giving the power iteration an initialization that correlates positively with the signal.
\end{enumerate}

For cases~(a) and~(c), where the prior is not symmetric, the construction generalizes to any number $m \ge 1$ of roots, i.e., estimating $\theta_1 \cdots \theta_m$.
Case~(b), where the prior is symmetric, requires an even number of roots (each tree has an odd number of leaves and these need to be paired). This reflects how nontrivial estimation is impossible for odd $m$ in case (b), due to the symmetry between $\btheta$ and $-\btheta$. In all cases, for simplicity we will only present the formal proof for $m=2$.

The proof involves explicit calculations of the first two moments of the polynomial. These calculations elucidate the connection between the combinatorics of AMP's expansion into trees and the state evolution formulas. The details can be found in \cref{sec:upper-bound}.

\newcommand{\he}{\mathrm{height}}

\section{Lower bound}
\label{sec:dual-analysis}
In this section, we prove the following low-degree lower bound for estimating $\btheta_R=\theta_1\theta_2 \cdots \theta_m$.
\begin{theorem}
\label{thm:lowdeg:m:roots}
Consider the spiked Wigner model in \cref{def:wig} with a fixed subgaussian prior~$\pi$. Fix $m \geq 1$ and let $R=\{1, \hdots, m\}$. The following hold as $n \rightarrow \infty$.
\begin{enumerate}[label=(\roman*)]
    \item Suppose that $\Psi''(q_*;\pi)>0$. Then, for every sequence $(D_n)_{n\geq 1}$ with $D_n=o(n)$, 
    \[
        \MMSE_{\leq D_n}(\btheta_R) \geq \E[\btheta_R^2]-q_*^m-o(1)\,.
    \]
    \item In the complementary case $\Psi''(q_*;\pi)=0$, the same conclusion holds for $D_n=o(n^{1/5})$.
\end{enumerate}
\end{theorem}
As a consequence of~\cref{thm:lowdeg:m:roots}, we immediately deduce~\cref{thm:lowdeg}-\ref{item:lower}.
\begin{proof}[Proof of \cref{thm:lowdeg}-\ref{item:lower}]
By~\eqref{eq:matrix:scalar:MMSE}, for any sequence $(D_n)_{n \geq 1}$, we have
\[
    \MMSE_{\leq D_n}(\btheta\btheta^\top)=\MMSE_{\leq D_n}(\theta_1\theta_2)+o(1)\,.
\]
The claimed lower bound then follows directly from~\cref{thm:lowdeg:m:roots} with $m=2$.
\end{proof}
We will prove Theorem~\ref{thm:lowdeg:m:roots} using the duality framework outlined in Section~\ref{subsec:duality} and the history certificate illustrated in Section~\ref{subsec:certificate:construction}. Section~\ref{subsec:lower:preliminaries} formally defines the history certificate and verifies the dual constraints. Section~\ref{subsec:lower:positive-excess} bounds the norm of the history coefficients, separating out the leading-order forest terms from the remaining negligible positive-excess terms. Finally, Section~\ref{sec:lowdeg:proof} combines these bounds to prove Theorem~\ref{thm:lowdeg:m:roots}.

Throughout this section, $E(\alpha)$ denotes the multiset of edges of a multigraph $\alpha$, and $V(\alpha)$ consists of the non-isolated vertices together with all roots in $R$. We call $\alpha$ \emph{root-connected} if each connected component contains at least one root; a component may consist of an isolated root.

\subsection{History coefficients and history certificate}
\label{subsec:lower:preliminaries}

We first prove a general version of the Stein--cumulant identity~\eqref{eq:Stein:cumulant:star:q} and use it to establish Proposition~\ref{prop:history:certifies:dual:constraints}. Let $\Theta\sim\pi$ be defined on a probability space $(\Omega,\sG,\Pb)$, and $\sF\subseteq\sG$ be a sub-$\sigma$-algebra. One may think of $\sF=\sigma(\xi)$, where $\xi$ is an auxiliary observation of $\Theta$. Our main application is the Gaussian channel $\xi=q\Theta+\sqrt q\,G$, where $G\sim\cN(0,1)$ is independent of $\Theta$.

For random variables $X_1, \hdots, X_k$ on $\Omega$, define their \textit{joint cumulant given $\sF$} by
\begin{equation}\label{eq:def:cumulant}
    \kappa(X_1, \hdots, X_k \mid \sF) := \sum_{\sigma \in \cP([k])} (-1)^{|\sigma|-1}(|\sigma|-1)! \prod_{I \in \sigma} \E\Big[\prod_{i \in I} X_i \,\Big|\, \sF \Big]\,,
\end{equation}
where $\sigma \in \cP([k])$ is a partition of $[k]=\{1, \hdots, k\}$ and $|\sigma|$ is the number of blocks in a partition. In the case that $X_1=\hdots=X_r=\Theta$, define
\begin{equation}\label{eq:cumulant:k:copies}
    \kappa_r(\Theta \mid \sF) := \kappa(\underbrace{\Theta, \hdots, \Theta}_{\text{$r$ copies}} \mid \sF)\,, \quad A_r(y; \sF)=\E[\kappa_{r+1}(\Theta \mid \sF) \mid \Theta=y]\,.
\end{equation}
The following lemma is a conditional version of a classical identity used in normal approximation literature~\cite{barbour86asymptotic,nourdin2010cumulants}. We include a proof for completeness.
\begin{lemma}[Stein--cumulant identity]
\label{lem:scalar:Stein:cumulant:star}
Let $d \geq 0$ be an integer. If $\E[|\Theta|^{d+1}]<\infty$, then
\begin{equation}\label{eq:Stein:cumulant:star:F}
    \sum_{r=0}^d \binom{d}{r} \E\left[A_{r}(\Theta; \sF)\Theta^{d-r}\right]=\E[\Theta^{d+1}]\,.
\end{equation}
\end{lemma}

\begin{proof}
We use the moment--cumulant formula (or equivalently, Möbius inversion on the lattice of set partitions of $[k]$); see, e.g., the display preceding~(2.9) in
\cite[Section~2.3.4, Pg.~38]{mccullagh2018tensor}. For random variables $X_1, \hdots, X_k$, with $\sup_{I \subseteq [k]} \E[\prod_{i \in I} |X_i|]<\infty$, the moment--cumulant formula asserts that
\begin{equation}\label{eq:moment:cumulant:conditional}
    \E\left[\prod_{i=1}^k X_i \,\middle|\, \sF\right]
    =
    \sum_{\sigma\in\cP([k])}
    \prod_{I\in\sigma}
    \kappa\left((X_i)_{i\in I}\mid\sF\right)\,.
\end{equation}
Observe that~\eqref{eq:moment:cumulant:conditional} implies the identity
\begin{equation}\label{eq:moment:cumulant:conditional:block}
    \E\left[\prod_{i=1}^k X_i \,\middle|\, \sF\right]=\sum_{\substack{I\subseteq[k]\\1\in I}}
    \E\left[ \kappa\left( (X_i)_{i\in I}\mid\sF \right)
        \prod_{j\in[k]\setminus I}X_j
        \,\middle|\, \sF
    \right]\,.
\end{equation}
To prove~\eqref{eq:Stein:cumulant:star:F}, use~\eqref{eq:moment:cumulant:conditional:block} with $X_1=\hdots=X_{d+1}=\Theta$. The LHS is $\E[\Theta^{d+1} \mid \sF]$. On the RHS, there are $\binom{d}{r}$ ways to choose the block $I$ such that $1 \in I$ and $|I \setminus \{1\}|=r$.
The cumulant corresponding to this block is $\kappa((X_{i})_{i \in I} \mid \sF)=\kappa_{r+1}(\Theta \mid \sF)$ and the product outside of this block is $\Theta^{d-r}$. Thus, 
\[
    \E[\Theta^{d+1} \mid \sF]=\sum_{r=0}^d \binom{d}{r} \E[\kappa_{r+1}(\Theta \mid \sF) \Theta^{d-r} \mid \sF]\,.
\]
Taking expectations on both sides, the tower property gives
\[
    \E[\Theta^{d+1}]=\sum_{r=0}^d \binom{d}{r} \E\left[ \kappa_{r+1}(\Theta \mid \sF)\Theta^{d-r} \right]=\sum_{r=0}^d \binom{d}{r} \E[A_{r}(\Theta; \sF)\Theta^{d-r}]\,,
\]
where the last equality follows by conditioning on $\Theta$ and inserting the definition of $A_r(\cdot; \sF)$. This completes the proof.
\end{proof}

The next lemma shows that the coefficients $A_r(y;q)$ defined in~\cref{subsec:certificate:construction} are a special case of~\eqref{eq:cumulant:k:copies}. We defer its proof to Appendix~\ref{subsec:omitted:lower:proofs}. Let $\sF_q=\sigma(q\Theta+\sqrt{q}G)$, where $G \sim \cN(0, 1)$ is independent of $\Theta$. Note that $A_r(y; 0)=A_r(y; \sF_0)=\kappa_{r+1}(\Theta)$ by definition.

\begin{lemma}\label{lem:Gaussian:field}
For all $q>0$, integers $r \geq 0$, and $\pi$-almost all $y \in \R$ it holds that
\[
    \E\left[\eta(qy+\sqrt{q}G; q) \frac{\He_r(G)}{q^{r/2}}\right]=\E[\kappa_{r+1}(\Theta \mid \sF_q) \mid \Theta=y]\,.
\]
\end{lemma}

We now formally define the history certificate illustrated by the triangle example in~\cref{subsec:outline:history:certificate}. We begin by specifying the allowed edge orientations and the corresponding exploration process. To state the definitions, we impose an arbitrary global ordering on the edges of a multigraph. This ordering will be left implicit since the resulting certificate is the same regardless of the ordering.
Given a multigraph $\alpha$, distinguish the $\alpha_{ij}$ copies of the edge $\{i, j\}$ by assigning a label $\ell=1, \hdots, \alpha_{ij}$. The edges of $\alpha$ are ordered lexicographically according to this labelling.

\begin{definition}[Bi-orientation and search forest]
\label{def:bi-orientation}
Let $\alpha$ be a multigraph with edge set $E(\alpha)$. A \textit{bi-orientation} of $\alpha$ assigns each edge $\{v,w\}\in E(\alpha)$ an orientation $v\rightarrow w$, $v\leftarrow w$, or $v\leftrightarrow w$. A singly directed edge $v\rightarrow w$ contributes once to the out-degree of $v$ and once to the in-degree of $w$, and a \textit{bidirected} edge contributes once to the out-degree of each endpoint and nothing to their in-degrees. For a self-loop, the two singly directed orientations are considered distinct.

For a bi-orientation $O$ of $\alpha$, its \textit{search forest} is defined by the following procedure. Start with $\cE_0=R\equiv \{1,\ldots, m\}$ and $F_0=\emptyset$. At each step $t=1, 2, \hdots$, call the vertices in $\cE_{t-1}$ \textit{discovered}. Choose the smallest $v \in \cE_{t-1}$ having an outgoing or bidirected edge to an undiscovered vertex $w$, and let $\cV_t$ be the set of all such undiscovered vertices. For each $w\in\cV_t$, let $e(w)$ be the first (according to the ordering defined above) copy of the edge $\{v,w\}$ oriented $v\rightarrow w$ or $v\leftrightarrow w$. Call $v$ the \textit{parent} of $w$ and $e(w)$ its \textit{parent edge}, and set
\[
    \cE_t=\cE_{t-1}\cup\cV_t,
    \qquad
    F_t=F_{t-1}\cup\{e(w):w\in\cV_t\}\,.
\]
This process terminates when no such vertex $v$ exists. If $T$ is the total number of completed steps, the search forest is $F(O)=F_T$, with vertex set $\cE_T$.
\end{definition}

\begin{definition}[History]
A \textit{history} of shape $\alpha$ is a bi-orientation of $\alpha$ whose search forest contains every vertex in $V(\alpha)$ and has no bidirected parent edges. The set of histories of $\alpha$ is denoted by $\fH(\alpha)$. For a history $h\in\fH(\alpha)$, write 
\begin{equation}\label{eq:rooted:in-out:degree}
    i_h(v)=\deg_h^-(v)-\In_{\{v \not \in R\}}\,, \quad o_h(v)=\deg_h^+(v)\,.
\end{equation}
For a $\sigma$-algebra $\sF$, define
\[
    c_h(\btheta; \sF)
    :=
    \prod_{v\in V(\alpha)}
        A_{o_h(v)}(\theta_v; \sF)\theta_v^{i_h(v)}\,.
\]
\end{definition}

\begin{definition}[History certificate]\label{def:history-cert}
For a $\sigma$-algebra $\sF$, the \textit{history coefficients} are
\begin{equation}
\label{eq:history:coefficient}
    c_\alpha(\btheta; \sF)
    :=
    \sum_{h\in\fH(\alpha)}c_h(\btheta; \sF)\,,
\end{equation}
and the \textit{history certificate} is
\begin{equation}
\label{eq:history:certificate}
    u(\btheta,\bZ; \sF)
    :=
    \sum_{\alpha:|\alpha|\leq D}
        \frac{c_\alpha(\btheta; \sF)h_\alpha(\bZ)}
             {\sqrt{\alpha!}n^{|\alpha|/2}}\,.
\end{equation}
\end{definition}

\begin{remark}
\label{rem:history:support:property}
\label{rem:history:excess:zero}
Note that, if $\alpha$ is non-empty and not root-connected (one of its connected components does not contain a root), then $c_\alpha \equiv 0$. Indeed, in this case, the search forest of any bi-orientation of $\alpha$ cannot span all vertices of $\alpha$, since a root-less component can never be reached.

Also note that, if $\alpha$ is root-connected and has excess $\ex(\alpha)=|E(\alpha)|-|V(\alpha)|+m=0$, then there is only one possible history. This is because, for any history $h \in \fH(\alpha)$, its search forest always has exactly $|F(h)|=|E(\alpha)|-\ex(\alpha)$ edges, all of which must be directed away from the root of its component. If $\ex(\alpha)=0$, then every edge of $\alpha$ belongs to the search forest, and the orientations of all edges in $h$ are forced. This includes the case of the empty graph.
\end{remark}

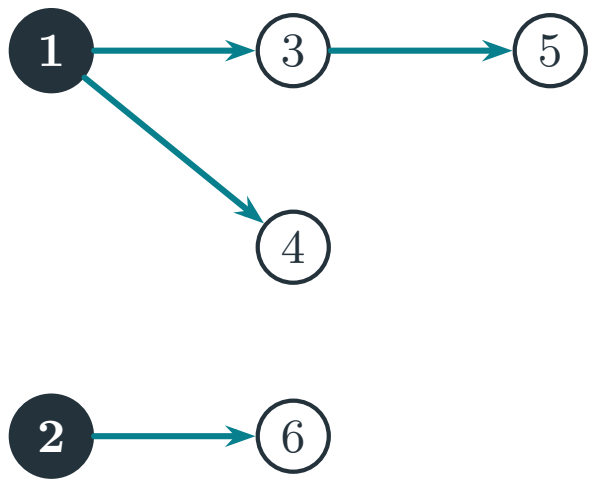
\begin{figure}[t]
    \centering
    \begingroup
\renewcommand{\small}{\fontsize{9}{11}\selectfont}

\begin{tikzpicture}[x=1cm, y=1cm, font=\small, text=ldink]
\node[ldroot]   (r1) at (0.00, 2.55) {1};
\node[ldvertex] (v3) at (1.60, 2.55) {3};
\node[ldvertex] (v4) at (1.60, 1.25) {4};
\node[ldvertex] (v5) at (3.30, 2.55) {5};

\draw[ldparent] (r1) -- (v3);
\draw[ldparent] (r1) -- (v4);
\draw[ldparent] (v3) -- (v5);

\node[ldroot]   (r2) at (0.00, 0.00) {2};
\node[ldvertex] (v6) at (1.60, 0.00) {6};

\draw[ldparent] (r2) -- (v6);

\end{tikzpicture}
\endgroup
    \caption{A forest with zero excess with root set $\{1, 2\}$. Each component contains exactly one root, so every edge is a parent edge directed away from the root.}
    \label{fig:history:example:B}
\end{figure}

\begin{proposition}
\label{prop:history:certifies:dual:constraints}
For any $\sigma$-algebra $\sF$, the coefficients $(c_\alpha(\btheta; \sF))_{\alpha:|\alpha| \leq D}$ in~\eqref{eq:history:coefficient} satisfy the dual constraints~\eqref{eq:dual:constraints:hermite:functional} with $x=\btheta_R$.
\end{proposition}

\begin{proof}
Fix the $\sigma$-algebra $\sF$ and root set $R=\{1, \hdots, m\}$. We show that, for every multigraph $\alpha$, 
\begin{equation}\label{eq:history:dual:target:F}
    \E\left[ \btheta_R(\btheta\btheta^\top)^\alpha \right]
    =
    \sum_{\beta\leq\alpha}\binom{\alpha}{\beta}
    \E\left[
        c_\beta(\btheta;\sF)
        (\btheta\btheta^\top)^{\alpha-\beta}
    \right]\,,
\end{equation}
where we recall the definition of the history coefficient~\eqref{eq:history:coefficient}. In the following, recall that we adopt the convention that $V(\alpha)$ always includes $R$, allowing roots to be isolated vertices. We first prove~\eqref{eq:history:dual:target:F} for simple graphs, then extend the argument to general multigraphs.

\paragraph{Simple graphs.}
Suppose first that $\alpha$ has no parallel edges or self-loops. We expand
\begin{equation}\label{eq:Stein:cumulant:expansion:simple}
    \E[\btheta_R(\btheta\btheta^\top)^\alpha]=\E\left[ \prod_{v \in V(\alpha)} \theta_v^{\deg_\alpha(v)+\In_{\{v \in R\}}} \right]
\end{equation}
by repeated applications of the Stein--cumulant identity (Lemma~\ref{lem:scalar:Stein:cumulant:star}). A vertex $v$ is called \textit{expanded} if we have used the Stein--cumulant identity with $\Theta=\theta_v$. Initially, at step $t=0$, declare every root vertex \textit{discovered} and no vertex expanded. At each step $t \geq 0$, choose the smallest discovered vertex that has not yet been expanded and use the Stein--cumulant identity. Whenever a non-root vertex $i$ is discovered by the procedure below, it is assigned a parent edge $e'(i)$.

For the vertex $v$ currently being expanded, define
\[
    \alpha_v := \begin{cases}
        \{e \in E(\alpha): v \in e\} & \text{if $v \in R$} \\
        \{e \in E(\alpha): v \in e\,,\,e \neq e'(v)\} & \text{if $v \not \in R$}
    \end{cases}
\]
Observe that, before a vertex $v$ is expanded, the exponent on $\theta_v$ in the monomial product $(\btheta\btheta^\top)$ is $\deg_\alpha(v)+\In_{\{v \in R\}}=|\alpha_v|+1$. Conditioning on all other coordinates, independence of the prior and using the identity~\eqref{eq:Stein:cumulant:star:F} with $\Theta=\theta_v$ and $d=|\alpha_v|$ gives
\[
    \E[\theta_v^{|\alpha_v|+1}]=\sum_{r=0}^{|\alpha_v|} \binom{|\alpha_v|}{r} \E[A_r(\theta_v; \sF)\theta_v^{|\alpha_v|-r}]=\sum_{\beta_v \leq \alpha_v} \E[A_{|\beta_v|}(\theta_v; \sF)\theta_v^{|\alpha_v|-|\beta_v|}]\,,
\]
where the second equality follows because there are $\binom{|\alpha_v|}{r}$ subgraphs $\beta_v \leq \alpha_v$ with $|\beta_v|=r$. At each step $t$, pick one choice of $\beta_v$. Draw an arrow $v \rightarrow w$ at every edge $\{v, w\} \in \beta_v$. If an edge is already directed $w \rightarrow v$, this edge becomes bidirected. Every edge in $\beta_v$ deletes exactly one occurrence of $\theta_v$ from the monomial product, and an additional factor of $\theta_v$ is deleted, corresponding to either the factor contributed by $\btheta_R$ if $v$ is a root, or the factor belonging to the parent edge $e'(v)$ if $v$ is a non-root. If an arrow is directed to a previously undiscovered vertex $w$, declare $w$ to be discovered and assign $\{v, w\}$ as its parent edge. Finally, declare $v$ as expanded and continue. This procedure stops if every discovered vertex has been expanded. Upon termination, the result is a sequence of subgraphs $(\beta_v)$ and a bi-orientation $h'$ of
\[
    \beta=\bigcup_{\text{$v$ is expanded}} \beta_v\,.
\]
By construction, it holds that
\[
    \E[\btheta_R(\btheta\btheta^\top)^\alpha]=\sum_{(\beta_v)} \E\left[ \prod_{\text{$v$ is expanded}} A_{|\beta_v|}(\theta_v; \sF)\theta_v^{|\alpha_v|-|\beta_v|} \prod_{\text{$v$ not expanded}} \theta_v^{\deg_\alpha(v)} \right]\,,
\]
where the summation is over all possible tuples of graphs $(\beta_v)$ which arise from this procedure.

\paragraph{Correspondence with histories.}
We show that
\begin{equation}\label{eq:history:expansion:proof:A}
    \sum_{(\beta_v)} \E\left[ \prod_{\text{$v$ expanded}} A_{|\beta_v|}(\theta_v; \sF)\theta_v^{|\alpha_v|-|\beta_v|} \prod_{\substack{\text{$v$ not} \\ \text{expanded}}} \theta_v^{\deg_\alpha(v)} \right]=\sum_{\beta\subseteq\alpha} \sum_{h\in\fH(\beta)} \E\left[ c_h(\btheta;\sF)(\btheta\btheta^\top)^{\alpha-\beta} \right]\,.
\end{equation}
The dual constraints for simple graphs then follow from the definition of the history coefficient in~\eqref{eq:history:coefficient}. To prove this, we will show that there is a one-to-one correspondence between pairs $(\beta, h')$ obtained from the iterative process described above, and pairs $(\beta, h)$ with $h \in \fH(\beta)$.

\begin{enumerate}
    \item[(i)] Consider a sequence of subgraphs $(\beta_v)$ giving rise to $(\beta, h')$ and recall the definition of search forest (Definition~\ref{def:bi-orientation}). We prove $h' \in \fH(\beta)$ by showing that the search forest $F(h')$ spans $V(\beta)$ and contains no bidirected edges. By construction, every non-root vertex of $\beta$ is discovered by some parent edge and eventually expanded. Observe that $F'=\{e'(i): i \in V(\beta) \setminus R\}$, the graph comprised of all such parent edges, is a forest that spans every vertex of $\beta$. Since every parent edge is singly directed, it suffices to show that this forest is precisely the search forest.

    Omit from the expansion procedure every step which does not discover any new vertices, and compare this procedure with the search procedure in Definition~\ref{def:bi-orientation}. Both start with the root vertices $R$ discovered. Suppose inductively that their discovered sets and parent edges agree up to step $t$, and let $v$ be the next vertex to be expanded (if none exist, then both processes terminate, and we are done). Every previously expanded vertex $i$ has no outgoing edges (an edge $i \rightarrow j$ or $i \leftrightarrow j$) to any undiscovered vertex. By definition of the expansion procedure, $v$ is the smallest discovered vertex with an outgoing edge to an undiscovered vertex. Thus, the search procedure also chooses $v$ as the next vertex to search. Since the outgoing edges of $v$ in $h'$ are exactly $\beta_v$, it follows that the two procedures discover exactly the same set of new vertices, and moreover, both procedures assign the same parent edges. This completes the induction, and hence proves that $F'=F(h')$, as required.

    \item[(ii)] Conversely, fix $\beta \leq \alpha$ and consider any history $h \in \fH(\beta)$. We show that a unique sequence of subgraphs $(\beta_v)$ gives rise to $(\beta, h)$. For $v \in V(\beta)$, prescribe
    \begin{equation}\label{eq:history:beta:prescription}
        \beta_v=\Big\{ \{v, w\} \in E(\beta): \text{$v \rightarrow w$ or $v \leftrightarrow w$ in $h$} \Big\}\,.
    \end{equation}
    To show that this sequence can be realized, we must check that $e'(v) \not \in \beta_v$, where $e'(v)$ is the parent edge assigned to $v$ by following the expansion procedure. Suppose, for a contradiction, that $t$ is the first step for which an expanded vertex $v$ has $e'(v) \in \beta_v$. Because the previous steps have $e'(v) \not \in \beta_v$, the same search--comparison argument used in (i) applies to all steps before $t$. In particular, since $v$ is discovered before step $t$, the parent edge $e'(v)$ is the same as the parent edge $e(v)$ in the search forest of $h$. But since $e(v)$ is singly directed into $v$, it cannot possibly be the case that $e(v)=e'(v) \in \beta_v$.

    By construction, the bi-orientation $h'$ obtained from the sequence $(\beta_v)$ is precisely $h$. Finally, since vertices are expanded only once, and no orientations can be removed at a future step, any expansion which produces $h$ must choose each outgoing edge $v \rightarrow w$ or $v \leftrightarrow w$ in $h$ when expanding $v$. This proves uniqueness.
\end{enumerate}
It remains to establish Eq.~\eqref{eq:history:expansion:proof:A}. For $\beta \leq \alpha$ and $h \in \fH(\beta)$, take $(\beta_v)$ to be the unique sequence of subgraphs~\eqref{eq:history:beta:prescription} which realizes $h$. By construction, $|\beta_v|=o_h(v)$. The identities $|\alpha_v|=\deg_\alpha(v)-\In_{\{v \not \in R\}}$ and $\deg_\beta(v)=i_h(v)+\In_{\{v \not \in R\}}+o_h(v)$ imply that $|\alpha_v|-|\beta_v|=\deg_\alpha(v)-\deg_\beta(v)+i_h(v)$. Also, 
\[
    \deg_{\alpha-\beta}(v)=\begin{cases}
        \deg_\alpha(v)-\deg_\beta(v) & \text{if $v \in V(\beta)$} \\
        \deg_\alpha(v)         & \text{if $v \in V(\alpha) \setminus V(\beta)$}
    \end{cases}\,.
\]
Combining the preceding together gives
\[
\begin{aligned}
    \prod_{\text{$v$ expanded}} A_{|\beta_v|}(\theta_v; \sF)\theta_v^{|\alpha_v|-|\beta_v|} \prod_{\substack{\text{$v$ not} \\ \text{expanded}}} \theta_v^{\deg_\alpha(v)} &= \prod_{v \in V(\beta)} A_{o_h(v)}(\theta_v; \sF)\theta_v^{i_h(v)} \prod_{v \in V(\alpha-\beta)} \theta_v^{\deg_{\alpha-\beta}(v)} \\
    &=c_h(\btheta; \sF)(\btheta\btheta^\top)^{\alpha-\beta}\,.
\end{aligned}
\]
This proves~\eqref{eq:history:expansion:proof:A}. In summary, we conclude that
\[
    \E\left[ \btheta_R(\btheta\btheta^\top)^\alpha \right]=\sum_{\beta\subseteq\alpha} \sum_{h\in\fH(\beta)} \E\left[ c_h(\btheta;\sF)(\btheta\btheta^\top)^{\alpha-\beta} \right]=\sum_{\beta\subseteq\alpha} \E\left[ c_\beta(\btheta;\sF)(\btheta\btheta^\top)^{\alpha-\beta} \right]\,.
\]
Since $\binom{\alpha}{\beta}=1$ for simple $\alpha$, this
proves~\eqref{eq:history:dual:target:F} for simple graphs.

\paragraph{Parallel edges and self-loops.}
If $\alpha$ has parallel edges (multiple copies of the same edge), then the same proof applies; this time, there are $\binom{\alpha}{\beta}$ ways to choose parallel edges to produce the same subgraph $\beta$. By symmetry, any choice of the parallel edges leads to the same history coefficient, and~\eqref{eq:history:dual:target:F} holds.

For self-loops, treat the two occurrences of $\theta_v$ in the term $\theta_v^2$ coming from a self-loop as distinct i.e., treat the self-loop as consisting of two distinct half-edges. The Stein--cumulant identity~\eqref{eq:Stein:cumulant:star:F} can be understood as choosing among these distinguished half-edges. Choosing neither half-edge leaves the loop unoriented, and choosing exactly one half-edge gives one of the two singly directed orientations.
With these conventions, the same proof as above applies, and we conclude that~\eqref{eq:history:dual:target:F} holds for general multigraphs.
\end{proof}

\subsection{Norm of the history certificate}
\label{subsec:lower:positive-excess}

Below, we specialize to the family of conditioning $\sigma$-fields $\sF_q=\sigma(q\Theta+\sqrt{q}G)$. Write $c_h(\btheta; q) \equiv c_h(\btheta; \sF_q)$ and similarly $c_\alpha(\btheta; q) \equiv c_\alpha(\btheta; \sF_q)$.

Although the history certificate can be constructed for any $q \geq 0$, a careful choice of $q$ will be needed to bound the norm.
We will work with $q=q_*+\eps$ for small enough $\eps>0$ in the stable case $\Psi''(q_*;\pi)>0$, and $q=q_*$ in the critical case $\Psi''(q_*; \pi) = 0$. Recall that because $\Psi'(q)=\frac{1}{2}(q-\Phi(q))$ (see~\cref{eq:Psi-prime}), we have 
\[
    q_*=\min \{q \geq 0: q=\Phi(q)\,, \Phi'(q) \leq 1\}\,.
\]
Recalling the conventions established earlier, fix $m \geq 1$ and the set of roots $R=\{1, \hdots, m\}$. We say a multigraph $\alpha$ is \emph{root-connected} if each of its connected components contains at least one root. A component is allowed to consist of an isolated root.
We adopt the convention that $V(\alpha)$ includes $R$ by adding roots as isolated vertices.
The set of all root-connected multigraphs is denoted by $\cC$. The dependence on the fixed root set $R$ is left implicit. Let $\cF \subset \cC$ denote the subset of all forests whose component trees contain exactly one root. We use subscripts to denote a degree restriction, so that $\cC_D$ and $\cF_D$ are the set of all root-connected multigraphs and forests with at most $D$ edges, respectively.

By Remark~\ref{rem:history:support:property}, only root-connected multigraphs contribute to the norm of the dual certificate. Thus, our goal is to show that, under the assumptions of \cref{thm:lowdeg:m:roots},
$$ \sum_{\alpha \in \cF_D} \frac{\E_\btheta[c_\alpha(\btheta; q)^2]}{n^{|E(\alpha)|}} \leq q^m + o(1)\,, \qquad \sum_{\al \in \cC_D \setminus \cF_D} \frac{\E_\btheta[c_\al(\btheta; q)^2]}{\al! n^{|E(\al)|}}=o(1)\,.$$

\subsubsection{Norm of forest terms}
\label{subsec:lower:excess-zero}

Rooted forests in $\cF$ have a unique history; see Remark~\ref{rem:history:excess:zero}.
Let $\alpha \in \cF$ have component trees $(\tau_i)_{i=1}^m$, where $\tau_i$ is a tree rooted at vertex $i$. If $i$ is isolated, then $\tau_i=\{i\}$. For a component tree $\tau_i$, abuse notation and let $V(\tau_i)$ denote its ordinary vertex set so that $V(\tau_i)$ does not include any other roots besides $i$. Writing $k_v$ for the number of children of $v$, the history coefficient of $\alpha$ is:
\begin{equation}\label{eq:excess:zero:history}
    c_\alpha(\btheta; q)=\prod_{v \in V(\alpha)} A_{k_v}(\theta_v; q)=\prod_{i=1}^m \prod_{v \in V(\tau_i)} A_{k_v}(\theta_v; q)\,.
\end{equation}
From the factorization in~\eqref{eq:excess:zero:history}, it will be useful to first compute the norm of trees rooted at a fixed vertex. Let $\cT_i$ denote the set of trees rooted at vertex~$i$ and define
$$ w_\tau(q) := n^{-|E(\tau)|} \prod_{v \in V(\tau)} \E[A_{k_v}(\theta_v; q)^2] \,.$$
Here, $V(\tau)$ denotes the ordinary vertex set of $\tau$ and does not add any additional roots.
We have the following key lemma.
\begin{lemma}
\label{lem:tree:summation}
If $\Phi(q) \leq q$, then for all $i \in [n]$, 
\[
    \sum_{\tau \in \calT_i} w_\tau(q) \leq \Phi(q)\,.
\]
\end{lemma}

Since the tree coefficients $c_\tau(\btheta;q_*)$ match the Hermite coefficients of AMP, the sum of squared tree weights should approximately equal the state evolution parameter $q_*$ (see also \cref{lem:idealized-tree-variance}). \cref{lem:tree:summation} shows that $\Phi(q) \leq q$ suffices to upper bound the total weight of the tree terms by $\Phi(q) \leq q$, even when $q$ is not a fixed point of $\Phi$.

\begin{proof}
Without loss of generality, assume $i=1$ and write $\cT \equiv \cT_1$. For a rooted tree $\tau$, let $\he(\tau)$ denote the maximal distance to the root from any vertex of the tree. For an integer $h \geq 0$, define
\begin{equation}
\label{eq:tree:height:induction}
    W_{h} := \sum_{\substack{\tau \in \calT \\ \he(\tau) \leq h}} w_\tau(q)\,.
\end{equation}
Assuming that $\Phi(q) \leq q$, we show by induction on $h$ that $W_{h} \leq \Phi(q)$ for all $h$. This suffices to prove the lemma, by taking $h=n-1$, since every rooted tree on $[n]$ has height at most $n-1$. For the base-case $h=0$, the only tree is an isolated root $\{1\}$, so
\[
    W_{0}=w_{\{1\}}(q)=\E[A_0(\Theta;q)^2] \leq \sum_{k \geq 0} \E[A_k(\Theta; q)^2] \frac{q^k}{k!}=\E[\eta(q\Theta+\sqrt{q}G; q)^2]=\Phi(q)\,,
\]
where the second-to-last step used Parseval's identity, recalling the definition of $(A_k(\Theta; q))_{k \geq 0}$ as the Hermite coefficients of the Bayes denoiser conditioned on $\Theta$.

Now let $h \geq 1$ and assume the induction hypothesis: $W_{h-1} \leq \Phi(q)$. Consider a tree $\tau \in \calT$ rooted at 1 and let $k$ be the number of children of the root.
Let $(\tau_j)_{j=1}^k$ be the subtrees after removing vertex 1 and its incident edges.
The $\tau_j$ are vertex disjoint, do not contain vertex 1, and each has height at most $h-1$. Letting $a_k=\E[A_k(\Theta; q)^2]$, the weights factor across the subtrees as
\[
    w_\tau(q)=\frac{a_k}{n^k} \prod_{j=1}^k \left\{n^{-|E(\tau_j)|} \prod_{v \in V(\tau_j)} a_{k_v}\right\} = \frac{a_k}{n^k} \prod_{j=1}^k w_{\tau_j}(q).
\]
We bound the summation over all possible choices of the subtrees $\tau_1, \hdots, \tau_k$ by
\[
    \sum_{(\tau_j)} \frac{a_k}{n^k} \prod_{j=1}^k w_{\tau_j}(q) \leq \frac{a_k}{n^k} \prod_{j=1}^k \left(\sum_{\tau_j} w_{\tau_j}(q)\right) \leq \frac{a_k}{n^k} \Phi(q)^k \leq \frac{a_k}{n^k}\, q^k\,,
\]
where the first inequality comes from dropping the constraints that the subtrees are disjoint and do not contain vertex 1, and the second inequality uses the induction hypothesis. Finally, there are at most $\binom{n-1}{k} \leq n^k/k!$ total choices of the root children. Taking summation over $k$, the preceding gives
\[
    W_{h} \leq \sum_{k \geq 0} \E[A_k(\Theta; q)^2] \frac{q^k}{k!}=\Phi(q)\,,
\]
where the last step used Parseval's identity.
This completes the induction and concludes the proof.
\end{proof}

\begin{corollary}
\label{prop:excess:zero:terms}
If $\Phi(q) \leq q$, then
\[
    \sum_{\alpha \in \cF_D} \frac{\E[c_\alpha(\btheta; q)^2]}{n^{|E(\alpha)|}} \leq \Phi(q)^m\,.
\]
\end{corollary}

\begin{proof}
Taking summation over all possible choices of the component trees $(\tau_i)_{i=1}^m$, we obtain
\[
    \sum_{\alpha \in \cF_D} \frac{\E[c_\alpha(\btheta; q)^2]}{n^{|E(\alpha)|}} =\sum_{(\tau_i)} \prod_{i=1}^m w_{\tau_i}(q) \leq \prod_{i=1}^m \sum_{\tau \in \calT_i} w_\tau(q) \leq \Phi(q)^m\,,
\]
where the second sum is over tuples of rooted trees $(\tau_i)$ such that (i) $\tau_i$ has root $i$; (ii) $\tau_i$ is vertex-disjoint from $\tau_j$ whenever $i \neq j$; (iii) there are at most $D$ edges in the union of the $\tau_i$. The first inequality comes from dropping condition (ii) and replacing (iii) by each $\tau_i$ individually having at most $D$ edges.
The second inequality uses Lemma~\ref{lem:tree:summation}. This concludes the proof.
\end{proof}

The calculation of the remaining non-forest terms is more involved,
but simple arguments can be used to already conclude the lower bound for $D = O(1)$.
In the case $\mu \neq 0$, this re-proves the lower bound in \cite{montanari2022equivalence} (which uses $m = 1$).
\begin{lemma}
    In the setup of \cref{thm:lowdeg:m:roots}, for all $m \geq 1$ and $D = O(1)$,
    regardless of $\Psi''(q_*;\pi)$,
    \[
        \MMSE_{\le D}(\btheta_R) \geq \E[\btheta^2_R] - q_*^m - o(1)\,.
    \]
\end{lemma}
\begin{proof}
    Let $u=u(\btheta, \bZ; q_*)$ be the history certificate~\eqref{eq:history:certificate} at parameter $q_* \geq 0$. Separate the forest and non-forest component of $u$: 
    \[
        u = \sum_{\tau \in \calF_D} \frac{c_{\tau}(\btheta;q_*)}{n^{|E(\tau)|/2}} \bZ^{\tau} + \sum_{\alpha \in \cC_D \setminus \calF_D} \frac{c_\alpha(\btheta;q_*)h_\alpha(\bZ)}{\sqrt{\alpha!}n^{|E(\alpha)|/2}} =: u_{\calF} + u_{\calC \setminus \calF}\,.
\]
    By~\cref{prop:history:certifies:dual:constraints}, $u$ satisfies the dual constraints.
    It is not essential to use the history certificate;
    one can use the forest component $u_{\calF}$
    then choose any (moment-bounded) coefficients $c_{\alpha}(\btheta)$ on the non-forest terms which satisfy the dual constraints.
    \cref{prop:MMSE:duality} gives the lower bound
    \[
        \MMSE_{\leq D}(\btheta_R) \geq \E[\btheta_R^2]-\E[u^2] = \E[\btheta_R^2] - \E[u_{\calF}^2] - \E[u_{\calC \setminus \calF}^2]\,.
    \]
    By \cref{prop:excess:zero:terms}, it holds
    \[
        \E[u_{\calF}^2] = \sum_{\tau \in \cF_D} \frac{\E[c_\tau(\btheta; q_*)^2]}{n^{|E(\tau)|}} \leq \Phi(q_*)^m = q_*^m\,.
    \]
    For the non-forest terms, let $\sC_D$ be the collection of isomorphism classes of $\cC_D$, and likewise for $\sF_D$.
    For an isomorphism class $\alpha \in \sC_D$, represented by an unlabeled rooted graph, let $H_\alpha(\bZ)$ be the sum of all Hermite polynomials of type $\alpha$,
    \[
        H_\alpha(\bZ) = \sum_{\beta \in \alpha}h_\beta(\bZ)\,.
    \]
    The number of summands is at most $n^{|V(\alpha)|-m}$. Therefore, we have $\E[H_\alpha(\bZ)^2] = O(n^{|V(\alpha)|-m})$ for constant-sized $\alpha$.
    Calculate $\E[u_{\calC \setminus \cF}^2]$ by grouping terms according to their isomorphism class (the value $\E[c_\alpha(\btheta;q_*)^2]$ is the same for isomorphic graphs),
    \[
        \E[u_{\calC \setminus \cF}^2] = \sum_{\alpha \in \sC_D \setminus \sF_D} \frac{\E[c_\alpha(\btheta;q_*)^2] \E[H_\alpha(\bZ)^2]}{\alpha! n^{|E(\alpha)|}} \leq \sum_{\alpha \in \sC_D \setminus \sF_D} \frac{\E[c_\alpha(\btheta;q_*)^2]}{\alpha!} \cdot O(n^{-|E(\alpha)|+|V(\alpha)|-m})\,.
    \]
    The excess $|E(\alpha)| - |V(\alpha)|+ m$ is positive for all non-forests,
    therefore each term $\alpha$ is individually $O(n^{-1})$.
    There are $|\sC_D| = O(1)$ isomorphism classes when $D = O(1)$, thus the summation is $O(n^{-1})$.
    That is, the non-forest terms contribute $O(n^{-1})$ towards $\E[u^2]$, which shows the claim.
    \end{proof}

\subsubsection{Preliminaries for positive-excess terms}
\label{sec:positive-excess}
Let $\ex(\alpha)$ be the rooted excess of a multigraph $\alpha$ with root set $R=\{1, \hdots, m\}$:
\[
    \ex(\alpha)=|E(\alpha)|-|V(\alpha)|+m\,,
\]
where recall $V(\alpha)$ always contains $R$ by convention. The elements of $\cF_D$ have excess zero and all other graphs in $\cC_D \setminus \cF_D$ have positive excess.
To bound the contribution of the positive-excess multigraphs, we partition by their excess.
The main goal for the remainder of this section is to prove the following upper bound.

\begin{proposition}
\label{prop:cyclic:norm:bound}
Assume that $\pi$ is subgaussian. Let $\ell \geq 1$ be an integer.

\begin{enumerate}
\item[(i)] Assume that $\Phi'(q_*)<1$ and let $\eps_0$ be defined in Corollary~\ref{cor:perturbation:of:Psi:Phi}. For every $\eps \in (0, \eps_0)$, there is a constant $C_1=C_1(\eps, \pi, m)>0$ such that, for $q=q_*+\eps$,
\[
    \sum_{\substack{\substack{\alpha \in  \cC_D: \\ \ex(\alpha)=\ell}}} \frac{1}{\alpha!n^{|E(\alpha)|}} \sum_{h_1, h_2 \in \fH(\alpha)} \E[c_{h_1}(\btheta; q)c_{h_2}(\btheta; q)] \leq \left(\frac{C_1D}{n}\right)^\ell\,.
\]
\item[(ii)] In general, there is a constant $C_2=C_2(\pi, m)$ such that, for $q=q_*$, 
\[
\sum_{\substack{\substack{\alpha \in  \cC_D: \\ \ex(\alpha)=\ell}}} \frac{1}{\alpha!n^{|E(\alpha)|}} \sum_{h_1, h_2 \in \fH(\alpha)} \E[c_{h_1}(\btheta; q)c_{h_2}(\btheta; q)] \leq \left(\frac{C_2D^5}{n}\right)^\ell\,.
\]
\end{enumerate}
\end{proposition}

\begin{figure}[t]
    \centering
    \resizebox{0.9\linewidth}{!}{\begingroup
\renewcommand{\small}{\fontsize{11}{13}\selectfont}
\renewcommand{\footnotesize}{\fontsize{10}{12}\selectfont}
\renewcommand{\normalsize}{\fontsize{11}{13}\selectfont}
\renewcommand{\large}{\fontsize{12}{14}\selectfont}
\normalsize

\tikzset{
    ldhead/.append style={font=\large},
    ldroot/.append style={minimum size=6mm}
}

\begin{tikzpicture}[x=1cm, y=1cm, font=\small, text=ldink]
\path[use as bounding box] (0, 0) rectangle (17, 5.25);

\node[ldhead, anchor=center] at (2.300, 4.900)
    {(a) $\alpha \in \cC$};
\node[ldhead, text width=5.5cm, anchor=center] at (8.300, 4.900)
    {(b) $\mathrm{Core}(\alpha)$ removes all non-root leaf vertices};
\node[ldhead, text width=6cm, anchor=center] at (14.525, 4.900)
    {(c) $\mathrm{RedCore}(\alpha)$ contracts non-root degree-2 vertices};

\fill[ldpaleorange, rounded corners=3pt]
    (0.650, 0.150) -- (2.250, 0.150) --
    (1.720, 1.280) -- (1.180, 1.280) -- cycle;
\fill[ldpaleorange, rounded corners=5pt]
    (0.550, 2.950) rectangle (1.050, 3.450);

\node[ldroot]   (A1) at (0.450, 2.200) {1};
\node[ldsmall]  (Av) at (1.450, 2.200) {};
\node[ldsmall]  (Ab) at (2.300, 3.100) {};
\node[ldsmall]  (Ad) at (2.300, 1.300) {};
\node[ldsmall]  (Aw) at (3.150, 2.200) {};
\node[ldroot]   (A2) at (4.150, 2.200) {2};
\node[ldsmall]  (Ag) at (3.150, 3.550) {};
\node[ldsmall]  (Ah) at (4.000, 3.100) {};

\foreach \s/\t in {A1/Av,Av/Ab,Ab/Aw,Av/Ad,Ad/Aw,
    Aw/A2,Aw/Ag,Ag/Ah,Ah/Aw}
    \draw[ldcore] (\s) -- (\t);

\node[ldtreevertex] (At1)  at (1.450, 1.000) {};
\node[ldtreevertex] (At11) at (0.900, 0.400) {};
\node[ldtreevertex] (At12) at (2.000, 0.400) {};
\node[ldtreevertex] (At2)  at (0.800, 3.200) {};

\foreach \s/\t in {Av/At1,At1/At11,At1/At12,Av/At2}
    \draw[ldedge, draw=ldorange] (\s) -- (\t);

\node[ldnote, text=ldorange] at (0.300, 0.650) {$\tau_2$};
\node[ldnote, text=ldorange] at (0.800, 3.700) {$\tau_1$};

\node[ldroot]   (G1) at (6.450, 2.200) {1};
\node[ldsmall]  (Gv) at (7.450, 2.200) {};
\node[ldsmall]  (Gb) at (8.300, 3.100) {};
\node[ldsmall]  (Gd) at (8.300, 1.300) {};
\node[ldsmall]  (Gw) at (9.150, 2.200) {};
\node[ldroot]   (G2) at (10.150, 2.200) {2};
\node[ldsmall]  (Gg) at (9.150, 3.550) {};
\node[ldsmall]  (Gh) at (10.000, 3.100) {};

\foreach \s/\t in {G1/Gv,Gv/Gb,Gb/Gw,Gv/Gd,Gd/Gw,
    Gw/G2,Gw/Gg,Gg/Gh,Gh/Gw}
    \draw[ldcore] (\s) -- (\t);

\node[ldroot]   (R1) at (12.650, 2.200) {1};
\node[ldsmall]  (Rv) at (13.700, 2.200) {};
\node[ldsmall]  (Rw) at (15.350, 2.200) {};
\node[ldroot]   (R2) at (16.400, 2.200) {2};

\draw[ldcore] (R1) -- (Rv);
\draw[ldcore] (Rv) to[bend left=28] (Rw);
\draw[ldcore] (Rv) to[bend right=28] (Rw);
\draw[ldcore] (Rw) -- (R2);
\draw[ldcore] (Rw) edge[loop above, min distance=13mm] (Rw);

\end{tikzpicture}
\endgroup
    }
    \caption{Decomposition of a root-connected multigraph $\alpha$ into its core, reduced  core, and hanging trees. $\alpha$ has two hanging trees, $\tau_1$ and $\tau_2$. Note that $\tau_1,\tau_2$ do not include the edge attaching them to the core.}
    \label{fig:core:graph}
\end{figure}
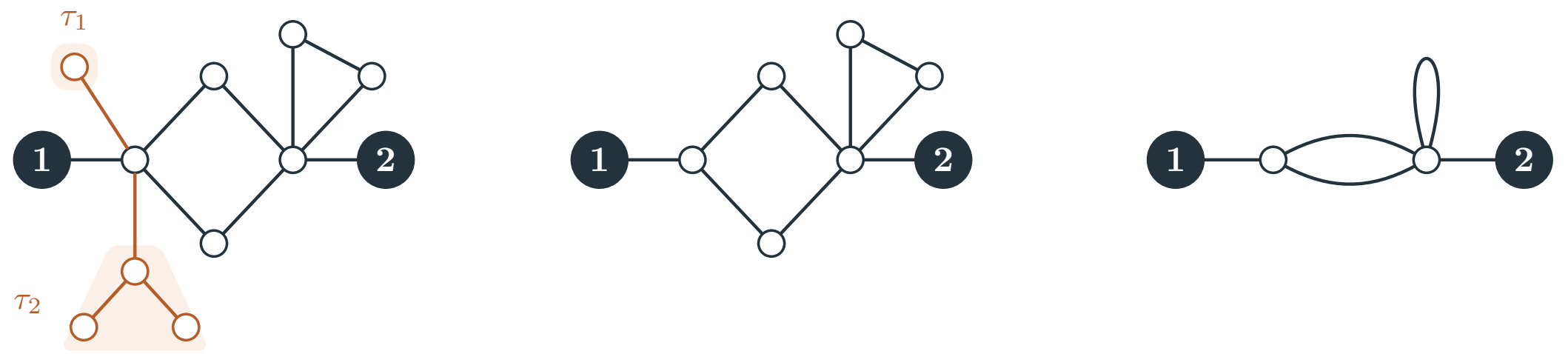

We group the terms $\alpha$ based on their \emph{core}.
There are two types of structures that can make a graph bigger,
but do not increase its excess: hanging trees, and paths of vertices with degree 2.
Removing these structures reduces $\alpha$ to its \emph{core} and \emph{reduced core}, respectively, depicted in Figure~\ref{fig:core:graph}.
\begin{definition}[Core and reduced core]
\label{def:core}\label{def:reduced:core}
Let $\al$ be a root-connected multigraph.
Define $\Core(\alpha)$ by iteratively deleting degree $1$ non-root vertices and their unique incident edge, until no such vertex remains.
For $v \in V(\Core(\alpha))$,
the trees attached to $v$ and removed by this process are called \textit{hanging trees}, written $\cT_v$. An edge connecting a hanging tree to $v$ is not included in the tree.

The \textit{reduced core}, denoted $\RedCore(\alpha)$, is further obtained by replacing (contracting) each maximal path $v_0, \hdots, v_k$ whose internal vertices are distinct, have degree~$2$, and are not roots, with a single edge $\{v_0, v_k\}$. A closed path (the case $v_0=v_k$) is contracted into a self-loop.
Observe that $\Core(\alpha)$ and $\RedCore(\alpha)$ are still root-connected.

A multigraph is called a (reduced) core if it is the (reduced) core of some root-connected multigraph.
\end{definition}

It is clear that $\Core(\alpha)$ and $\RedCore(\alpha)$ have the same excess as $\alpha$ because each step in the construction decreases both the vertex count and edge count by the same amount.
Moreover, there is a natural way to connect histories of $\alpha$ to histories of its core, as described in the lemma below.

Extend the history coefficient notation to a vertex $v$ by letting
    \[
        c_h(\theta_v;q) := \theta_v^{i_h(v)} A_{o_h(v)}(\theta_v;q)\,.
    \]

\begin{lemma}\label{lem:history-correspondence}
    Let $\al \in \cC$ and $\gam = \Core(\al)$.
    There is a one-to-one correspondence between $\fH(\al)$ and $\fH(\gam)$ such that, 
if $h^\al \in \fH(\al)$ corresponds to $h^\gam \in \fH(\gam)$, then
\begin{align}\label{eq:history:core:bijection}
    c_{h^\al}(\btheta;q)
    &=
\left(\prod_{v \in V(\gam)} c_{h^\al}(\theta_v;q)\right)
    \prod_{v\in V(\gamma)}
    \prod_{\tau\in\cT_v(\alpha)}
    \left\{ \prod_{w \in V(\tau)} A_{k_w}(\theta_w; q) \right\}\,.
\end{align}
\end{lemma}
\begin{proof}
Let $\alpha$ be a root-connected multigraph with $\Core(\alpha)=\gamma$. We show that restricting a history of $\alpha$ to the core $\gamma$ induces a bijection from $\fH(\alpha)$ to $\fH(\gamma)$. It will be useful to recall the procedure used to define the search forest (\cref{def:bi-orientation}).

Observe that deleting any edge in $\alpha \setminus \gamma$ from $\alpha$ separates a subtree from the rest of its connected component. Thus, the search forest of any history of $\alpha$ must contain this edge in order to reach vertices of this subtree, and any edge outside $\gamma$ must be singly directed away from the core. In particular, a non-core vertex cannot discover a core vertex. Therefore, if we ignore all non-core discoveries in a search on $\alpha$, then we obtain a search on $\gamma$, with the same set of parent edges. It follows that every history of $\alpha$ restricts to a history of $\gamma$.

Conversely, directing every edge outside of $\gamma$ away from the core extends any history of $\gamma$ to a unique history of $\alpha$, since the search among the core vertices is unchanged, whereas every non-core vertex is discovered through its unique singly directed edge.

Finally, it remains to verify~\eqref{eq:history:core:bijection}. Every hanging tree is directed away from the core, so each of its vertices contributes the factor $A_{k_w}(\theta_w;q)$. Separating the contributions from the core and the hanging trees gives the stated identity.
\end{proof}

The remaining tool we will need in order to prove \cref{prop:cyclic:norm:bound} is an analytic estimate on the coefficients of the state evolution map $\Phi$.
For integers $p,s \geq 0$, define
\[
    M(p,s) = \sum_{r \geq 0} \frac{\Phi(q)^r}{r!}\E[\Theta^{2p} A_{r+s}(\Theta;q)^2]\,.
\]
To motivate the following estimates, consider the power series
\[
    \Gamma_q(y)=\sum_{r \geq 0} \frac{y^r}{r!} \E[A_r(\Theta; q)^2]\,.
\]
Note that plugging in $p=0$ gives $M(0, s)=\Gamma_q^{(s)}(\Phi(q))$, where $\Gamma_q^{(s)}$ is the formal $s$-th derivative of the power series $\Gamma_q$. In the stable case, when $\Phi'(q_*)<1$,~\cref{cor:perturbation:of:Psi:Phi} implies that, for $q=q_*+\eps$, we have $\Gamma_q(\Phi(q)+\delta)<\infty$ for some small $\delta>0$. A standard bound on the derivatives of power series then yields $M(0, s) \leq C^{s+1}s!$. This estimate is a key ingredient in our proof of degree $o(n)$ hardness, and is the main technical reason for working above the AMP fixed point in the stable case; in general, $\Gamma_{q_*}(q)$ may be unbounded for all $q>q_*$.

This is also the point where the subgaussian property of the prior will become relevant, since using Cauchy--Schwarz produces moments of the form $\sqrt{\E[|\Theta|^{4p}]} \leq C^pp!$.

\begin{lemma}\label{lem:analytic-1}
Assume that $\pi$ is subgaussian. The following estimates apply.
    \begin{enumerate}[label=(\roman*)]
        \item If $q=q_*$, there is a constant $C_1=C_1(\pi)>0$ such that $M(p,s) \leq C_1^{p+s+1}((p+s)!)^3$ and moreover, $M(0,1) \leq \frac{1}{2}(1+\Phi'(q_*))$.
        \item If $\Phi'(q_*)<1$ and $q=q_*+\eps$, where $\eps \in (0, \eps_0)$ and $\eps_0$ is defined in~\cref{cor:perturbation:of:Psi:Phi}, then there is a constant $C_2=C_2(\pi, \eps)>0$ such that $M(p,s) \leq C_2^{p+s+1}(p+s)!$ and $M(0,1)\leq\frac12(1+\Phi'(q_*))$.
    \end{enumerate}
\end{lemma}
The proof of \cref{lem:analytic-1} is deferred to Appendix~\ref{subsec:omitted:lower:proofs}.
As shorthand, define
\[
    \rho := \frac{1+\Phi'(q_*)}{2}\,.
\]
The stability assumption will play a crucial role when analyzing the paths contracted in the reduced core: under the strict stability assumption $\Phi'(q_*)<1$ we have $\rho \in (0, 1)$; otherwise, $\rho=1$.

\subsubsection{Norm of positive-excess terms\texorpdfstring{: proof of \cref{prop:cyclic:norm:bound}}{}}

\begin{remark}[Constant $C$]
\label{rem:const-C}
Throughout this section, $C$ will denote a positive constant depending only on $\pi, \eps, m$. This constant may change from line to line. For instance, we will often have multiple factors of the form $C^\ell$ and absorb these into a single $C^\ell$.
\end{remark}

To prove \cref{prop:cyclic:norm:bound}, we start by using Cauchy--Schwarz twice:
\begin{align*}
    \sum_{\substack{\al \in \cC_D: \\ \ex(\alpha) = \ell}} \sum_{h_1,h_2 \in \fH(\alpha)} \frac{\E[c_{h_1}(\btheta; q) c_{h_2}(\btheta; q)]}{\alpha! n^{|E(\al)|}} &\leq \sum_{\substack{\al \in \cC_D: \\ \ex(\alpha) = \ell}} \frac{1}{\alpha! n^{|E(\al)|}}\left(\sum_{h \in \fH(\alpha)} \E[c_{h}(\btheta; q)^2]^{1/2} \right)^2\\
    &\leq \sum_{\substack{\al \in \cC_D: \\ \ex(\alpha) = \ell}} \frac{|\fH(\al)|}{\alpha! n^{|E(\alpha)|}}\sum_{h \in \fH(\alpha)} \E[c_{h}(\btheta; q)^2]
\end{align*}

    \paragraph{Step 1: sum hanging trees.}
    Fix a core $\gamma$ and $h \in \fH(\gam)$.
    For $\al$ such that $\Core(\al) = \gam$, let $h^\al$ be the history extended to $\alpha$. By \cref{lem:history-correspondence} and the definition of $w_\tau(q)$ from \cref{subsec:lower:excess-zero},
    \begin{align}\label{eq:hangingtrees1}
        \sum_{\substack{\al \in \cC_D: \\ \Core(\al) = \gam}} \frac{|\fH(\alpha)|}{\alpha! n^{|E(\alpha)|}}\E[c_{h^\al}(\btheta;q)^2] = \frac{|\fH(\gam)|}{\gam!n^{|E(\gam)|}}\sum_{(\tau_j)} \prod_{v \in V(\gam)} \E[c_{h^\al}(\theta_v;q)^2] \prod_{j} \frac{w_{\tau_j}(q)}{n}
    \end{align}
    where the sum is over any arrangement of hanging trees $\tau_j$ attached to the core $\gamma$.
    The extra factors of $\frac{1}{n}$ come from the edges $\{v, r_j\}$ attaching a hanging tree to the core.
    Let $v \in V(\gam)$. We upper bound the summation over hanging trees at $v$ as follows:
    \begin{align*}
        &\quad \sum_{(\tau_j)} \E[c_{h^\al}(\theta_v;q)^2] \prod_{\tau_j \in \calT_v(\al)} \frac{w_{\tau_j}(q)}{n} \\
        &= \sum_{(\tau_j)} \E\left[A_{o_{h^\al}(v)}(\theta_v;q)^2 \theta_v^{2i_{h^\al}(v)}\right] \prod_{\tau_j \in \calT_v(\al)} \frac{w_{\tau_j}(q)}{n} \\
        \intertext{While the $\tau_j$ must be vertex-disjoint, drop this constraint in order to sum each $\tau_j$ separately:}
        &\leq  \sum_{k = 0}^\infty \E\left[A_{k+o_{h}(v)}(\theta_v;q)^2 \theta_v^{2i_{h}(v)}\right] \frac{1}{k!n^k} \sum_{r_1, \hdots, r_k} \sum_{\tau_1, \dots, \tau_k} \prod_{j=1}^k w_{\tau_j}(q) \\
        \intertext{Here, the sum over $r_j$ is over the possible choices of root label for the hanging trees. The factor $\frac{1}{k!}$ appears above because we have introduced an ordering on the $k$ subtrees. Summing the trees using \cref{lem:tree:summation}, we obtain the quantity $M(p,s)$ defined in \cref{lem:analytic-1}:}
        &\leq \sum_{k = 0}^\infty \E\left[A_{k+o_{h}(v)}(\theta_v;q)^2 \theta_v^{2i_{h}(v)}\right] \frac{\Phi(q)^k}{k!} = M(i_h(v), o_h(v))\,.
    \end{align*}
    Plug this bound back into \cref{eq:hangingtrees1} to get 
    \begin{equation}\label{eq:hangingtrees2}
        \sum_{\substack{\al \in \cC_D: \\ \Core(\al) = \gam}} \frac{|\fH(\al)|}{\alpha! n^{|E(\alpha)|}}\E[c_{h^\al}(\btheta;q)^2] \leq \frac{|\fH(\gam)|}{\gam! n^{|E(\gam)|}} \prod_{v \in V(\gam)} M(i_h(v), o_h(v))\,.
    \end{equation}
 
    \paragraph{Step 2: sum maximal paths of degree~$2$.}
    By \cref{eq:hangingtrees2}, we have
\begin{align*}
        &\quad \sum_{\substack{\al \in \cC_D: \\ \RedCore(\al) = \delta}} \sum_{h \in \fH(\al)} \frac{|\fH(\al)|}{\al! n^{|E(\al)|}}  \E[c_h(\btheta;q)^2]\\
        &= \sum_{\substack{\text{core }\gam: \\ \RedCore(\gam) = \delta \\ |E(\gam)|\leq D}} \sum_{h \in \fH(\gam)}  \sum_{\substack{\al \in \cC_D: \\ \Core(\al) = \gam}} \frac{|\fH(\al)|}{\al! n^{|E(\al)|}}\E[c_{h^\al}(\btheta;q)^2]\\
        &{\le} \sum_{\substack{\text{core }\gam: \\ \RedCore(\gam) = \delta \\ |E(\gam)|\leq D}} \frac{|\fH(\gam)|}{\gam!n^{|E(\gam)|}} \sum_{h \in \fH(\gam)} \prod_{v \in V(\gam)} M(i_h(v), o_h(v))\,. \numberthis\label{eq:path1}
    \end{align*}
    The cores $\gam$ corresponding to $\delta$ are obtained by expanding each edge
    of $\delta$ into a path. Since every internal vertex has positive in-degree, the path is either consistently directed (all edges singly directed in the same direction), directed towards one internal sink, or directed towards one bidirected edge.
    We upper bound the total contribution of these paths as follows.
    \begin{itemize}
    \item All interior vertices of the path contribute a factor $M(i_h(v), o_h(v)) = M(0,1) \leq \rho$, using \cref{lem:analytic-1}, except for at most one factor of $M(1,0) \leq C$ per path.
    
    \item For a path of length $t$, there are at most $n^{t-1}$ ways to label the interior vertices of the path, which cancels with the factor of $n^{-t+1}$ for this path from $n^{-|E(\gam)|}$. Doing this cancellation for all paths gives a factor $n^{-|E(\delta)|}$.

\item Given $h \in \fH(\gamma)$, construct a bi-orientation of $\delta$ as follows. For each contracted path, if all edges of the path are singly directed in the same direction, orient the contracted edge in that direction; otherwise, orient the edge bidirectionally. The resulting bi-orientation $h'$ need not be a history of $\delta$, but this construction preserves the out/in-degree of each vertex in $\delta$.

    A single path of length $t$ can be assigned $\leq 3t$ valid orientations. Thus, letting $t_e$ denote the length of the path $e \in E(\delta)$ in $\gamma$, there are at most $\prod_{e \in E(\delta)} 3t_e$ histories which contract to $h'$.

    \item The preceding also implies $|\fH(\gamma)| \leq \prod_{e \in E(\delta)} 3t_e$.

    \item The factorial $\gam!$ will be counted separately.
    \end{itemize}
    In total, we obtain an upper bound on the summand in \eqref{eq:path1}, assuming that $\delta$ has no multiedges
    \begin{align}
        &\leq \frac{1}{n^{|E(\delta)|}} \sum_{h'} \prod_{v \in V(\delta)} M(i_{h'}(v), o_{h'}(v)) \sum_{(t_e)}\prod_{e \in E(\delta)}C t_e^2 \rho^{t_e} \label{eq:path2}
    \end{align}
    where the sums range over bi-orientations $h'$ of $\delta$ obtained by contracting paths, and over the choices of path lengths $(t_e)$.

Suppose now that $\delta$ has a multiedge. Define $\overline E(\delta):=\{(i,j):1\leq i\leq j\leq n,\ \delta_{ij}>0\}$, the set of distinct endpoint pairs occurring in $\delta$, including self-loops when $i=j$.
For a fixed $e \in \overline{E}(\delta)$, expand each parallel edge into a path and suppose there are $k_t$
    paths of length $t$.
    Then in the second bullet above, the number of labelings overcounts by a factor of $\prod_{t > 1} k_t!$. Mod out this overcount, and note that the additional factor $\gam!$ appearing in the denominator of \eqref{eq:path1} contributes at least $k_1!$ for each multiedge,
    giving an overall denominator $\prod_{t \geq 1}k_t!$. We obtain the following upper bound on \eqref{eq:path1}:
    \begin{align}\label{eq:path3a}
        &\leq \frac{1}{n^{|E(\delta)|}} \sum_{h'} \prod_{v \in V(\delta)} M(i_{h'}(v), o_{h'}(v)) \sum_{(k_{e,t})} \prod_{e \in \overline E(\delta)} \prod_{t\geq 1}  \frac{(Ct^2 \rho^{t})^{k_{e,t}}}{k_{e,t}!}  \intertext{where the sum over $(k_{e,t})$ is over choices for how many paths of length $t$ appear in multiedge $e$.
The total number of edges in $\gam$ is at most $D$,
        but we upper bound the sum over $(k_{e,t})$ by allowing $t$ to range up to $D$ separately for each $e$.
        Now we can upper bound the sum using the multinomial theorem: for integer $r \geq1$ and $x_1, \dots, x_D \in \R$,}
        & \qquad \qquad \sum_{\substack{k_t \geq 0\\\sum_{t=1}^D k_t = r}} \binom{r}{k_1,\dots, k_D} \prod_{t = 1}^D x_t^{k_t} = \left(\sum_{t = 1}^D x_t\right)^r\,. \notag
        \intertext{Apply the multinomial theorem with $x_t = Ct^2\rho^t$ and $r$ equal to the multiplicity of the edge, dividing through by $r!$ to get that \eqref{eq:path3a} is at most:}
        &\leq \frac{1}{\delta!n^{|E(\delta)|}} \sum_{h'} \prod_{v \in V(\delta)} M(i_{h'}(v), o_{h'}(v))\cdot \left(\sum_{t = 1}^D C t^2 \rho^t\right)^{|E(\delta)|}\,. \label{eq:path3}
    \end{align}

\paragraph{Step 3: sum reduced cores.}
The final step is to enumerate the reduced cores. Let $\cR$ denote the set of all reduced cores.
Note that the size of a reduced core is bounded by its excess.
\begin{lemma}
\label{lem:reduced-core-size}
If $\delta$ is a reduced core, then $|V(\delta)| \le 2\ell + m$ and $|E(\delta)| \le 3\ell$.
\end{lemma}
\begin{proof}
Recall $\ell = |E(\delta)| - |V(\delta)| + m$. Since every non-root vertex of $\delta$ has degree $\ge 3$,
\[ 2|E(\delta)|=\sum_{v \in V(\delta)} \deg_\delta(v) \ge 3(|V(\delta)|-m) = 3(|E(\delta)| - \ell), \]
implying $|E(\delta)| \le 3\ell$ and $|V(\delta)| = |E(\delta)| - \ell + m \le 2\ell+m$.
\end{proof}

We enumerate the reduced cores using the ``configuration model''. Consider the following process to construct a multigraph on $[n]$ with roots $R = \{1,2,\dots, m\}$ and $v$ vertices, $e$ edges:
\begin{itemize}
\item First choose $v-m$ non-root vertices among the $n-m$ possibilities. There are $\binom{n-m}{v-m}$ options for this step.
\item Then choose the degrees of all $v$ vertices, such that the degrees sum to $2e$. Using ``stars and bars'', there are $\binom{2e+v-1}{v-1}$ options for this step.
\item Finally, choose a perfect matching on the $2e$ ``half-edges.'' There are $(2e-1)!!$ options for this step. (For counting purposes, we treat these $2e$ half-edges as distinguishable, even if they belong to the same vertex.)
\end{itemize}
Every reduced core has at least one configuration that generates it (but some configurations do not generate a valid reduced core, e.g., if a non-root vertex has degree~$2$).

\begin{lemma}\label{lem:enumeration-overcount}
    Every reduced core $\delta \in \cR$ is enumerated exactly $\frac{\prod_{v \in V(\delta)} \deg_{\delta}(v)!}{\delta!2^{s}}$ times, where $s = \sum_{i \in V(\delta)}\delta_{ii}$
    is the number of self-loops.
\end{lemma}
\begin{proof}
    Given $\delta$, the choice of vertices and degrees is fixed, and the overcounting arises from the matching on half-edges. Suppose a particular vertex $v \in V(\delta)$ has $s$ self-loops along with $p$ neighbors (not including itself) connected to $v$ by edges of multiplicity $k_1,\ldots,k_p \ge 1$ respectively. To produce a matching on half-edges that generates $\delta$, we first have a choice of how $v$ allocates its half-edges (of which there are $\deg_\delta(v) = 2s + \sum_i k_i$) among its neighbors (including itself, for self-loops). The number of ways to make this allocation is the multinomial coefficient
    \[ \frac{\deg_\delta(v)!}{(2s)! \, k_1! \cdots k_p!} \, . \]
    The analogous choice is made at every vertex. Second, for each multiedge in $\delta$ with multiplicity $r$, there are $r!$ ways to match the associated half-edges; for a self-loop of multiplicity $s$, the number of choices is instead $(2s-1)!!$. In total, the overcounting factor for $\delta$ contains:
    \begin{itemize}
    \item a factor of $\deg_\delta(v)!$ for each vertex,
    \item a factor of $\frac{r!}{(r!)^2} = \frac{1}{r!}$ for each edge of multiplicity $r$ (not including self-loops), and
    \item a factor of $\frac{(2s-1)!!}{(2s)!} = \frac{1}{s! \, 2^s}$
    for each self-loop of multiplicity $s$.
    \end{itemize}
    This concludes the proof.
\end{proof}

\begin{corollary}\label{lem:config:count}
For all integers $\ell\geq1$ and $m\leq v\leq2\ell+m$,
\[
    \sum_{\substack{\textnormal{reduced core }\delta:\\
        \ex(\delta)=\ell,\ |V(\delta)|=v}}
    \frac{1}{\delta!}
    \prod_{i\in V(\delta)}\deg_\delta(i)!
    \leq
    C^\ell n^{v-m}
    \frac{(\ell+v-m)!}{(v-m)!}
    \leq n^{v-m}(C\ell)^\ell\,.
\]
\end{corollary}

\begin{proof}
Let $e = \ell + v -m$.
By \cref{lem:enumeration-overcount}, the left-hand summation is at most $2^\ell N(v,e)$
where $N(v,e)$ is the total number of executions of the configuration model.
The factor $2^\ell$ accounts for the term $2^s$ which appears in \cref{lem:enumeration-overcount} but not here,
but which satisfies $2^s \leq 2^\ell$ since every self-loop is an excess edge.

We have $N(v,e) = \binom{n-m}{v-m} \binom{2e+v-1}{v-1} (2e-1)!!$. Then, using $v \leq C\ell$ and $e \leq C\ell$ by \cref{lem:reduced-core-size} (recall we allow $C$ to change as per Remark~\ref{rem:const-C}),
\[
    N(v,e) \le \frac{(n-m)^{v-m}}{(v-m)!} \, 2^{2e+v-1} 2^{e}e! \le C^\ell n^{v-m}\frac{e!}{(v-m)!} \,.
\]
This proves the first inequality. For the second inequality, we have
\[
    \frac{e!}{(v - m)!} = \ell! \binom{e}{\ell} \leq \ell^\ell 2^{e} \leq (C\ell)^\ell\,. \qedhere
\]
\end{proof}

Now we perform the overall summation.
Partition the sum over all $\al$ by their reduced core and use \eqref{eq:path3}:
    \begin{align*}
        &\quad\sum_{\substack{\al \in \cC_D: \\ \ex(\alpha) = \ell}} \sum_{h \in \fH(\alpha)} \frac{|\fH(\al)|}{\alpha! n^{|E(\alpha)|}}\E[c_{h}(\btheta; q)^2]\\
        &= \sum_{\substack{\text{reduced core }\delta: \\ \ex(\delta) = \ell}} \sum_{\substack{\al \in \cC_D: \\ \RedCore(\al) = \delta}} \sum_{h \in \fH(\al)} \frac{|\fH(\al)|}{\alpha! n^{|E(\al)|}} \E[c_h(\btheta;q)^2]\\
        &\leq \sum_{\substack{\text{reduced core }\delta: \\ \ex(\delta) = \ell}} \frac{1}{\delta! n^{|E(\delta)|}} \sum_{h'} \prod_{v \in V(\delta)} M(i_{h'}(v), o_{h'}(v)) \left(\sum_{t = 1}^D Ct^2 \rho^t\right)^{|E(\delta)|} \,.\numberthis\label{eq:cores1}\\
        \intertext{We specialize the proof for the two cases from \cref{prop:cyclic:norm:bound}. In the stable case (i), then $\rho < 1$ implies that $\sum_{t \geq 1} t^2 \rho^t = C$ is a finite constant. Additionally, $M(i_{h'}(v), o_{h'}(v)) \leq C^{1+\deg_{\delta}(v)}\deg_{\delta}(v)!$ by \cref{lem:analytic-1}, and the number of $h'$ is at most $3^{|E(\delta)|}$, the number of bi-orientations of $\delta$.
        Combining these bounds,}
        &\leq \sum_{\substack{\text{reduced core }\delta: \\ \ex(\delta) = \ell}} \frac{C^{|V(\delta)| + |E(\delta)|}}{\delta!n^{|E(\delta)|}} \prod_{v \in V(\delta)} \deg_{\delta}(v)!\\
        &\overset{\text{Lem.~\ref{lem:reduced-core-size}}}{\leq} C^\ell \sum_{\substack{\text{reduced core }\delta: \\ \ex(\delta) = \ell}} \frac{1}{\delta!n^{|E(\delta)|}} \prod_{v \in V(\delta)} \deg_{\delta}(v)!\\
        &\overset{\text{Cor.~\ref{lem:config:count}}}{\leq} C^\ell \sum_{v=m}^{2\ell+m} \frac{(C\ell)^\ell}{n^{\ell}} \\
    &= (2\ell + 1)\left(\frac{C\ell}{n}\right)^\ell\,.
    \end{align*}
    Use the upper bound $2\ell+1 \leq C^\ell$ then $\ell \le D$ to get the upper bound $(CD/n)^\ell$.
    This completes the proof of case (i).

\medskip

For case (ii) of \cref{prop:cyclic:norm:bound}, we have $\rho \leq 1$ and $M(i_{h'}(v), o_{h'}(v)) \leq C^{1+\deg_{\delta}(v)} (\deg_{\delta}(v)!)^3$ by \cref{lem:analytic-1}.
The argument of case (i) proves the upper bound $(CD^{k}/n)^\ell$ for an integer $k \leq 20$
after using $\sum_{t = 1}^D t^2 \leq D^3$ in \eqref{eq:path3}
and bounding the factorial terms by $\prod_{v \in V(\delta)}\deg_{\delta}(v)! \leq (CD)^{2\ell}$ (to be proven in \cref{lem:factorial-bound}). We optimize the value of $k$ starting from \eqref{eq:hangingtrees2}. We have:
\[
    \frac{|\fH(\gam)|}{\gam! n^{|E(\gam)|}} \prod_{v \in V(\gam)} M(i_h(v), o_h(v)) \leq \frac{C^{|V(\delta)| + |E(\delta)|}|\fH(\gam)|}{\gam! n^{|E(\gam)|}} \prod_{v \in V(\delta)} (\deg_{\delta}(v)!)^3\,. 
\]
This bound is uniform over the history $h$. Using \cref{lem:reduced-core-size} the factors of $C$ are at most $C^\ell$. The total value of a single reduced core $\delta$ is then
\begin{align}
    &\quad \sum_{\substack{\al \in \cC_D: \\ \RedCore(\al) = \delta}} \sum_{h \in \fH(\al)} \frac{|\fH(\al)|}{\al! n^{|E(\al)|}}  \E[c_h(\btheta;q)^2] {\le} \sum_{\substack{\text{core }\gam: \\ \RedCore(\gam) = \delta \\ |E(\gam)|\leq D}} \frac{C^{\ell}|\fH(\gam)|^2}{\gam!n^{|E(\gam)|}} \prod_{v \in V(\delta)} (\deg_{\delta}(v)!)^3\,.\label{eq:case2:1}
\end{align}
Use the uniform upper bound,
\begin{align}\label{eq:hist-count}
    |\fH(\gam)| \leq 3^{\ell} \binom{|E(\gam)|}{\ell} \leq \frac{(3D)^{\ell}}{\ell!}\,.
\end{align}
To prove the first inequality in \eqref{eq:hist-count}, every history determines an oriented spanning forest of $\gam$;
we choose the remaining $\ell$ edges then choose their orientations. This proves \eqref{eq:hist-count}. Then \eqref{eq:case2:1} is at most
\begin{align}
    &\leq \frac{(CD)^{2\ell}}{(\ell!)^2} \prod_{v \in V(\delta)} (\deg_{\delta}(v)!)^3 \sum_{\substack{\text{core }\gam: \\ \RedCore(\gam) = \delta \\ |E(\gam)|\leq D}} \frac{1}{\gam!n^{|E(\gam)|}} \\
    &\leq \frac{(CD)^{2\ell}}{(\ell!)^2n^{|E(\delta)|}} \prod_{v \in V(\delta)} (\deg_{\delta}(v)!)^3 \sum_{(k_{e,t})} \prod_{e \in \overline E(\delta)} \prod_{t \geq 1}\frac{1}{k_{e,t}!}.\label{eq:case2:2}
\end{align}
The second inequality follows the same reasoning as \eqref{eq:path3a}.
The sum over $(k_{e,t})$ is over the number of paths of each length $t$ in each multiedge $e$, such that the total length of all paths is at most $D$.
In the current case, the sum over $(k_{e,t})$ can be evaluated exactly. We claim:
\[
    \delta!\sum_{(k_{e,t})} \prod_{e \in \overline E(\delta)} \prod_{t \geq 1}\frac{1}{k_{e,t}!}=\sum_{(k_{e,t})}\prod_{\substack{\text{multiedge }e \in \overline E(\delta)\\\text{with multiplicity }r_e}} \binom{r_e}{k_{e,1},k_{e,2},\dots}=\binom{D}{|E(\delta)|}\,.
\]
The first equality is the definition of the multinomial coefficient. The second equality
counts the number of ways to expand $E(\delta)$ into paths with total length at most $D$:
either $\binom{D}{|E(\delta)|}$ by ``stars and bars'', or by choosing the number of paths of each length in each multiedge.
Dividing through by $\delta!$ we obtain that \eqref{eq:case2:2} is equal to
\begin{align}\label{eq:case2:3}
    = \frac{(CD)^{2\ell}}{(\ell!)^2 \delta!n^{|E(\delta)|}} \binom{D}{|E(\delta)|} \prod_{v \in V(\delta)} (\deg_{\delta}(v)!)^3 \leq \frac{(CD)^{2\ell + |E(\delta)|}}{(\ell!)^2 |E(\delta)|! \delta!n^{|E(\delta)|}} \prod_{v \in V(\delta)} (\deg_{\delta}(v)!)^3.
\end{align}
Next, we bound the factorial terms as follows.
\begin{lemma}\label{lem:factorial-bound}
$\prod_{v \in V(\delta)} \deg_{\delta}(v)! \leq C^\ell (2\ell)!$
\end{lemma}
\begin{proof}
    The total number of terms among the factorials is $2|E(\delta)| = 2|V(\delta)| + 2\ell - 2m$. 
    Every non-root $v \in V(\delta)$ has degree at least 3, so pull out its top two factors:
    \begin{align*}
        \prod_{v \in V(\delta) \setminus R} \deg_{\delta}(v)! &= \prod_{v \in V(\delta) \setminus R} \deg_{\delta}(v) \cdot (\deg_{\delta}(v) - 1) \cdot (\deg_{\delta}(v) - 2)!\\
        &\leq \prod_{v \in V(\delta) \setminus R} C^{\deg_{\delta}(v) - 2} \cdot (\deg_{\delta}(v) - 2)!\,.
    \end{align*}
    The total exponent on $C$ is at most $2\ell$.
    There are now only $2\ell$ factorial terms remaining. Use $\prod_{k} a_k! \leq (\sum_k a_k)!$ for non-negative integers $a_k$ to upper bound them by $(2\ell)!$.
\end{proof}
Apply \cref{lem:factorial-bound} to 2 of the 3 factorial terms in \eqref{eq:case2:3}. Then \eqref{eq:case2:3} is at most:
\begin{align}\label{eq:case2:4}
    &\leq \frac{((2\ell)!)^2(CD)^{2\ell + |E(\delta)|}}{(\ell!)^2 |E(\delta)|! \delta!n^{|E(\delta)|}} \prod_{v \in V(\delta)} \deg_{\delta}(v)! \,.\end{align}
Perform the final summation over all $\alpha$:
\begin{align*}
    &\quad\sum_{\substack{\al \in \cC_D: \\ \ex(\alpha) = \ell}} \sum_{h \in \fH(\alpha)} \frac{|\fH(\al)|}{\alpha! n^{|E(\alpha)|}}\E[c_{h}(\btheta; q)^2]\\
    &\overset{\text{\eqref{eq:case2:4}}}{\leq} \frac{((2\ell)!)^2 C^\ell}{(\ell!)^2} \sum_{\substack{\text{reduced core }\delta: \\ \ex(\delta) = \ell}} \frac{D^{2\ell + |E(\delta)|}}{|E(\delta)|!\delta!n^{|E(\delta)|}} \prod_{v \in V(\delta)} \deg_{\delta}(v)! \\
    \intertext{Define the notation $e = \ell+v-m$.}
    &= \frac{((2\ell)!)^2 C^\ell}{(\ell!)^2} \sum_{v = m}^{2\ell + m} \frac{D^{2\ell + e}}{e!n^{e}} \sum_{\substack{\text{reduced core }\delta:\\\ex(\delta) = \ell,\, |V(\delta)| = v}} \frac{1}{\delta!}\prod_{v \in V(\delta)} \deg_{\delta}(v)!\\
    &\overset{\text{Cor.~\ref{lem:config:count}}}{\leq} \frac{((2\ell)!)^2 C^\ell}{(\ell!)^2} \sum_{v = m}^{2\ell + m} \frac{D^{2\ell + e}}{(v-m)! n^{\ell}}
    \end{align*}
Bound the factorial terms as $\frac{(2\ell)!}{(\ell!)^2} = \binom{2\ell}{\ell} \leq 4^\ell$ and
\[
    \frac{(2\ell)!}{(v-m)!} \leq (2\ell)^{2\ell - v + m} \leq C^\ell D^{2\ell - v + m}
\]
The total exponent on $D$ is $5\ell$ with no more dependence on $v$. Continuing from above,
\begin{align*}
    &\leq \frac{C^\ell D^{5\ell}}{n^\ell} \sum_{v = m}^{2\ell + m} 1 \le \frac{C^\ell D^{5\ell}}{n^\ell}(2\ell+1) \leq \left(\frac{CD^5}{n}\right)^\ell\,,
\end{align*}
using $2\ell+1 \leq C^\ell$ in the last step. This completes the proof for case (ii).

\subsection{Proof of Theorem~\ref{thm:lowdeg:m:roots}}
\label{sec:lowdeg:proof}

\begin{proof}[Proof of Theorem~\ref{thm:lowdeg:m:roots}]
As shown in Section~\ref{sec:se}, the stability condition $\Psi''(q_*; \pi)>0$ is equivalent to $\Phi'(q_*)<1$. Let $\eps_0$ be defined in~\cref{cor:perturbation:of:Psi:Phi}. For $\eps \in (0, \eps_0)$, take $q=q_*+\eps$. Let $u=u(\cdot; q)$ be the history certificate~\eqref{eq:history:certificate} with conditioning $\sigma$-field $\sF_q=\sigma(q\Theta+\sqrt{q}G)$. By~\cref{prop:history:certifies:dual:constraints}, $u$ satisfies the dual constraints and therefore~\cref{prop:MMSE:duality} gives the lower bound, 
\[
    \MMSE_{\leq D_n}(\btheta_R) \geq \E[\btheta_R^2]-\E[u^2]\,.
\]
For $1\leq\ell\leq D_n$, define the positive-excess contribution
\[
    N_\ell:=\sum_{\substack{\alpha\in \cC_{D_n}\\\ex(\alpha)=\ell}}
    \frac{\E[c_\alpha(\btheta;q)^2]}{\alpha!n^{|\alpha|}}\,.
\]
Bounding the excess-zero contribution to the norm using~\cref{prop:excess:zero:terms}, 
\[
    \E[u^2]=\sum_{\ell=0}^{D_n} \sum_{\substack{\alpha:|\alpha| \leq D_n \\ \ex(\alpha)=\ell}} \frac{\E[c_\alpha(\btheta; \sF_q)^2]}{\alpha!n^{|\alpha|}} \leq q^m+\sum_{\ell=1}^{D_n}N_\ell\,.
\]
By~\cref{prop:cyclic:norm:bound}-(i), there is a constant $C=C(\eps, \pi, m)$ such that, for all sufficiently large $n$, 
\[
    \sum_{\ell=1}^{D_n} N_\ell \leq \sum_{\ell \geq 1} \left( \frac{CD_n}{n} \right)^\ell \leq 2CD_n/n \rightarrow 0\, \quad \text{as $n \rightarrow \infty$}\,,
\]
provided $D_n=o(n)$. Together, the preceding implies that 
\[
    \liminf_{n\rightarrow\infty} \MMSE_{\leq D_n}(\btheta_R) \geq \E[\btheta_R^2]-(q_*+\eps)^m\,.
\]
Letting $\eps \downarrow 0$ concludes the proof for the stable case. In the critical case $\Psi''(q_*; \pi)=0$, then we have $\Phi'(q_*)=1$ by~\eqref{eq:Psi-prime}.
Taking $q=q_*$, the lower bound for degree $D_n=o(n^{1/5})$ follows the same argument, using~\cref{prop:cyclic:norm:bound}-(ii) to bound the positive-excess contribution instead.
\end{proof}

\section{Upper bound}
\label{sec:upper-bound}

We prove the low-degree upper bound in \Cref{thm:lowdeg}-\ref{item:upper}. In fact, the upper bound holds under the following weaker assumption than subgaussian, which we assume throughout this section:
\begin{equation}
 \E_{\Theta \sim \pi}[e^{\alpha|\Theta|}]<\infty\qquad\text{for some }\alpha>0.
 \label{eq:upper-moment-assumption}
\end{equation}
Throughout this section, $\pi$ is fixed as $n\to\infty$, and we write $\mu:=\E\Theta$ and $\lambda:=\E[\Theta^2]$.
In this section it will be convenient to use the rescaled observation
\begin{equation}\label{eq:Y-bar}
    \overline{\bY}:=\frac{\bY}{\sqrt n}
    =\frac{\btheta\btheta^\top}{n}+\overline{\bZ},
    \qquad \overline{\bZ}:=\frac{\bZ}{\sqrt n},
\end{equation}
where $\bY,\bZ$ are as in \cref{def:wig}. Thus the upper triangular
entries of $\overline{\bZ}$ are independent $\mathcal N(0,1/n)$
variables.

We construct an explicit polynomial approximation to AMP. There are three cases that need to be handled differently, depending on the prior $\pi$ (see \cref{fig:upper-bound-constructions}). \Cref{sec:upper-preliminaries} develops some common ingredients: polynomial approximation of the scalar denoisers and the corresponding tree polynomial approximation for AMP.
\Cref{sec:upper-symmetric} gives the proof in the case where $\pi$ is symmetric around zero. \Cref{sec:nonzero-mean,sec:zero-asymmetric} explain how to adapt the proof to the other two cases: nonzero mean and non-symmetric zero mean (the nonzero mean case of \Cref{thm:lowdeg}-\ref{item:upper} was already proven by \cite{montanari2022equivalence}, but we include it for expository purposes). Finally, \Cref{sec:upper-main-theorem} combines these cases to prove \Cref{thm:lowdeg}-\ref{item:upper}.

\subsection{Preliminaries}
\label{sec:upper-preliminaries}

\subsubsection{Polynomial approximation of the Bayes denoiser}

The behavior of AMP is governed by the state evolution map in \eqref{eq:phi}. We repeat its definition here, using the notation of \cref{eq:denoiser}. For \(a>0\), define
\begin{equation}
 X_a\coloneqq a\Theta+\sqrt a G,\qquad
 \eta(x;a)\coloneqq\E[\Theta\mid X_a=x],\qquad
 \Phi(a)\coloneqq\E[\Theta\eta(X_a;a)]=\E[\eta(X_a;a)^2],
 \label{eq:scalar-channel}
\end{equation}
where the Gaussian noise \(G\sim\mathcal N(0,1)\) is independent
of \(\Theta \sim \pi\).
For a fixed initialization $q_1> 0$, define the sequence of iterates
\begin{equation}
 q_{t+1}\coloneqq\Phi(q_t).
 \label{eq:nonzero-q}
\end{equation}

Towards approximating AMP by a polynomial, first we approximate $\eta(\cdot;a)$ by a constant-degree polynomial.
For \(a>0\) and an integer \(d\ge1\), let
\begin{equation}
 \hat\eta(\cdot;a,d)\coloneqq\operatorname*{arg\,min}_{\deg g\le d}
                 \E[(g(X_a)-\Theta)^2].
 \label{eq:least-squares}
\end{equation}
Note that without the degree constraint, the minimizer is the posterior mean $\eta(\cdot;a)$. Thus, $\hat{\eta}(\cdot;a,d)$ should be thought of as a degree-$d$ polynomial approximation to $\eta(\cdot;a)$.
Since \(X_a\) has a density, the objective is strictly convex in the polynomial coefficients, so the minimizer is unique.

This choice of $\hat{\eta}(\cdot;a,d)$ has a convenient Nishimori-type identity.
The proof of \cref{lem:calibration} gives the usual Nishimori identity $\E[\Theta\eta(X_a;a)] = \E[\eta(X_a;a)^2]$ when the degree constraint on $\hat \eta$ is dropped.

\begin{lemma}[Nishimori identity for $\hat\eta$]
\label{lem:calibration}
\[ \E[\Theta \hat\eta(X_a;a,d)]
 =\E[\hat\eta(X_a;a,d)^2].
 \]
\end{lemma}

\begin{proof}
The directional derivative of \eqref{eq:least-squares} in the direction of a polynomial $p$ of degree at most $d$ is
\[
 2\E[(g(X_a)-\Theta)p(X_a)].
\]
When $g = \hat \eta(\cdot;a,d)$ is the minimizer, this is 0 for all $p$. Choosing \(p=\hat\eta(\cdot;a,d)\) gives the equality.
\end{proof}

Analogous to $q_t$ for $\eta$, we can define a sequence $a_t$
for the state evolution for AMP using the denoisers $\hat \eta$.
Fix an initialization $a_1 > 0$ and a degree bound \(d\ge1\),
and define $a_t$ by the recurrence
\begin{equation}
  \hat\eta_t(x) \coloneqq \hat\eta(x;a_t,d),\qquad a_{t+1}\coloneqq\E[\Theta \hat \eta_t(X_{a_t})] = \E[\hat\eta_t(X_{a_t})^2].
 \label{eq:zero-a}
\end{equation}
As expected, the sequences $a_t$ and $q_t$ have similar behavior for large $d$.

\begin{lemma}[Approximation of the state evolution]
\label{lem:approx-scalar}
Assume~\eqref{eq:upper-moment-assumption}. Let $q_1 = a_1 > 0$ and \(d\ge1\). For each fixed $t$, the iterates $a_t$ defined by~\eqref{eq:zero-a} and $q_t$ defined by \eqref{eq:nonzero-q} satisfy
\begin{equation}
 a_t\longrightarrow q_t \quad\mathrm{as}\quad d\to\infty.
 \label{eq:zero-approximation}
\end{equation}
\end{lemma}

\noindent The proof is deferred to \Cref{app:scalar-approximation}.

\subsubsection{Hermite expansion}

Recall the probabilists' Hermite polynomials from \cref{subsec:notation}. As derived in \cref{sec:Y-basis}, the polynomial expansion of AMP involves the Hermite expansion of the denoisers. For the sequence of denoisers $\hat \eta_t$ defined by \eqref{eq:zero-a},
write the Hermite expansion,
\begin{equation}
 \hat\eta_t(x)=\hat\eta(x;a_t,d)=\sum_{k=0}^d b_{t,k}\He_k(x;a_t).
 \label{eq:weights}
\end{equation}
We will use the Hermite coefficients $b_{t,k}$ to define our polynomial approximation.
We also need the identities
\begin{align}
 \E_G[\He_k(a\theta+\sqrt aG;a)]&=a^k\theta^k,
 \label{eq:hermite-first}\\
 \E_G[\He_k(a\theta+\sqrt aG;a)\He_\ell(a\theta+\sqrt aG;a)]
 &=\sum_{j=0}^{k\wedge \ell}\binom kj\binom \ell j j!
       a^{k+\ell-j}\theta^{k+\ell-2j}.
 \label{eq:hermite-second}
\end{align}
These follow, for instance, by the Hermite translation formula
\[
 \He_k(a\theta+\sqrt aG;a)
 =\sum_{j=0}^k\binom kj(a\theta)^{k-j}a^{j/2}\He_j(G)
\]
and the orthogonality relation \(\E[\He_j(G)\He_\ell(G)]
=j!\,\In_{\{j=\ell\}}\).

\subsubsection{Weighted tree polynomials}

It will be convenient to work with \emph{rooted ordered trees} where each vertex has a specified order on the children (also called \emph{plane trees}). This will eliminate some factorial factors that would occur without the ordering. We will define $\mathscr T_t$, the set of such trees (up to isomorphism) with depth $\le t$ and branching factor $\le d$. Formally, let \(\mathscr T_0\coloneqq\{\bullet\}\) where $\bullet$ denotes the empty tree consisting of a single root vertex. Recursively, a tree in \(\mathscr T_{t+1}\) consists of a root with an ordered list of
\(k\in\{0,1,\ldots,d\}\) children, each of which is a tree in
\(\mathscr T_t\). For $\tau \in \mathscr T_{t+1}$ we write \(\tau=(\tau_1,\ldots,\tau_k)\),
where \(\tau_j \in \mathscr T_t\) is the subtree rooted at the \(j\)th child of the root.
Write $k_u$ for the number of children of vertex $u$.

A tree \(\tau \in \mathscr T_t\) has depth at most \(t\). Its boundary $\partial_t\tau$
consists of the vertices at depth exactly \(t\); write \(\ell_t(\tau) = |\partial_t\tau|\)
for their number. {We write $\ell(\tau)$ and $\partial \tau$ for notational simplicity, although a leaf at depth less than $t$ is not included in $\partial \tau$.}
Using the weights $b_{t,k}$ from \eqref{eq:weights}, define the weight of a tree recursively {as suggested by the derivation in \cref{sec:Y-basis}:}
\begin{equation}
\begin{gathered}
 B_{0,\bullet}\coloneqq1,\qquad
 B_{t,\tau}\coloneqq b_{t,k}\prod_{j=1}^kB_{t-1,\tau_j}
 \\
 \text{when } \tau \in \mathscr T_t \text{ with } \tau=(\tau_1,\ldots,\tau_k).
\end{gathered}
 \label{eq:tree-weights}
\end{equation}
Empty products equal one by convention. In other words, for \(\tau\in\mathscr T_t\), each vertex at depth \(h<t\) with \(k\)
children contributes a factor of \({b_{t-h, \, k}}\) while a boundary leaf has no factor in $B_{t,\tau}$.

For a fixed depth $T$, the finite family \(\mathscr T_T\) only depends on \(T,d\).
Every tree in it has at most \(d^T\) boundary vertices and at most
\begin{equation}
 D_T\coloneqq d+d^2+\cdots+d^T
 \label{eq:degree}
\end{equation}
edges, with \(D_0\coloneqq0\).

{Next, we give tree formulas for the state evolution parameters, one for ``$m_t$'' and the second for ``$q_t$'', although due to the Nishimori identity \cref{lem:calibration} both formulas are equal to $a_t$.
See \cref{lem:idealized-tree-variance} for the analogous formula in the ``$Z$ basis''.}
We make the following combinatorial definitions involving overlapping pairs of trees $\tau, \tau' \in \mathscr T_t$, whose vertex sets are denoted $V(\tau),V(\tau')$. A partial matching between $V(\tau)$ and $V(\tau')$ is a subset
$\chi\subseteq V(\tau)\times V(\tau')$ such that each vertex occurs in at most one pair, and the roots are paired.
Let $\operatorname{PM}(\tau,\tau')$ be the set of partial matchings.
Each $\chi \in \operatorname{PM}(\tau, \tau')$ has an associated multigraph $(V(\chi), E(\chi))$
obtained by identifying matched vertices.
Write $\partial \chi$ for the image of $\partial\tau \cup \partial\tau'$ after identifying matched vertices ($t$ is left implicit).
$E(\chi)$ is a multiset whose elements have multiplicity 1 or 2 since the identification may produce double-edges.
Let $E_\cap(\chi)$ be a set containing one copy of each double-edge.
Let $V_\cap(\chi)\subseteq V(\chi)$ be the shared vertices, namely the images of matched pairs.
Define $\overline{\operatorname{PM}}(\tau,\tau')$ to consist of those
$\chi\in\operatorname{PM}(\tau,\tau')$ such that the graph
$(V_\cap(\chi),E_\cap(\chi))$ is root-connected (we allow it to consist of just the matched roots).
Equivalently, whenever two non-root vertices are matched, their parents must also be matched.
Starting from each root, the matching can split but not rejoin.

\begin{lemma}\label{lem:at-fourier}{
For every integer $t\ge0$,}
\begin{equation}\label{eq:at-fourier}
    \begin{aligned}
    a_{t+1} &= \sum_{\tau \in {\mathscr T_t}} {B_{t,\tau}} a_1^{{\ell}(\tau)} \prod_{u \in V(\tau) \setminus \partial\tau} \E[\Theta^{k_u + 1}]\\
    &= \sum_{\tau,\tau' \in {\mathscr T_t}}\sum_{\chi \in \overline{\operatorname{PM}}(\tau,\tau')} {B_{t,\tau}}{B_{t,\tau'}} a_1^{|\partial \chi|} \prod_{u \in V(\chi) \setminus \partial \chi} \E[\Theta^{\deg_{E(\chi)}(u) - 2 \deg_{E_\cap(\chi)}(u)}].
    \end{aligned}
\end{equation}
\end{lemma}
\begin{proof}
    By \cref{eq:zero-a,eq:weights} and the Hermite identities~\eqref{eq:hermite-first} and~\eqref{eq:hermite-second},
    \begin{align}\label{eq:a-recurrence}
        a_{t+1} = \E[\Theta \hat \eta_t(X_{a_t})] &= \sum_{k = 0}^d b_{t,k} \E[\Theta \He_{k}(X_{a_t}; a_t)] = \sum_{k=0}^d b_{t,k} \E[\Theta^{k+1}] a_t^k, \\
        a_{t+1} = \E[\hat \eta_t(X_{a_t})^2] &= \sum_{k,k'=0}^d b_{t,k}b_{t,k'} \E[\He_k(X_{a_t}; a_t) \He_{k'}(X_{a_t}; a_t)] \notag\\
        &= \sum_{k,k'=0}^d b_{t,k}b_{t,k'} \sum_{j=0}^{k \wedge k'} \binom kj \binom{k'} j j! \E[\Theta^{k+k'-2j}] a_t^{k+k'-j}\,.\notag
    \end{align}

    Let $M_{t+1}, N_{t+1}$ be the two right-hand sides of \cref{eq:at-fourier}.
    We claim the combinatorial summations $M_t, N_t$ solve the recurrence \eqref{eq:a-recurrence}.
    Indeed, we show the combinatorial recurrences
    \begin{equation}\label{eq:a-b-recurrence}
    \begin{aligned}
        M_{t+1}&=\sum_{k=0}^d b_{t,k}\E[\Theta^{k+1}]M_t^k,\\
        N_{t+1}&=\sum_{k,k'=0}^d b_{t,k}b_{t,k'}
        \sum_{j=0}^{k\wedge k'}\binom kj\binom{k'}j j! \,
        \E[\Theta^{k+k'-2j}] M_t^{k+k'-2j} N_t^j.
    \end{aligned}
    \end{equation}
    In the base case, ${M_1 = N_1 = a_1}$ for the singleton tree.
    Assuming \cref{eq:a-b-recurrence}, then
    $M_t = N_t = a_t$ holds by induction and \cref{eq:a-recurrence}. It remains to justify \cref{eq:a-b-recurrence}.
    
    The recurrence for $M_t$ is obtained by summing over the degree of the root of $\tau$.
    $\E[\Theta^{k+1}]$ accounts for the moment term at the root,
    while each branch can be summed separately to obtain $M_t$ per child.
    
    The recurrence for $N_{t}$ is obtained by first summing over the degree of the root $k, k'$ in $\tau$ and $\tau'$. There are three possible types of children of the root in $\chi$: those that belong to $\tau$ only, those that belong to $\tau'$ only, and those that are matched in $\chi$. Suppose $j$ pairs are matched. Since we consider trees with ordered children, there are
    $\binom{k}{j}\binom{k'}{j} j!$ ways to choose this partial matching. The moment term contributed by the root is \(\E[\Theta^{k+k'-2j}]\). For each of the \(k+k'-2j\) unmatched children, all descendants must remain unmatched due to the structural property of $\overline{\operatorname{PM}}(\tau, \tau')$, so the contribution from this branch is \(M_t\) from the first-moment calculation. On the other hand, each matched pair of children
    contributes \(N_t\) recursively. This yields the recurrence.
\end{proof}

\subsection{The case of symmetric prior}
\label{sec:upper-symmetric}

For now we focus on the case where $\pi$ is symmetric around zero. We assume that for \(\Theta\sim\pi\),
\begin{equation}
 \Theta\stackrel{\mathrm d}= -\Theta,\qquad \lambda\coloneqq\E[\Theta^2]>1,\qquad
 \E e^{\alpha|\Theta|}<\infty\quad\text{for some }\alpha>0.
 \label{eq:zero-assumptions}
\end{equation}
In particular, \(\E\Theta=0\), and all moments of $\Theta$ are finite.
Note that we assume $\lambda > 1$, meaning we are above the BBP threshold. The goal is to estimate \(\theta_1\theta_2\) given $\overline{\bY}$, where $\overline{\bY}$ is defined in \cref{eq:Y-bar}.

Recall from \cref{sec:se} that in this case, the state evolution recurrence is initialized at $q_1 := \lambda-\lambda^{-1}$, and $q_*$ is defined as the limit of the iteration $q_{t+1}=\Phi(q_t)$.

\begin{theorem}[Polynomial estimation with symmetric prior]
\label{thm:zero}
Assume \eqref{eq:zero-assumptions}. For every \(\eps>0\), there is
a constant \(C=C(\eps,\pi) > 0\) and a polynomial \(f(\overline{\bY}) = f_n(\overline{\bY})\) such that
\begin{equation}
 \deg(f)\le C\log n,\qquad
 \E[(f(\overline{\bY})-\theta_1\theta_2)^2]
 \le\lambda^2-q_*^2+\eps+o(1).
 \label{eq:zero-risk}
\end{equation}
\end{theorem}

\subsubsection{Constructing the polynomial}

Fix integers $T,d\ge1$. Let $a_1 := \lam - \lam^{-1}$ (matching $q_1$) and define the sequence $a_t$ by \eqref{eq:zero-a} along with the associated coefficients $b_{t,k}$. Due to symmetry of the prior, \(\hat\eta(\cdot;a,d)\) is an odd function. Since \(\He_k(\cdot;a)\) has the same parity as \(k\),
\[
 b_{t,k}=0\qquad\text{for every even }k.
\]
Therefore any tree $\tau \in \mathscr T_T$ with nonzero weight $B_{T,\tau}$ must have: an odd number of children at every
internal vertex, no leaf before depth \(T\), and an odd number $\ell(\tau)$ of boundary leaves.

Let \(\tau_1,\tau_2\in\mathscr T_T\), both with nonzero weights, and write $o_1,o_2$ for their respective roots. Let \(F=(\tau_1,\tau_2)\)
be their disjoint union, a forest with ordered roots \(o_1,o_2\). From above, $\ell(\tau_1)$ and $\ell(\tau_2)$ are odd, so the total number of boundary leaves $\ell(F) := \ell(\tau_1) + \ell(\tau_2)$ is even.

Pair up the $\ell(F)$ boundary vertices arbitrarily and connect each pair with a length-$L$ path for some even $L$ that will eventually be chosen to scale as $L = \Theta(\log n)$. There are $\ell(F)/2$ paths in total. We call these \emph{connector paths}. Each connector path adds $L$ edges and $L-1$ new vertices. Note that due to parity of the leaves, at least one path must connect $\tau_1$ to $\tau_2$. Write $H_L(F)$ for the resulting graph, including $F$ together with the added paths.

For any simple graph \(H\) with ordered roots $o_1,o_2$, define
\begin{equation}
\overline{Y}^H\coloneqq \sum_{\substack{\phi:V(H)\hookrightarrow[n]\\
\phi(o_1)=1,\ \phi(o_2)=2}}
\prod_{uv\in E(H)}\overline{Y}_{\phi(u)\phi(v)}.
\label{eq:zero-graph-sum}
\end{equation}
Here \([n]:=\{1,\ldots,n\}\), and \(V(H),E(H)\) denote the vertex
and edge sets. The symbol $\hookrightarrow$ means the labeling $\phi$ is injective.
Note that distinct labelings may give rise to the same monomial, if $H$ has non-trivial automorphisms.

For a choice of parameters \(T,d,L\), our polynomial estimator is
\begin{equation}
f(\overline{\bY})\coloneqq
\sum_{F=(\tau_1,\tau_2)\in\mathscr T_T^2}
B_{T,\tau_1}B_{T,\tau_2}
\left(\frac{n(\lambda-\lambda^{-1})^2}{\lambda^{L+1}}\right)^{\ell(F)/2}
\overline{Y}^{H_L(F)}.
\label{eq:zero-estimator}
\end{equation}
The sum omits pairs with \(B_{T,\tau_1}B_{T,\tau_2}=0\). Recall that $\tau_1,\tau_2$ were defined as \emph{ordered} trees, but the ordering is irrelevant for $\overline{Y}^{H_L(F)}$, where only the graph structure matters. This causes some identical terms to be repeated in the sum.  The polynomial $f(\overline{\bY})$ has degree at most \(2D_T+d^T L\) with \(D_T\) defined in \eqref{eq:degree}, as this is an upper bound on the total number of edges in $H_L(F)$. We will take $T,d$ constant and $L = \Theta(\log n)$, so the degree is $\Theta(\log n)$.

In the following, our goal is to show
\[ \EE[f(\overline{\bY}) \cdot \theta_1 \theta_2] = {a_{T+1}^2} + o(1) \qquad \text{and} \qquad \EE[f(\overline{\bY})^2] = {a_{T+1}^2} + o(1). \]
Together, these will imply the desired mean-squared error value (this will be shown in \Cref{sec:completing}).

\subsubsection{First moment}
\label{sec:zero-first-moments}

We compute $\EE[f(\overline{\bY}) \cdot \theta_1 \theta_2]$. For a simple graph \(H\)
with roots \(o_1,o_2\),
\begin{equation}
 \E[\overline{Y}^H\cdot\theta_1\theta_2]
 =\frac{(n-2)_{|V(H)|-2}}{n^{|E(H)|}}
   \E[\Theta^{\deg_H(o_1)+1}]\E[\Theta^{\deg_H(o_2)+1}]
   \prod_{v\notin\{o_1,o_2\}}\E[\Theta^{\deg_H(v)}],
 \label{eq:zero-first-graph}
\end{equation}
where $(n)_k := n(n-1)(n-2) \cdots (n-k+1)$ has been used to count the injective labelings.

For \(F=(\tau_1,\tau_2)\), let \(\partial F\) be its set of boundary
vertices. Substitute~\eqref{eq:zero-first-graph} into~\eqref{eq:zero-estimator} and note that the connector paths in $H_L(F)$ contain a total of $(L+1)\ell(F)/2$ degree-2 vertices (including the endpoints), which contributes $\EE[\Theta^2]^{(L+1)\ell(F)/2} = \lambda^{(L+1)\ell(F)/2}$ and cancels with the identical factor in~\eqref{eq:zero-estimator}:
\begin{align*}
&\E[f(\overline{\bY})\cdot\theta_1\theta_2]\\
&=\sum_{F=(\tau_1,\tau_2)\in\mathscr T_T^2}
B_{T,\tau_1}B_{T,\tau_2}\frac{(n-2)_{|V(H_L(F))|-2}}{n^{|E(H_L(F))|-\ell(F)/2}} (\lambda - \lambda^{-1})^{\ell(F)} \prod_{v\in V(F)\setminus\partial F}\E[\Theta^{k_v+1}] \\
&= o(1) + \sum_{F=(\tau_1,\tau_2)\in\mathscr T_T^2}
B_{T,\tau_1}B_{T,\tau_2} (\lambda - \lambda^{-1})^{\ell(F)} \prod_{v\in V(F)\setminus\partial F}\E[\Theta^{k_v+1}],
\intertext{where the second step used $|V(H_L(F))| - 2 = |E(H_L(F))| - \ell(F)/2$ together with $L = O(\log n)$ and $|\mathscr T_T| = O(1)$ to remove the $n$-dependence. The expression above factors across the two trees:}
&= o(1) + \left(\begin{aligned}
&\sum_{\tau\in\mathscr T_T}
B_{T,\tau}(\lambda - \lambda^{-1})^{\ell(\tau)} \prod_{v\in V(\tau)\setminus\partial\tau}\E[\Theta^{k_v+1}]
\end{aligned}\right)^2 \\
&= o(1) + {a_{T+1}^2}.
\end{align*}
The last equality is by \cref{lem:at-fourier}.

\subsubsection{Second moment}
\label{sec:zero-second-moments}

Now we move on to prove the second moment is $\EE[f(\overline{\bY})^2] = {a_{T+1}^2} + o(1)$.

\paragraph{Moment expansion.}
For two forests \(F=(\tau_1,\tau_2)\), \(F'=(\tau'_1,\tau'_2)\) of nonzero weight in~\eqref{eq:zero-estimator}, write
\[
 B_{F}\coloneqq B_{T,\tau_1}B_{T,\tau_2},\qquad
 H\coloneqq H_L(F),\qquad H'\coloneqq H_L(F').
\]
Expand the second moment as
\[ \E[f(\overline{\bY})^2] = \sum_{F,F'} B_{F} B_{F'} \left(\frac{n(\lambda-\lambda^{-1})^2}{\lambda^{L+1}}\right)^{\frac{\ell(F)+\ell(F')}{2}} \E[\overline{Y}^H \overline{Y}^{H'}], \]
and further expand
\begin{equation}\label{eq:EYY}
\E[\overline{Y}^H \overline{Y}^{H'}] = \sum_{\phi,\phi'} \E\left[\left(\prod_{uv\in E(H)}\overline{Y}_{\phi(u)\phi(v)}\right)\left(\prod_{uv\in E(H')}\overline{Y}_{\phi'(u)\phi'(v)}\right)\right],
\end{equation}
where $\phi,\phi'$ are injective labelings of $H,H'$ respectively, as in~\eqref{eq:zero-graph-sum}. It will be helpful to group the terms $(\phi,\phi')$ based on how the graphs $\phi(H)$ and $\phi'(H')$ intersect.
A pair $(\phi,\phi')$ induces the partial matching
\[ \chi := \{(u,u') \,:\, u \in V(H), \, u' \in V(H'), \, \phi(u) = \phi'(u')\}. \]

To sum over pairs $(\phi,\phi')$, we will first enumerate the possible partial matchings $\chi$, denoted $\chi \in \operatorname{PM}(H,H')$ following the notation of \cref{lem:at-fourier}. Then we count the number of pairs $(\phi,\phi')$ that correspond to a given $\chi$: there are $(n-2)_{|V(\chi)|-2}$ such pairs. The moment term $\E[\cdots]$ in~\eqref{eq:EYY} depends only on $\chi$. Specifically, conditional on the signal $\btheta$, each single-edge in $\chi$'s multigraph contributes a factor of the form
\[ \E[\overline{Y}_{ij}\mid\btheta]=\frac{\theta_i\theta_j}{n} \]
while each double-edge contributes
\[ \E[\overline{Y}_{ij}^2\mid\btheta]
 =\frac{\theta_i^2\theta_j^2}{n^2}+\frac1n. \]
Expanding the product over double-edges amounts to choosing $J \subseteq E_\cap$ where $E_\cap = E_\cap(\chi)$ contains one copy of each double-edge, and $J$ represents the double-edges that use the ``noise'' term $\frac 1n$ rather than the ``signal'' term $\frac{\theta_i^2\theta_j^2}{n^2}$. Putting it together,
\[ \E[\overline{Y}^H \overline{Y}^{H'}] = \sum_{\chi \in \operatorname{PM}(H,H')} \sum_{J \subseteq E_\cap} \frac{(n-2)_{|V(\chi)|-2}}{n^{|E(\chi)|-|J|}} \prod_{u \in V(\chi)}\E[\Theta^{\deg_{E(\chi)}(u) - 2 \deg_J(u)}], \]
and so
\begin{equation}\label{eq:second-summands}
\E[f(\overline{\bY})^2] = \sum_{F,F'} \sum_{\chi \in \operatorname{PM}(H,H')} \sum_{J \subseteq E_\cap} s(F,F',\chi,J)
\end{equation}
where
\[ \begin{aligned}
&s(F,F',\chi,J) := {}\\
&\quad B_{F}B_{F'}\left(\frac{n(\lambda-\lambda^{-1})^2}{\lambda^{L+1}}\right)^{\frac{\ell(F)+\ell(F')}{2}} \frac{(n-2)_{|V(\chi)|-2}}{n^{|E(\chi)|-|J|}} \prod_{u \in V(\chi)}\E[\Theta^{\deg_{E(\chi)}(u) - 2 \deg_J(u)}].\end{aligned} \]
We will next classify which of the summands $(F,F',\chi,J)$ in~\eqref{eq:second-summands} have a leading-order contribution and which are negligible.

\paragraph{Grouping by type.}

As previously, $V_\cap=V_\cap(\chi)$ denotes the shared vertices.
Regard $J$ as a graph on vertex set $V_\cap$, including isolated vertices.
Each connector
path in $H$ or $H'$ has a unique midpoint, at distance \(L/2\) from either end. In $H \sqcup H'$ there are $(\ell(F)+\ell(F'))/2$ such vertices, which we call \emph{connector midpoints}. A connector midpoint is called \emph{root-connected} if the corresponding vertex in $V(\chi)$ has a path in $J$ to one (or both) of the roots $o_1,o_2$.
With the conventions above, define the \emph{type} of a summand $(F,F',\chi,J)$ from~\eqref{eq:second-summands} to be
\[ (v,e,c,m), \]
where
\begin{itemize}[nosep]
 \item \(v=|V_\cap|\) is the number of shared vertices{, including both roots};
 \item \(e=|J|\) is the number of noise edges;
 \item \(c\) is the number of connected components of \(J\),
 including isolated vertices;
 \item \(m\) is the number of root-connected connector midpoints.
\end{itemize}
We declare the types
\begin{equation}
 (e{+2},e,2,0)
 \label{eq:zero-leading-types}
\end{equation}
to be \emph{good}, with \(e \ge 0\) allowed to vary. Equivalently,
the noise edges form two rooted trees connecting all shared vertices to the two
roots, and the trees do not reach any connector midpoint.
All other types are called \emph{bad}. We will show that the good terms have a constant-order contribution, asymptotically equal to ${a_{T+1}^2}$, while the bad terms contribute $o(1)$ in total.

\begin{lemma}[Leading contribution from good terms]
If $L = o(\sqrt{n})$ and $L \to \infty$, then
\label{lem:second-leading}
\[ \sum_{\mathrm{good}\;(F,F',\chi,J)} s(F,F',\chi,J) = {a_{T+1}^2} + o(1). \]
\end{lemma}

\begin{lemma}[Negligible contribution from bad terms]
\label{lem:second-neg}
There is a constant $C_0=C_0(T,d,\pi)>0$ such that, for every constant
$\kappa\ge C_0$, choosing $L=2\lceil \kappa\log n/2\rceil$ gives
\[ \sum_{\mathrm{bad}\;(F,F',\chi,J)} |s(F,F',\chi,J)| = o(1). \]
\end{lemma}

\noindent The following two subsections are devoted to proving these lemmas. Together they immediately imply the desired second moment bound: $\E[f(\overline{\bY})^2] = {a_{T+1}^2} + o(1)$.

A key ingredient for both proofs will be to understand the $n$-dependence in $s(F,F',\chi,J)$. To this end, we have
\[ 1 - O(L^2/n) \le \frac{(n-2)_{|V(\chi)|-2}}{n^{|V(\chi)|-2}} \le 1 \]
and, due to the structure of $H,H'$,
\begin{align*}
|E(\chi)| = |E(H)| + |E(H')| &= |V(H)| + |V(H')| - 4 + \frac{\ell(F)+\ell(F')}{2} \\
&= |V(\chi)| + |V_\cap| - 4 + \frac{\ell(F)+\ell(F')}{2}.
\end{align*}
Together, these imply that $s$ scales as $n^{e-v{+2}}$, namely
\begin{equation}\label{eq:n-dep}
(1-O(L^2/n))n^{e-v{+2}} \le n^{(\ell(F)+\ell(F'))/2} \frac{(n-2)_{|V(\chi)|-2}}{n^{|E(\chi)|-|J|}} \le n^{e-v{+2}}.
\end{equation}
This bound holds uniformly over all $F,F',\chi,J$.
The exponent $e-v{+2}$ may be positive or negative.

\subsubsection{Leading terms}

\begin{proof}[Proof of Lemma~\ref{lem:second-leading}]
For a good term $(F,F',\chi,J)$, $J$ is determined by $\chi$, namely $J = E_\cap(\chi)$. Furthermore, $\chi$ can be equivalently specified by the following data:
\begin{itemize}
    \item partial matchings between the trees: $\chi_1 \in \operatorname{PM}(\tau_1,\tau_1')$ and $\chi_2 \in \operatorname{PM}(\tau_2,\tau_2')$;
    \item for each leaf in $V_\cap(\chi)$ --- meaning $u \in V_\cap(\chi)$ at depth $T$ from one of the roots --- a value $\zeta_u \in \{0,1,\ldots,L/2-1\}$ representing the length of the path in $J$ that emanates from this leaf along the connector paths.
\end{itemize}
Note that $\zeta_u < L/2$ since $J$ is not allowed to reach a connector midpoint in a good term. Since $E_\cap(\chi)$ forms two rooted trees, we have $\chi_1\in \overline{\operatorname{PM}}(\tau_1,\tau_1')$ and $\chi_2 \in \overline{\operatorname{PM}}(\tau_2, \tau_2')$.
Write $\partial_\cap \chi_i := \partial \chi_i \cap V_\cap(\chi_i)$ for the shared leaves, and $\ell_\cap := |\partial_\cap \chi_1|+|\partial_\cap \chi_2|$ for their number. Let ${\boldsymbol{\zeta}} = (\zeta_u)_{u \in \partial_\cap \chi_1 \sqcup \partial_\cap \chi_2}$ denote the collection of all $\zeta_u$ values, and $|{\boldsymbol{\zeta}}| := \sum_u |\zeta_u|$. Recall from~\eqref{eq:n-dep} that the good terms have constant-order $n$-dependence, since $v=e{+2}$. Now we can write
\begin{align*}
&\sum_{\mathrm{good}\;(F,F',\chi,J)} s(F,F',\chi,J) \\
& \quad = \begin{multlined}[t][0.9\linewidth]
\sum_{F,F',\chi_1,\chi_2,{\boldsymbol{\zeta}}} B_{F}B_{F'}\left(\frac{\lambda-\lambda^{-1}}{\lambda^{(L+1)/2}}\right)^{\ell(F)+\ell(F')} \frac{(n-2)_{|V(\chi)|-2}}{n^{|V(\chi)|-2}}\\
\quad{}\times\prod_{u \in V(\chi)}\E[\Theta^{\deg_{E(\chi)}(u) - 2 \deg_{E_\cap(\chi)}(u)}],
\end{multlined}
\end{align*}
where the sum is over $\chi_i \in \overline{\operatorname{PM}}(\tau_i,\tau_i')$ and $\zeta_u \in \{0,1,\ldots,L/2-1\}$. We will first focus on evaluating the sum over ${\boldsymbol{\zeta}}$. The connector vertices (i.e., the vertices on the connector paths, including endpoints) have the following contribution to the moment term $\E[\cdots]$: each vertex such that both incident edges are shared contributes a factor of $1$, and every other connector vertex contributes a factor of $\lambda$. We can separate out the contribution from the connector paths as
\[ \begin{aligned}
&\prod_{u \in V(\chi)}\E[\Theta^{\deg_{E(\chi)}(u) - 2 \deg_J(u)}]\\
&\quad = \lambda^{(L+1)\frac{\ell(F)+\ell(F')}{2} - \ell_\cap - 2|{\boldsymbol{\zeta}}|} \hspace{-10pt} \prod_{u \in V(\chi_1) \setminus \partial \chi_1 \sqcup V(\chi_2) \setminus \partial \chi_2} \hspace{-10pt} \E[\Theta^{\deg_{E(\chi)}(u) - 2 \deg_{E_\cap(\chi)}(u)}].\end{aligned} \]
For fixed $F,F',\chi_1,\chi_2$, the ${\boldsymbol{\zeta}}$-dependent part of the good terms is
\begin{align*}
\sum_{\boldsymbol{\zeta}} \frac{(n-2)_{|V(\chi)|-2}}{n^{|V(\chi)|-2}} \cdot \lambda^{-2|{\boldsymbol{\zeta}}|} &= (1-O(L^2/n))\sum_{\boldsymbol{\zeta}} \lambda^{-2|{\boldsymbol{\zeta}}|} \\
&= (1-O(L^2/n))\left(\sum_{k=0}^{L/2-1} \lambda^{-2k}\right)^{\ell_\cap} = \left(\frac{\lambda}{\lambda-\lambda^{-1}}\right)^{\ell_\cap}+o(1).
\end{align*}
Notice that the term $q_1 = \lam - \lam^{-1}$ has appeared. The $\lambda$-dependent part of the good terms becomes
\[ \left(\frac{\lambda-\lambda^{-1}}{\lambda^{(L+1)/2}}\right)^{\ell(F)+\ell(F')} \lambda^{(L+1)\frac{\ell(F)+\ell(F')}{2} - \ell_\cap} \left(\frac{\lambda}{\lambda-\lambda^{-1}}\right)^{\ell_\cap} = (\lambda-\lambda^{-1})^{|\partial\chi_1|+|\partial\chi_2|}. \]
Now that we have evaluated the sum over ${\boldsymbol{\zeta}}$, the number of remaining terms $(F,F',\chi_1,\chi_2)$ is a constant. Combining the arguments above,
\begin{align*}
&\sum_{\mathrm{good}\;(F,F',\chi,J)} s(F,F',\chi,J) = o(1) \; + \\
& \qquad \sum_{F,F',\chi_1,\chi_2} B_{F}B_{F'}(\lambda-\lambda^{-1})^{|\partial \chi_1| + |\partial \chi_2|} \prod_{u \in V(\chi_1) \setminus \partial \chi_1 \sqcup V(\chi_2) \setminus \partial \chi_2} \E[\Theta^{\deg_{E(\chi)}(u) - 2 \deg_{E_\cap(\chi)}(u)}].
\intertext{Similar to the first moment analysis, this factors across the two roots:}
& \qquad = o(1) + \left(\begin{aligned}&\sum_{\tau,\tau' \in \mathscr T_T}\sum_{\chi \in \overline{\operatorname{PM}}(\tau,\tau')} B_{T,\tau}B_{T,\tau'}\\
&\quad{}\times(\lambda-\lambda^{-1})^{|\partial \chi|} \prod_{u \in V(\chi) \setminus \partial \chi} \E[\Theta^{\deg_{E(\chi)}(u) - 2 \deg_{E_\cap(\chi)}(u)}]\end{aligned}\right)^2 \\
& \qquad = o(1) + {a_{T+1}^2}.
\end{align*}
The last line is by \cref{lem:at-fourier}.
\end{proof}

\subsubsection{Negligible terms}

For fixed \(F,F'\), let
\(\mathcal S_{v,e,c,m}(F,F')\) denote the set of pairs $(\chi,J)$ with
\(\chi\in\operatorname{PM}(H,H')\), \(J\subseteq E_\cap(\chi)\)
such that \((F,F',\chi,J)\) has type \((v,e,c,m)\). Note there are a constant number of pairs $(F,F')$. We will bound the number of pairs
\((\chi,J)\) of each type (\Cref{lem:zero-count}) and bound the magnitude of each summand $s(F,F',\chi,J)$ of a given type (\Cref{lem:zero-size}).

\begin{lemma}[Number of terms of a given type]
\label{lem:zero-count}
There are constants \(C,K>0\), only depending on \(T,d\), such that
\leavevmode{
for every $F,F'\in\mathscr T_T^2$ with $B_{F}B_{F'}\ne0$ and every feasible type $(v,e,c,m)$,
}
\begin{equation}
 |\mathcal S_{v,e,c,m}(F,F')|\le (CL^K)^c.
 \label{eq:zero-type-number}
\end{equation}
\end{lemma}

\begin{proof}
Recall that we view $J$ as a graph on vertex set $V_\cap(\chi)$. We say a vertex in $V_\cap(\chi)$ is \emph{from the forests} if at least one of its two preimages in $H \sqcup H'$ belongs to a forest $F$ or $F'$, otherwise it is \emph{from the connector paths}. The number of vertices from the forests and their degrees with respect to $J$
are $O(1)$ depending on \(T,d\). Every
vertex from the connector paths has degree exactly two in both $H$ and $H'$, so it has degree
at most two in \(J\).

Every cycle in \(J\) includes a
vertex from the forests, since the connector paths alone contain no
cycles. This means we can split the edges of \(J\) into segments by cutting at every forest-vertex and every vertex of \(J\)-degree different from two (the two endpoints of a segment might be equal).
The interior vertices of each segment are all from the connector paths and have $J$-degree two. Record an oriented segment by its initial vertex
pair \((u,u') \in V(H) \times V(H')\), the neighbor-pair \((z,z')\) of the initial pair along the segment, and its length.
After the first edge, the segment must continue along internal connector vertices in both graphs, so these data determine the segment.
$V_\cap(\chi)$ might also include isolated vertices, which we record directly as \((u,u') \in V(H) \times V(H')\).
Together, this data recovers all of \(\chi\) and \(J\).

The number of isolated vertices is $\le c$ since, by definition, $c$ counts the connected components of $J$ including isolated vertices. We claim that the number of segments is also at most \(c + O(1)\).
Indeed, the number of segments incident to a forest vertex is $O(1)$ depending on $T,d$, but also every segment which is purely from the connector paths is its own connected component.

Now we count the ways to choose the overall encoding. Each graph $H$ or $H'$ has \(O(L)\) vertices and bounded maximum degree.
There are $O(L^3)$ possible encodings of a single segment:
the initial vertex pair has \(O(L^2)\) choices, its
first-neighbor pair has $O(1)$ choices, and the
segment length has \(O(L)\) choices. Each isolated vertex has \(O(L^2)\) choices. Bound the number of these objects by $O(c)$ to get the claimed bound
\((CL^K)^c\) on pairs \((\chi,J)\). Restricting $(\chi, J)$ to the specified values of
\(v,e,m\) can only decrease this count.
\end{proof}

\begin{lemma}[Magnitude of a summand of a given type]
\label{lem:zero-size}
There is a constant \(C>0\), only depending on \(T,d,\pi\), such
that every summand of type \((v,e,c,m)\) satisfies
\begin{equation}
 |s(F,F',\chi,J)|\le C^c n^{e-v{+2}}\lambda^{-2e}.
 \label{eq:zero-type-size}
\end{equation}
\end{lemma}

\begin{proof}
Using the definition of \(s(F,F',\chi,J)\) following
\eqref{eq:second-summands}, along with the upper bound on the $n$-dependent part from
\eqref{eq:n-dep}, it suffices to prove the bound
\begin{equation}\label{eq:neg-goal}
 \lambda^{-(L+1)(\ell(F)+\ell(F'))/2}
 \prod_{u\in V(\chi)}
 \left|\E[\Theta^{\deg_{E(\chi)}(u)-2\deg_J(u)}]\right|
 \le C^c\lambda^{-2e}.
\end{equation}
Here we have dropped the factor $B_{F}B_{F'}(\lambda-\lambda^{-1})^{\ell(F)+\ell(F')}$ because it is bounded by a constant, which can be absorbed into $C^c$ since $c \ge 1$.

Let \(U\subseteq V(\chi)\) consist of the vertices for which at least one of the $\le 2$ preimages in $H \sqcup H'$ is a forest vertex (i.e., belonging to $F$ or $F'$). For each \(u\in V(\chi)\), write
\[
 j_u\coloneqq\deg_J(u),\qquad
 M_u\coloneqq
 \left|\E[\Theta^{\deg_{E(\chi)}(u)-2j_u}]\right|,
\]
where \(j_u=0\) when \(u\) is unshared. All constants below only
depend on \(T,d,\pi\), and may increase from line to line.

We first handle the vertices in \(U\). Since \(F,F'\) have bounded
size and \(H,H'\) have bounded maximum degree,
\[
 \prod_{u\in U}M_u\le C,\qquad
 \sum_{u\in U}j_u\le C.
\]

Now we handle the vertices outside $U$. Every such $u \in V(\chi) \setminus U$ has either one or two preimages, all
in connector interiors. Each preimage has degree two, giving
the following moment factors:
\[
 M_u=
 \begin{cases}
  \lambda,&u\text{ is unshared},\\
  \E[\Theta^4],&u\text{ is shared and }j_u=0,\\
  \lambda,&u\text{ is shared and }j_u=1,\\
  1,&u\text{ is shared and }j_u=2.
 \end{cases}
\]
A shared vertex with \(j_u=0\) is isolated in \(J\), so there
are at most \(c\) such vertices. Bound \(\E[\Theta^4]\le C \le C\lambda^2\). Now in every case above,
the exponent of \(\lambda\) is the number of
preimages of \(u\) minus \(j_u\), with an extra constant factor $C^c$ in front.

To sum these exponents, let \(N_0\) be the number of preimages
of vertices outside \(U\).
We have
\[
 \prod_{u\notin U}M_u
 \le C^c\lambda^{N_0-\sum_{u\notin U}j_u}.
\]
All preimages contributing to $N_0$ are interior connector vertices, whose total number in \(H\sqcup H'\) is
\[
 N\coloneqq\frac{(L-1)(\ell(F)+\ell(F'))}{2},\qquad N_0\le N.
\]
Now \(\sum_{u\in V(\chi)}j_u=2e\), since each noise edge has
two endpoints. The vertices in \(U\) account for only a constant
part of this sum, so \(\sum_{u\notin U}j_u\ge2e-C\).
Using \(\lambda>1\) we obtain, after increasing \(C\),
\[
 \prod_{u\in V(\chi)}M_u\le C^c\lambda^{N-2e}.
\]

Returning to our goal in~\eqref{eq:neg-goal}, we have
\[ \lambda^{-(L+1)(\ell(F)+\ell(F'))/2}
 \prod_{u\in V(\chi)}
 M_u
 \le C^c \lambda^{-(L+1)(\ell(F)+\ell(F'))/2 + N - 2e} \le C^c\lambda^{-2e} \]
as desired, where we have used the definition of $N$ and increased $C$ in the last step.
\end{proof}

\begin{proof}[Proof of \Cref{lem:second-neg}]
Combine \Cref{lem:zero-count,lem:zero-size} and increase
\(C,K\) when necessary to get
\begin{equation}
 \sum_{(\chi,J)\in\mathcal S_{v,e,c,m}(F,F')}
       |s(F,F',\chi,J)|
 \le (CL^K)^c n^{e-v{+2}}\lambda^{-2e}.
 \label{eq:zero-type-contribution}
\end{equation}
We now sum this bound over bad types $(v,e,c,m)$.
Classify the bad types into two disjoint cases:
\begin{enumerate}[label=(\roman*),nosep]
 \item \(m=0\), \(J\) is a forest, and \(c\ge3\);
 \item \(m\ge1\) or $J$ contains a cycle.
\end{enumerate}
To see that this exhausts the possibilities, observe that if $m = 0$, then the two roots are in different components of $J$, so $c \geq 2$. The case $c=2$ has either $J$ a forest (the good case) or $J$ containing a cycle (case ii).
The case $c \geq 3$ has either $J$ a forest (case i) or $J$ contains a cycle (case ii).
The remaining case $m \geq 1$ goes into case (ii).

Every type in class (i) is a forest and so has
\begin{equation}\label{eq:floating-components}
    e-v = \,{-c} < 0
\end{equation}
Fix $F, F'$ of which there are only constantly many choices.
The bound in \cref{eq:zero-type-contribution} satisfies
\begin{align*}
 &\quad\sum_{\substack{(v,e,c,m)\\\text{class (i)}}} \sum_{(\chi, J) \in \calS_{v,e,c,m}(F,F')}|s(F,F',\chi,J)|\\
 &\leq \sum_{\substack{(v,e,c,m)\\\text{class (i)}}} (CL^K)^c n^{2-c} \leq \sum_{c \geq 3} (CL^K)^c n^{2-c} = O(L^{3K}/n)\,.
 \numberthis\label{eq:class-i-bound}
\end{align*}
To pass from the summation over types to the summation over $c$, we have observed
that the number of possibilities for $v,e,m$ is at most $O(L)$ which can be handled by increasing the constant $K$.

For every type in class (ii) we use the upper bound, for a constant $C$,
\begin{equation}
 e-v\le C -c.
 \label{eq:zero-edge-vertex-bound}
\end{equation}
Indeed, we can remove the constantly many edges of $F$ to make $J$ a collection of paths,
at which point the previous equality for forests $e-v = \,{-c}$ applies.

We claim $e \geq L/2$ for every type in class (ii).
If \(m\ge1\), a path in \(J\) from a root to a connector
midpoint has at least \(L/2\) edges. If \(J\) contains a cycle,
that cycle must traverse an entire connector path, giving at least
\(L\) edges. Thus we have $e \geq L/2$.
Summing the bound \cref{eq:zero-type-contribution} in this case,
\begin{align*}
 &\quad\sum_{\substack{(v,e,c,m)\\\text{class (ii)}}} \sum_{(\chi, J) \in \calS_{v,e,c,m}(F,F')}|s(F,F',\chi,J)|\\
 &\leq \sum_{\substack{(v,e,c,m)\\\text{class (ii)}}} (CL^K)^c n^{C + 2-c}\lam^{-2e} \leq O(n^C L^{K} \lam^{-L})\,.
 \numberthis\label{eq:class-ii-bound}
\end{align*}
With $L=2\lceil\kappa\log n/2\rceil$, we have $L=O(\log n)$ and
$\lambda^{-L}\le n^{-\kappa\log\lambda}$.
Thus \cref{eq:class-i-bound} is $o(1)$, and \cref{eq:class-ii-bound} is also $o(1)$
for every sufficiently large fixed $\kappa$, since $\lambda>1$.
This proves the lemma.
\end{proof}

\subsubsection{Completing the proof}
\label{sec:completing}

\begin{proof}[Proof of \Cref{thm:zero}]
Recall that $q_t\to q_*$ for the initialization
$q_1=\lambda-\lambda^{-1}$ (see \Cref{lem:zero-fixed-point}).

Fix \(\eps>0\). First choose a finite \(T\ge1\) such that
\({q_{T+1}^2}\ge q^2_*-\eps/2\). Next, by
Lemma~\ref{lem:approx-scalar} and $a_1 = q_1$, choose a finite \(d\ge1\) such that
\(|{a_{T+1}^2}-{q_{T+1}^2}|\le\eps/2\). Then
\[
 {a_{T+1}^2}\ge q^2_*-\eps.
\]
With $T,d$ fixed, choose a sufficiently large fixed $\kappa$ and set
$L=2\lceil\kappa\log n/2\rceil$ so that
\Cref{lem:second-leading,lem:second-neg} apply.
The degree of the estimator $f$ in \eqref{eq:zero-estimator} is \(\deg f\le2D_T+d^T L=O(\log n)\).
As shown in \cref{sec:zero-first-moments,sec:zero-second-moments},
\[
 \E[f(\overline{\bY}) \cdot \theta_1\theta_2] = {a_{T+1}^2} + o(1) \qquad \E[f(\overline{\bY})^2]={a_{T+1}^2}+o(1).
\]
Then we have, using that
\(\E[(\theta_1\theta_2)^2]=\lambda^2\),
\begin{align*}
 \E[(f(\overline{\bY})-\theta_1\theta_2)^2]
 &=\lambda^2-2\E[f(\overline{\bY})\cdot\theta_1\theta_2]
             +\E[f(\overline{\bY})^2]\\
 &=\lambda^2-{a_{T+1}^2}+o(1)\\
 &\le\lambda^2-q_*^2+\eps+o(1). \qedhere
\end{align*}
\end{proof}

\subsection{The case of nonzero mean}
\label{sec:nonzero-mean}

Now assume that the prior $\Theta \sim \pi$ satisfies
\begin{equation}
 \mu\coloneqq\E\Theta\ne0,\qquad
 \lambda\coloneqq\E[\Theta^2],\qquad
 \E e^{\alpha|\Theta|}<\infty\quad\text{for some }\alpha>0.
 \label{eq:nonzero-assumptions}
\end{equation}
There is no restriction on $\lambda$ in this case. Recall that the state evolution recurrence is initialized at $q_1:=\mu^2>0$,
and the iterates $q_{t+1}=\Phi(q_t)$ converge to $q_*$ (see \Cref{lem:zero-fixed-point}).

\begin{theorem}[Polynomial estimation with nonzero mean]
\label{thm:nonzero}
Assume~\eqref{eq:nonzero-assumptions}. For every $\eps>0$, there is
a constant $C=C(\eps,\pi)>0$ and a polynomial $f(\overline{\bY})=f_n(\overline{\bY})$ such that
\begin{equation}
 \deg(f)\le C,\qquad
 \E[(f(\overline{\bY})-\theta_1\theta_2)^2]\le\lambda^2-q_*^2+\eps+o(1).
 \label{eq:nonzero-risk}
\end{equation}
\end{theorem}

\subsubsection{Constructing the polynomial}

Fix integers $T,d\ge1$. Let $a_1:=\mu^2$ (matching $q_1$) and define the sequence
$a_t$ by \eqref{eq:zero-a} along with the associated coefficients $b_{t,k}$.
The estimator is again a sum of weighted forest polynomials:
\begin{equation}
 f(\overline{\bY})\coloneqq\sum_{F=(\tau_1,\tau_2)\in\mathscr T_T^2}
 B_{T,\tau_1}B_{T,\tau_2}\mu^{\ell(F)}\overline{Y}^F.
 \label{eq:nonzero-estimator}
\end{equation}
Compared to the previous case, the weights on the non-boundary vertices of $F$ are defined in the same way,
each boundary vertex of $F$ now has a factor of $\mu$ instead of $\lam - \lam^{-1}$,
and we do not need any paths extending from the boundary.
The degree of the estimator is at most $2D_T$
which is a constant when $T,d$ are fixed,
where $D_T$ is from~\eqref{eq:degree}.

\subsubsection{First and second moments}

We adapt the moment calculations from \cref{sec:upper-symmetric}, explaining only the differences. We will show
\begin{equation}
\E[f(\overline{\bY})\theta_1\theta_2]={a_{T+1}^2}+o(1),\qquad
 \E[f(\overline{\bY})^2]={a_{T+1}^2}+o(1).
 \label{eq:nonzero-moments}
\end{equation}

\paragraph{First moment.}
The calculation is exactly the same as \cref{sec:zero-first-moments}
except that the contribution of each leaf is now $\mu\E[\Theta] = \mu^2 = a_1$
instead of $\frac{\lam - \lam^{-1}}{\lam^{(L+1)/2}}\left(\E[\Theta^2]\right)^{(L+1)/2} = \lam - \lam^{-1}$.

\paragraph{Second moment.}
This calculation is significantly easier than the symmetric case.
Expand $\E[f(\overline{\bY})^2]$ as a sum over $F,F',\chi,J$ as in~\eqref{eq:second-summands}.
The exponent in $n$ for each term is again $|J| - |V_\cap(\chi)| + 2$.
Now, $J$ is acyclic since $F,F'$ are acyclic,
so this exponent is nonpositive
and is zero exactly if $\chi \in \overline{\operatorname{PM}}(F, F')$ and $J = E_\cap(\chi)$.
These leading terms evaluate to ${a_{T+1}^2}$ by \cref{lem:at-fourier}.
There are only $O(1)$ non-leading terms,
so the error is $O(n^{-1})$.

\paragraph{Completing the proof.}
With $q_1=a_1=\mu^2$, choose $T$ and $d$ as in \cref{sec:completing}.
The moment identities \eqref{eq:nonzero-moments} then imply the MSE bound
\eqref{eq:nonzero-risk}, proving \Cref{thm:nonzero}.

\subsection{The case of zero mean without symmetry}
\label{sec:zero-asymmetric}

In the final case, assume that
$\pi$ is not symmetric about zero and that
\begin{equation}
 \E\Theta=0,\qquad \lambda\coloneqq\E[\Theta^2]>1,\qquad
 \E e^{\alpha|\Theta|}<\infty\quad\text{for some }\alpha>0.
 \label{eq:asymmetric-assumptions}
\end{equation}

Recall that the state evolution recurrence is initialized at $q_1:=\lambda-\lambda^{-1}>0$,
and the iterates $q_{t+1}=\Phi(q_t)$ converge to $q_*$ (see \Cref{lem:zero-fixed-point}).

\begin{theorem}[Polynomial estimation with zero mean without symmetry]
\label{thm:asymmetric}
Assume~\eqref{eq:asymmetric-assumptions} and that $\pi$ is not
symmetric about zero. For every $\eps>0$, there is a constant
$C=C(\eps,\pi)>0$ and a polynomial $f(\overline{\bY})=f_n(\overline{\bY})$ such that
\begin{equation}
 \deg(f)\le C\log n,\qquad
 \E[(f(\overline{\bY})-\theta_1\theta_2)^2]\le\lambda^2-q_*^2+\eps+o(1).
 \label{eq:asymmetric-risk}
\end{equation}
\end{theorem}

\subsubsection{Constructing the polynomial}

Fix integers $T,d\ge1$. Let $a_1:=\lambda-\lambda^{-1}$ (matching $q_1$) and define the sequence
$a_t$ by \eqref{eq:zero-a} along with the associated coefficients $b_{t,k}$.
Unlike the symmetric case, all of the weights $b_{t,k}$ could potentially be nonzero,
and all trees $\tau \in \mathscr T_T$ could have nonzero weight $B_{T,\tau}$.

By the asymmetry of the prior, there is an integer $h\ge1$ such that \cite[Theorem 30.1]{billingsley1995probability}
\begin{equation}
 \nu\coloneqq\E[\Theta^{2h+1}]\ne0.
 \label{eq:asymmetric-odd-moment}
\end{equation}
Given a forest $F \in \mathscr T_T^2$, attach a path of length $L$
to each $u \in \partial F$
with a new ``terminal vertex'' $z_u$ at the other end of the path.
At $z_u$, attach $h$ cycles, each of length $L$.
Denote the resulting graph by $H_L(F)$.
For example, if $\E[\Theta^3]\ne0$, we take $h=1$, so each
boundary leaf receives a path ending at a single cycle.
The estimator is
\begin{equation}
 f(\overline{\bY})\coloneqq\sum_{F=(\tau_1,\tau_2)\in\mathscr T_T^2}
 B_{T,\tau_1}B_{T,\tau_2}
 \left(\frac{n^h(\lambda-\lambda^{-1})}
 {\nu\lambda^{(L-1)h + L}}\right)^{\ell(F)}\overline{Y}^{H_L(F)}.
 \label{eq:asymmetric-estimator}
\end{equation}
The degree of $f$ is at most
\[
 2D_T+2(h+1)d^T L=O(\log n)
\]
when $T,d,h$ are fixed and $L=\Theta(\log n)$.

\subsubsection{First and second moments}

We adapt the moment calculations from \cref{sec:upper-symmetric}, explaining only the differences.
For a sufficiently large fixed $\kappa=\kappa(T,d,\pi)$, take $L=2\lceil\kappa\log n/2\rceil$.
We will show
\begin{equation}
 \E[f(\overline{\bY})\theta_1\theta_2]={a_{T+1}^2}+o(1),\qquad
 \E[f(\overline{\bY})^2]={a_{T+1}^2}+o(1).
 \label{eq:asymmetric-moments}
\end{equation}

\paragraph{First moment.}
The calculation is exactly the same as \cref{sec:zero-first-moments}
except that the contribution of each leaf was previously $\frac{\lam - \lam^{-1}}{\lam^{(L+1)/2}}\left(\E[\Theta^2]\right)^{(L+1)/2} = \lam - \lam^{-1}$
and is now
\[
    \frac{\lam - \lam^{-1}}{\nu \lam^{(L-1)h + L}} \left(\E[\Theta^2]\right)^{(L-1)h + L}\E[\Theta^{2h+1}] = \lam - \lam^{-1}\,.
\]

\paragraph{Second moment: leading terms.}
Expand $\E[f(\overline{\bY})^2]$ as a sum over $F,F',\chi,J$ as in~\eqref{eq:second-summands}.
Group the terms by their type $(v,e,c,m)$ with $v,e,c$ the same as before: $v=|V_\cap(\chi)|$, $e=|J|$,
$c$ equals the number of components of $J$. Modify $m$ to be the number of terminal vertices $z_u$
connected by $J$ to either root, in place of midpoints of the connector paths.

The power of $n$ of a type is again $n^{e-v{+2}}$ and the good types
are again $(e{+2},e,2,0)$.
In a good type, $J$ forms two disjoint rooted trees,
and neither tree reaches a terminal vertex.
Up to multiplying $L$ by a factor 2, the good types are the
same as the good types for the symmetric prior, since in the two cases $J$ is not long enough to encounter the terminal vertex
or the midpoint of the connector path.
By the same analysis, the total contribution of the good terms is ${a_{T+1}^2} + o(1)$.

\paragraph{Second moment: negligible terms.}
The second moment analysis for a symmetric prior consisted of an
upper bound on the number of pairs with a given type (\Cref{lem:zero-count}),
an upper bound on the magnitude of each bad type (\Cref{lem:zero-size}),
and a final summation over the bad types.
The same qualitative bounds \Cref{lem:zero-count,lem:zero-size} continue to hold for the new definition of type:
\begin{itemize}
    \item For \Cref{lem:zero-count}, split $J$ into segments along the $O(1)$
    terminal vertices as well as the forest vertices. The degrees of all vertices
    remain constant. The same encoding now shows $|\calS_{v,e,c,m}(F,F')| \leq (CL^K)^c$.
    \item For \Cref{lem:zero-size}, add the
    $O(1)$ terminal vertices to the exceptional set $U$.
    Since their degree is constant, their moments contribute only a constant factor.
    The analysis of the degree-2 vertices along the paths and cycles
    is the same, showing $|s(F,F',\chi,J)|\le C^c n^{e-v{+2}}\lambda^{-2e}$.
\end{itemize}

\noindent
 Multiplying these together, the contribution of each bad type is at most $(CL^K)^c n^{e-v{+2}}\lambda^{-2e}$.
It remains to sum over the possible bad types.
Every bad term again belongs to one of two possible classes:
\begin{enumerate}[label=(\roman*),nosep]
    \item $m=0$, $J$ is a forest, and $c \geq 3$;
    \item $m\geq 1$ or $J$ contains a cycle.
\end{enumerate}
The second case is only possible when $J$ contains either a full path from the root to a terminal vertex, or a full cycle hanging off a terminal vertex. In both cases, $e \geq L$.
The same analysis as before now shows the bad types sum to $o(1)$.

\paragraph{Completing the proof.}
With $q_1=a_1=\lambda-\lambda^{-1}$, choose $T$ and $d$ as in \cref{sec:completing},
and then take $\kappa$ sufficiently large that \eqref{eq:asymmetric-moments} holds.
These moment identities imply the MSE bound \eqref{eq:asymmetric-risk}, proving \Cref{thm:asymmetric}.

\subsection{Completing the upper bound}
\label{sec:upper-main-theorem}

\begin{proof}[Proof of \Cref{thm:lowdeg}-\ref{item:upper}]
Let $\mu:=\E\Theta$, $\lambda:=\E\Theta^2$, and fix $\eps>0$.
We prove the upper bound under \eqref{eq:upper-moment-assumption}, which holds in particular for every subgaussian prior. \leavevmode{
It suffices to construct an estimator for $\theta_1\theta_2$ from $\overline{\bY}$; see \cref{eq:matrix:scalar:MMSE}.
}

If $\mu\ne0$, \Cref{thm:nonzero} supplies a polynomial of degree at most $C_{\eps,\pi}$ with MSE at most $\lambda^2-q_*^2+\eps+o(1)$.
If $\mu=0$ and $\lambda>1$, use \Cref{thm:zero} when $\pi$ is symmetric and \Cref{thm:asymmetric} otherwise. In either case, the degree is at most $C_{\eps,\pi}\log n$, and the MSE is again at most $\lambda^2-q_*^2+\eps+o(1)$.
If $\mu=0$ and $\lambda\le1$, \Cref{lem:zero-fixed-point} shows $q_*=0$. The trivial constant estimator $f(\overline{\bY})=0$ achieves MSE $\E[(\theta_1\theta_2)^2]=\lambda^2=\MSE_\AMP$ and has degree zero.
\end{proof}

\appendix

\section{Polynomial approximation to AMP}
\label{sec:trees}

\newcommand{\depth}{\textnormal{depth}}

Let $\overline{\bY} := \frac{\bY}{\sqrt{n}}$ and $\overline{\bZ} := \frac{\bZ}{\sqrt{n}}$. Recall the AMP iteration from \eqref{eq:amp},
\begin{align*}\numberthis\label{eq:amp2}
\bx^{(t+1)} &= \overline{\bY} \hat \btheta^{(t)} - \langle \eta_t'(\bx^{(t)})\rangle \cdot \hat \btheta^{(t-1)}\,, \\
    \hat \btheta^{(t+1)} &:= \eta_{t+1}(\bx^{(t+1)})\,.
\end{align*}
For simplicity, we will take the initialization to be $\hat\btheta^{(0)}=\omega\mathbf1$ and $\bx^{(1)}=\omega\overline{\bY}\mathbf1$ for $\omega \in \R$. This is most relevant for the case of nonzero-mean prior, where the case $\omega = \mu$ aligns with~\eqref{eq:amp-init}.

Let $\calT$ denote the set of trees on $[n]$ rooted at vertex 1.
Our goal in this section is to heuristically compute approximations to $\hat\btheta^{(t)}$ in the ``$\bZ$ basis'' and the ``$\bY$ basis'', meaning expansions of the form
\begin{equation}\label{eq:basis-approximations}
    \hat \btheta^{(t)}_1 \approx \sum_{\tau \in \calT} \widetilde A_{t,\tau}(\btheta) \overline{\bZ}^\tau, \qquad \qquad \hat\btheta^{(t)}_1 \approx \sum_{\tau \in \calT} \widetilde B_{t,\tau}\overline{\bY}^\tau\,.
\end{equation}
We allow arbitrary denoisers $\eta_t$, not necessarily the Bayes denoisers.

{The $\bZ$ basis expansion has orthogonal components, while the $\bY$ basis expansion does not.}
For each value of $t$ and $n$, assuming suitable regularity, we can obtain the $Z$ basis coefficients $\widetilde A_{t,\tau}(\btheta)$ by fixing $\theta$ and expanding $\hat \btheta^{(t)}$ in the Hermite polynomial basis with respect to $\bZ$.
But in general these coefficients are rather complicated functions of the AMP iteration up to time $t$.
In the limit $n \to \infty$, the coefficients simplify in terms of the state evolution parameters $m_t$ and $q_t$ from \eqref{eq:state-evolution},
and it is these simplified expressions that we wish to heuristically compute.

Let $\calT_{\le t}$ be the set of rooted trees with depth at most $t$ and write $\partial_t\tau$ for the vertices at depth exactly $t$.
For a fixed time $t$ and tree $\tau \in \calT_{\le t}$, for $v \in V(\tau)$ let $k_v$ be its number of children and $t_v = t - \depth(v)$ be the time associated with $v$.
Define $A_{t,k}(\Theta)$ and $B_{t,k}$ by the following rescaled Hermite expansions for the denoisers:
\begin{align}\label{eq:hermite-z}
    \eta_t(m_t \Theta + \sqrt{q_t}\,G) &= \sum_{k = 0}^\infty \frac{A_{t,k}(\Theta)}{k!} \He_k(\sqrt{q_t}\,G; q_t) = \sum_{k = 0}^\infty \frac{A_{t,k}(\Theta)}{k!} q_t^{k/2} \He_k(G)\\
    \eta_t(x) &= \sum_{k=0}^\infty \frac{B_{t,k}}{k!} \He_k(x;q_t)\,. \label{eq:hermite-y}
\end{align}
Then for fixed $t$ and $n \to \infty$, we will calculate the coefficients in \eqref{eq:basis-approximations} to be, for $\tau \in \calT_{\le t}$:
\begin{align}\label{eq:coefficient-answer}
    \widetilde A_{t,\tau}(\btheta) &= \omega^{|\partial_t\tau|}\prod_{v \in V(\tau) \setminus \partial_t \tau} A_{t_v,k_v}(\theta_v), \\
    \widetilde B_{t,\tau} &= \omega^{|\partial_t\tau|}\prod_{v \in V(\tau)\setminus \partial_t\tau} B_{t_v,k_v}\,. \notag
\end{align}
Each boundary leaf is associated with the initialization value $\omega$.

It should be noted that the tree terms $\calT$ make the dominant contribution in $\bx^{(t)}_1$ up to $o_n(1)$ error, hence we don't use cyclic terms in \eqref{eq:basis-approximations}.
To see this, for an \emph{unlabeled} rooted graph $\alpha$ with root $o$, define the polynomials \begin{equation}
\overline{\bZ}^\alpha \coloneqq \sum_{\substack{\phi:V(\al)\hookrightarrow[n]\\
\phi(o)=1}}
\prod_{uv\in E(\al)}\overline{Z}_{\phi(u)\phi(v)}, \qquad \overline{\bY}^\alpha \coloneqq \sum_{\substack{\phi:V(\al)\hookrightarrow[n]\\
\phi(o)=1}}
\prod_{uv\in E(\al)}\overline{Y}_{\phi(u)\phi(v)}.
\label{eq:z-polynomial}
\end{equation}
Here \([n]:=\{1,\ldots,n\}\) and the symbol $\hookrightarrow$ means the labeling $\phi$ is injective.
A direct calculation shows $\Var(\overline{\bZ}^\alpha) = \Theta(n^{-\ex(\alpha)})$ and $\Var(\overline{\bY}^\alpha) = \Theta(n^{-\ex(\alpha)})$ for each constant-size connected simple $\alpha$, where $\ex(\alpha) = |E(\alpha)| - |V(\alpha)| + 1$.
The only connected rooted graphs with $\ex(\tau) = 0$ are trees, hence these graphs are the dominant ones among connected rooted graphs.\footnote{
If $\al$ has a multiedge of multiplicity at least three, or if the double edges do not form ``hanging trees'' off a tree of single edges, then
the variance of $\overline{\bZ}^\al$ is again $o_n(1)$.
If the double edges form ``hanging trees'',
then $\overline{\bZ}^\alpha - \overline{\bZ}^\beta$
has variance $o_n(1)$ where $\beta$ equals $\alpha$ with the hanging trees removed.
This fully classifies the constant-size rooted connected $\alpha$ and shows that $\overline{\bZ}^\tau$ for trees $\tau$ spans the non-negligible space.
}

We demonstrate two heuristic methods for calculating \eqref{eq:coefficient-answer}: direct combinatorics on trees, and Gaussian integration by parts plus standard ``cavity method'' assumptions.
The first works for both the $\bZ$ and $\bY$ bases, while the second works only for the $\bZ$ basis, but it can also calculate the lower-order coefficients in the $\bZ$ expansion (i.e., the $n \to \infty$ limit of $\widetilde A_\alpha(\btheta)$ on graphs $\alpha$ which may have cycles).
To summarize, while we expect that \eqref{eq:basis-approximations} holds with the coefficients from \eqref{eq:coefficient-answer} in an appropriate rigorous sense,
our remaining goal is to heuristically derive the coefficients \eqref{eq:coefficient-answer}.

\subsection{Heuristic \texorpdfstring{$Z$}{Z} basis coefficients via tree combinatorics}
\label{app:expand-amp-Z}

We use a direct combinatorial approach following \cite{jones2025fourier}.
Maintaining the ``tree approximations'' $\hat \btheta^{(t)}_1 \approx \sum_{\tau \in \calT} \widetilde A_{t,\tau}(\btheta) \overline{\bZ}^\tau$ and $\bx^{(t)}_1 \approx \sum_{\tau \in \calT} \widetilde A'_{t,\tau}(\btheta) \overline{\bZ}^\tau$, we observe how AMP affects these representations.
The AMP operations act as follows.
\begin{enumerate}[label=\arabic*.]
\item The tree approximations to $\bx^{(t)}$ and $\hat \btheta^{(t)}$ contain trees of height at most $t$. The factor for a vertex $v$ is determined on iteration $t_v$, and the initialization corresponds to leaves at depth $t$.

\item The Onsager-corrected matrix multiplication is
\[
    \bx^{(t+1)} = \overline{\bY} \eta_t(\bx^{(t)}) - (\text{Onsager correction}) = \underbrace{\frac{\ip{\btheta}{\eta_t(\bx^{(t)})}}{n} \btheta}_{\text{term 1}} + \underbrace{\overline{\bZ} \eta_t(\bx^{(t)}) - (\text{Onsager correction})}_{\text{term 2}}\,.
\]
Term 1 is $\frac{\ip{\btheta}{\eta_t(\bx^{(t)})}}{n} \btheta \approx m_{t+1} \btheta$ for $m_{t+1}$ defined in \eqref{eq:state-evolution}.
Multiplying by $\overline{\bZ}$ in term 2 can extend the root of trees in $\eta_t(\bx^{(t)})$ by a single edge, or ``backtrack'' by one edge (all other possibilities create negligible cyclic graphs).
The Onsager correction is designed to equal the ``backtracking'' terms, up to $o_n(1)$ error.
Hence, the tree approximation for term 2 is formed
by extending the root of all trees in $\eta_t(\bx^{(t)})$ by a single edge.

\item Applying a function $\eta_t$ increases the branching factor
of the root, with coefficient $A_{t,k}(\theta_v)$ when the root $v$ has degree $k$.
The coefficient $A_{t,k}(\theta_v)$ will be explained below.
\end{enumerate}
Iterating these operations via \eqref{eq:amp2}, we inductively obtain formulas for the tree approximation of $\bx^{(t)}$ and $\hat \btheta^{(t)}$,
the latter having the $Z$ basis coefficients claimed in \eqref{eq:coefficient-answer},
\begin{align}\label{eq:tree-apx-Z}
    \bx^{(t)}_1 &\approx m_t \theta_1 + \sum_{\tau \in \calT_{\le t}: \deg_\tau(1) = 1}\omega^{|\partial_t\tau|}\prod_{v \in V(\tau) \setminus (\{1\} \cup \partial_t\tau)} A_{t_v,k_v}(\theta_v) \overline{\bZ}^\tau\\
    \hat \theta^{(t)}_1 = \eta_t(\bx^{(t)}_1) &\approx \sum_{\tau \in \calT_{\le t}} \omega^{|\partial_t\tau|}\prod_{v \in V(\tau) \setminus \partial_t\tau} A_{t_v,k_v}(\theta_v) \overline{\bZ}^\tau\,. \notag
\end{align}
Note that the two addends in the tree approximation for $\bx^{(t)}$ are an explicit decomposition of $\bx^{(t)}$ into its signal and noise components $\bx_1^{(t)} \overset{\text{Law}}{\approx} m_t \Theta + \sqrt{q_t}\,G$ for $(\Theta, G) \sim \pi \otimes\calN(0,1)$.

To see how the Hermite coefficient $A_{t,k}(\theta_v)$ arises,
suppose inductively that \eqref{eq:tree-apx-Z} holds for $\bx^{(t)}$. Expand the $t$th denoiser as in \eqref{eq:hermite-z},
\begin{align}
    \eta_t(\bx^{(t)}_1) \approx \eta_t(m_t\theta_1 + g^{(t)}) &= \sum_{k = 0}^\infty A_{t,k}(\theta_1) \frac{\He_k(g^{(t)}; q_t)}{k!}\,, \numberthis\label{eq:treeapx2}\\
    g^{(t)} &:= \sum_{\tau \in \calT_{\le t}: \deg_\tau(1) = 1}\omega^{|\partial_t\tau|}\prod_{v \in V(\tau) \setminus (\{1\}\cup\partial_t\tau)} A_{t_v, k_v}(\theta_v) \overline{\bZ}^\tau\,.
\end{align}
This is a scalar identity.
The key combinatorial property of \eqref{eq:treeapx2} is that the tree approximation of $\He_k(g^{(t)}; q_t)$ is:
\begin{equation}\label{eq:hermite}
    \frac{\He_k(g^{(t)}; q_t)}{k!} \approx (g^{(t)})^{\sqcup k}
\end{equation}
where $\sqcup k$ denotes restricting to collections of $k$ labeled trees from $g^{(t)}$ which are \emph{disjoint} except for their shared root vertex.
\eqref{eq:hermite} holds because (1)~the only non-negligible tree terms in $(g^{(t)})^k$ are those where the branches either fully overlap in pairs, or are disjoint; (2)~the non-leading terms in $\He_k(\cdot; q_t)$
cancel all terms with at least one such overlapping pair; (3)~the factor $k!$ removes the ordering in the product $(g^{(t)})^k$. Combinatorial observation (2) uses that $q_t$ is equal to the limiting variance of $g^{(t)}$, namely
\begin{align*}
    \E[(g^{(t)})^2] &= \sum_{\tau \in \calT_{\le t}: \deg_\tau(1) = 1} n^{-|E(\tau)|}\omega^{2|\partial_t\tau|}\prod_{v \in V(\tau) \setminus (\{1\}\cup\partial_t\tau)} \E[A_{t_v, k_v}(\theta_v)^2] \\
    &\approx \sum_{\tau \in \calT_{\le t-1}} n^{-|E(\tau)|}\omega^{2|\partial_{t-1}\tau|} \prod_{v \in V(\tau) \setminus \partial_{t-1}\tau} \E[A_{t_v, k_v}(\theta_v)^2]\\
    &= \sum_{\tau \in \calT_{\le t-1}} n^{-|E(\tau)|}\E[\widetilde A_{t-1,\tau}(\btheta)^2]\overset{\text{\cref{lem:idealized-tree-variance}}}{\approx} q_t\,.
\end{align*}
See \cite[Section 3]{jones2025fourier} for examples.

\begin{lemma}
\label{lem:idealized-tree-variance}
    Let $\sT_{\le t}$ be the set of unlabeled rooted trees with ordered children of depth at most $t$ and any size.
    Let $A_{t,k}(\Theta)$ be as in \eqref{eq:hermite-z} and $\widetilde A_{t,\tau}(\mathbf{\Theta}) = \omega^{|\partial_t\tau|}\prod_{v \in V(\tau)\setminus\partial_t\tau} A_{t_v,k_v}(\Theta_v)$
    with $\Theta_v \overset{\textnormal{i.i.d.}}{\sim} \pi$ independently for each $v$.
    For $t \geq 0$,
    \begin{equation*}\label{eq:qt-z-diagram}
        q_{t+1} = \sum_{\tau \in \sT_{\le t}} \frac{1}{\prod_{v\in V(\tau)}k_v!}\,\E[\widetilde A_{t,\tau}(\mathbf{\Theta})^2]\,. \end{equation*}
\end{lemma}
\begin{proof}
    For $t=0$, both sides equal $\omega^2=q_1$. For $t \geq 1$,
    we sum over the degree of the root vertex:
    \begin{align*}
        &\sum_{\tau \in \sT_{\le t}} \frac{1}{\prod_{v\in V(\tau)}k_v!}\,\E[\widetilde A_{t,\tau}(\mathbf{\Theta})^2] \\
        &\quad = \sum_{k = 0}^\infty \frac{1}{k!}\,\E_{\Theta \sim \pi}[A_{t,k}(\Theta)^2] \left(\sum_{\tau \in \sT_{\le t-1}} \frac{1}{\prod_{v\in V(\tau)}k_v!}\,\E[\widetilde A_{t-1,\tau}(\mathbf{\Theta})^2]\right)^{k} \\
        &\quad = \sum_{k = 0}^\infty \frac{1}{k!}\,\E [A_{t,k}(\Theta)^2] q_{t}^k
    \end{align*}
    Recall that $A_{t,k}(\Theta)$ is defined by the Hermite expansion of $\eta_t$ in \eqref{eq:hermite-z}.
    By the Parseval formula applied conditionally on $\Theta$,
    \[
        \sum_{k = 0}^\infty \frac{1}{k!}\,\E [A_{t,k}(\Theta)^2] q_{t}^k = \E[\eta_{t}(m_{t} \Theta + \sqrt{q_{t}} G)^2] = q_{t+1}. \qedhere
    \]
\end{proof}

To connect the calculated $Z$ basis coefficients in \eqref{eq:basis-approximations} with \cref{sec:stein-via-amp}, specialize to Bayes AMP with $\omega=\mu\ne0$. Then $m_t=q_t$ and $A_{t,k}(y)=A_k(y;q_t)$. For any fixed rooted tree $\tau$, once $t>\depth(\tau)$, the boundary $\partial_t\tau$ is empty, so
\[
    \widetilde A_{t,\tau}(\btheta)
    = \prod_{v\in V(\tau)}
      A_{k_v}(\theta_v;q_{t_v})
    \underset{t \to \infty}{\longrightarrow}
    \prod_{v\in V(\tau)} A_{k_v}(\theta_v;q_*)\,.
\]
Accounting for the normalization difference $\overline{\bZ}^{\tau}
=n^{-|\tau|/2}\bZ^\tau$, and since $h_\tau(\bZ) = \bZ^\tau$ for a tree, this heuristic recovers the Hermite coefficients in \eqref{eq:amp-hermite}.

\subsection{Heuristic \texorpdfstring{$Y$}{Y} basis coefficients}
\label{sec:Y-basis}

Now we heuristically calculate the coefficients $\widetilde B_{t,\tau}$ in \eqref{eq:coefficient-answer} for the $\bY$ basis.
These coefficients are related to the $\bZ$ basis coefficients by $\widetilde B_{t,\tau} = \widetilde A_{t,\tau}(\mathbf{0})$. This can be observed to hold in the case where $\bx^{(t)}$ is a polynomial function of the matrix $\bY$.
Write $\bx^{(t)}$ in the monomial basis in the $\bY$ variables (recall that $\calI$ indexes all possible monomials),
\[
    \bx^{(t)}_1 = \sum_{\alpha \in \calI} \widetilde B_{t,\alpha} \overline{\bY}^\alpha = \sum_{\alpha \in \calI} \widetilde B_{t,\alpha} \left(\frac 1n \btheta\btheta^\T + \overline{\bZ}\right)^\alpha\,.
\]
We see that when the RHS is expanded into the monomial basis in the $(\btheta, \bZ)$ variables, the coefficient on $\overline{\bZ}^\alpha$ is $\widetilde B_{t,\alpha}$.
On the other hand, if $\bx^{(t)}_1 = \sum_{\alpha \in \calI}\widetilde A_\alpha(\btheta)\overline{\bZ}^\alpha$, then this coefficient can also be obtained as $\widetilde A_\alpha(\mathbf{0})$. Hence $\widetilde B_{t,\tau} = \widetilde A_{t,\tau}(\mathbf{0})$.

The $\bY$ basis coefficients can also be derived directly by tracking the coefficients on the trees through the AMP iteration, like in the $\bZ$ basis.
Maintain ``tree approximations'' $\hat \btheta^{(t)} \approx \sum_{\tau \in \calT} \widetilde B_{t,\tau} \overline{\bY}^\tau$ and $\bx^{(t)} \approx \sum_{\tau \in \calT} \widetilde B'_{t,\tau}\overline{\bY}^\tau$.
The matrix multiplication $\overline{\bY}\hat \btheta^{(t)}$, after Onsager correction, has the effect
of extending the root by a single edge for all trees in $\hat \btheta^{(t)}$.
Applying $\eta_t$ joins any number $k$ of branches at the root with a coefficient $B_{t,k}$ to be derived momentarily.
The resulting tree approximations are
\begin{align}\label{eq:tree-apx-Y}
    \bx^{(t)}_1 &\approx \sum_{\tau \in \calT_{\le t}:\deg_\tau(1) = 1}\omega^{|\partial_t\tau|}\prod_{v \in V(\tau) \setminus (\{1\}\cup\partial_t\tau)} B_{t_v,k_v} \overline{\bY}^\tau\\
    \hat \theta^{(t)}_1 = \eta_t(\bx^{(t)}_1) &\approx \sum_{\tau \in \calT_{\le t}} \omega^{|\partial_t\tau|}\prod_{v \in V(\tau)\setminus\partial_t\tau} B_{t_v,k_v} \overline{\bY}^\tau\,. \notag
\end{align}
This recovers the predicted coefficients \eqref{eq:coefficient-answer} for $\hat \btheta^{(t)}$.
The coefficients $B_{t,k}$ were chosen by \eqref{eq:hermite-y},
\[
    \eta_t(\bx^{(t)}_1) = \sum_{k=0}^\infty B_{t,k}\frac{\He_k(\bx^{(t)}_1; q_t)}{k!}\,.
\]
Mirroring \cref{eq:hermite}, we claim the key combinatorial relation:
\begin{equation}\label{eq:y-hermite}
    \frac{\He_k(\bx^{(t)}_1; q_t)}{k!} \approx (\bx^{(t)}_1)^{\sqcup k}\,.
\end{equation}
As in the $Z$ basis, the tree terms in $(\bx^{(t)}_1)^k$ are obtained by overlapping pairs of subtrees,
and any three-way or higher intersection between subtrees is negligible.
Any such overlapping pairs are cancelled by lower-order terms in $\He_k(\cdot; q_t)$.
Analogously to \cref{lem:idealized-tree-variance}, proper cancellation of the paired terms requires that $q_t$ equals the diagrammatic variance of the trees in the $\bY$ basis; this is proven as the second formula in \cref{lem:at-fourier}.

To connect the calculated $Y$ basis coefficients in \eqref{eq:basis-approximations}
with the construction in
\cref{sec:nonzero-mean}, take $\omega=\mu\ne0$ and
$\eta_t=\hat\eta_t$ as in \eqref{eq:zero-a}.
Then $m_t=q_t=a_t$. The normalization of the Hermite coefficients in
\eqref{eq:hermite-y} and \eqref{eq:weights} differs by a factorial term,
$B_{t,k}=k!\,b_{t,k}$.
The normalizations differ because \cref{sec:upper-bound}
uses a sum over injective embeddings of unlabeled ordered trees,
whereas here we are using labeled trees on $[n]$;
each labeled tree appears $\prod_{v\in V(\tau)}k_v!$ times in \cref{sec:upper-bound}
which converts $b_{t,k}$ to $B_{t,k}$ at each vertex.
Each boundary leaf contributes a factor of $\mu$,
and multiplying the weights of two disjoint rooted
trees recovers the forest weights in
\eqref{eq:nonzero-estimator}.

\subsection{Heuristic \texorpdfstring{$Z$}{Z} basis coefficients via cavity method}
\label{sec:amp-triangle}

The $Z$ basis coefficients can also be heuristically calculated using a cavity method-style argument \cite{kz-notes}. Besides giving a second derivation of the tree coefficients in \eqref{eq:coefficient-answer},
this method also predicts the coefficients on cyclic graphs $\alpha$.
We demonstrate the method by heuristically computing $\widetilde A_{t,\triangle}(\btheta)$,
the coefficient on the triangle $\overline{Z}_{12}\overline{Z}_{23}\overline{Z}_{13}$.

Fix $t,n$ and suppose we have an $L^2(\btheta, \bZ)$ expansion,
\[
    \hat \btheta^{(t)}_1 = \sum_{\alpha \in \calI} \frac{\widetilde A_{t,\alpha}(\btheta)}{\alpha!} \He_\alpha(\overline{\bZ}; \tfrac 1n ).
\]
By orthogonality then Gaussian integration by parts, the triangle coefficient is
\[
    \widetilde A_{t,\triangle}(\btheta) = n^3\E_{\bZ}[\hat\btheta^{(t)}_1 \overline{\bZ}_{12}\overline{\bZ}_{23}\overline{\bZ}_{13}] = \E_{\bZ}[\partial_{\overline{Z}_{12}\overline{Z}_{23}\overline{Z}_{13}}\hat \btheta^{(t)}_1]\,.
\]
We calculate the mixed derivatives of $\hat\btheta^{(t)}$ and $\bx^{(t)}$ recursively.
Recall $\hat\btheta^{(t)} = \eta_t(\bx^{(t)})$. By the chain rule and product rule,
\begin{align*}\label{eq:chain-rule}
    \partial_{\overline{Z}_{12}\overline{Z}_{23}\overline{Z}_{13}}\hat \btheta^{(t)}_1 &= \eta_t'(\bx^{(t)}_1) \cdot \partial_{\overline{Z}_{12}\overline{Z}_{23}\overline{Z}_{13}} \bx^{(t)}_1 \\
    & + \eta_t''(\bx^{(t)}_1) \cdot\left((\partial_{\overline{Z}_{12}\overline{Z}_{23}} \bx^{(t)}_1)\cdot (\partial_{\overline{Z}_{13}} \bx^{(t)}_1) + (\partial_{\overline{Z}_{12}\overline{Z}_{13}} \bx^{(t)}_1)\cdot (\partial_{\overline{Z}_{23}} \bx^{(t)}_1) + (\partial_{\overline{Z}_{23}\overline{Z}_{13}} \bx^{(t)}_1)\cdot (\partial_{\overline{Z}_{12}} \bx^{(t)}_1)\right) \\
    &+ \eta_t'''(\bx^{(t)}_1) (\partial_{\overline{Z}_{12}} \bx^{(t)}_1) (\partial_{\overline{Z}_{23}} \bx^{(t)}_1) (\partial_{\overline{Z}_{13}} \bx^{(t)}_1) \numberthis
\end{align*}
The cavity method assumes key simplifications in the limit $n \to \infty$ for the quantity $\bx^{(t)}_1$ which is called the ``local field'' or ``cavity field'' at vertex 1 and time $t$.
In the expression for $\bx^{(t)}_1$,
\begin{align}\label{eq:cavity-x-1}
    \bx^{(t)} = \overline{\bY} \hat \btheta^{(t-1)} - (\text{Onsager correction}) = \frac{\langle\btheta, \hat \btheta^{(t-1)}\rangle}{n} \btheta + \overline{\bZ} \hat \btheta^{(t-1)} - (\text{Onsager correction})
\end{align}
we assume that the dependence of $\bx^{(t)}_1$ on terms such as $\overline{\bZ}_{23}$ which do not include
vertex 1 is asymptotically negligible. Under this approximation, \cref{eq:chain-rule}
simplifies to
\begin{align*}
    \partial_{\overline{Z}_{12}\overline{Z}_{23}\overline{Z}_{13}} \hat \btheta^{(t)}_1 &\approx \eta_t'(\bx^{(t)}_1) \cdot \partial_{\overline{Z}_{12}\overline{Z}_{23}\overline{Z}_{13}} \bx^{(t)}_1 \\
    & + \eta_t''(\bx^{(t)}_1) \cdot\left((\partial_{\overline{Z}_{12}\overline{Z}_{23}} \bx^{(t)}_1)\cdot (\partial_{\overline{Z}_{13}} \bx^{(t)}_1) + (\partial_{\overline{Z}_{23}\overline{Z}_{13}} \bx^{(t)}_1)\cdot (\partial_{\overline{Z}_{12}} \bx^{(t)}_1)\right) \numberthis\label{eq:chain-rule2}
\end{align*}
\cref{eq:chain-rule2} can be interpreted as taking the two derivatives $\partial_{\overline{Z}_{12}}$ and $\partial_{\overline{Z}_{13}}$ on a single copy of $\bx^{(t)}_1$
or splitting the derivatives among two copies of $\bx^{(t)}_1$. In general, a derivative with order $k$ at vertex 1 can be split among up to $k$ copies of $\bx^{(t)}_1$ and the number of copies used equals the order of the $\eta$ derivative.

A multiple derivative of $\bx^{(t)}_1$ in \eqref{eq:cavity-x-1} is evaluated by the cavity assumption and the product rule,
\[
    \partial_{\overline{Z}_{12}\overline{Z}_{23}\overline{Z}_{13}} \bx^{(t)}_1 \approx \partial_{\overline{Z}_{12}\overline{Z}_{23}\overline{Z}_{13}}[ (\overline{\bZ} \hat\btheta^{(t-1)})_1 - (\text{Onsager correction})] \approx \partial_{\overline{Z}_{23}\overline{Z}_{13}}\hat\btheta^{(t-1)}_2 + \partial_{\overline{Z}_{12}\overline{Z}_{23}}\hat \btheta^{(t-1)}_3\,.
\]
That is, one derivative moves the root vertex from 1 to elsewhere, with the remaining derivatives to be evaluated recursively.
Iterating these derivative computations leads to, assuming $t \geq 4$,
\begin{align*}
    &\quad \partial_{\overline{Z}_{12}\overline{Z}_{23}\overline{Z}_{13}} \hat \btheta^{(t)}_1 \\
    &\approx \eta_t'(\bx^{(t)}_1) \eta'_{t-1}(\bx^{(t-1)}_2) \eta'_{t-2}(\bx^{(t-2)}_3) \cdot \eta_{t-3}(\bx^{(t-3)}_1) + \eta'_t(\bx^{(t)}_1) \eta'_{t-1}(\bx^{(t-1)}_3) \eta'_{t-2}(\bx^{(t-2)}_2) \eta_{t-3}(\bx^{(t-3)}_1) \\
    & + \eta_t''(\bx^{(t)}_1) \eta'_{t-1}(\bx^{(t-1)}_2) \eta_{t-2}(\bx^{(t-2)}_3) \eta_{t-1}(\bx^{(t-1)}_3) + \eta''_t(\bx^{(t)}_1) \eta'_{t-1}(\bx^{(t-1)}_3)\eta_{t-2}(\bx^{(t-2)}_2)\eta_{t-1}(\bx^{(t-1)}_2) \numberthis\label{eq:chain-rule3}
\end{align*}

Finally, the limiting expectation of \cref{eq:chain-rule3}
is obtained by averaging $\bx^{(t)}_i$ in their state evolution limit.
Let $m_t, Q_t$ be the state evolution parameters for the AMP iteration where $Q_t \in \R^{t \times t}$ is the limiting covariance matrix of the Gaussian fields. These parameters are defined recursively by,
assuming the initialization covers when $t,r$, or $s$ equals 0,
\begin{align*}
    m_{t+1} = \E[\Theta \eta_t(m_t \Theta + G_t)], &\qquad (Q_{t+1})_{r+1,s+1} = \E[\eta_{r}(m_r \Theta + G_r)\eta_{s}(m_{s}\Theta + G_{s})] , \\
    \Theta \sim \pi, &\qquad (G_1, \dots, G_t) \sim \calN(\mathbf{0}, Q_t)\,.
\end{align*}
Let $X^t_v = m_t \theta_v + G_{t,v}$ for each vertex $v$ and time $t$ with $(G_{1,v}, \dots, G_{t,v}) \sim \calN(0, Q_t)$ sampled independently for each vertex $v$.
Taking the expectation of \cref{eq:chain-rule3}
with these variables, the predicted coefficient on the triangle is
\begin{align}
    \widetilde A_{t,\triangle}(\btheta) \approx\;& \E_{X}[A_{\triangle}(\btheta, X)]\notag\\
    A_{\triangle}(\btheta, X) = \;&\eta'_t(X^t_1) \eta'_{t-1}(X^{t-1}_2)\eta'_{t-2}(X^{t-2}_3) \eta_{t-3}(X^{t-3}_1) + \eta'_t(X^t_1) \eta'_{t-1}(X^{t-1}_3) \eta'_{t-2}(X_2^{t-2})\eta_{t-3}(X_1^{t-3}) \notag\\
    +& \eta''_t(X_1^t)\eta'_{t-1}(X_2^{t-1})\eta_{t-2}(X_3^{t-2})\eta_{t-1}(X^{t-1}_3) + \eta''_{t}(X^t_1)\eta'_{t-1}(X_3^{t-1})\eta_{t-2}(X_2^{t-2})\eta_{t-1}(X_2^{t-1})\,.
    \label{eq:triangle-amp-finite-time}
\end{align}

\section{Comparison of candidate dual constructions}
\label{app:cyclic-terms}

\newcommand{\MP}{\textnormal{MP}}
\newcommand{\BP}{\textnormal{BP}}

In this section, we compare the coefficients on cyclic terms for our history certificate with the certificate by~\cite{li2026almost}, and also with the coefficients of AMP.
While the coefficients on the tree terms are the same in all three settings,
the coefficients on the cyclic terms are different,
and we demonstrate the differences below through the example of the triangle.

We first describe the approach by~\cite{li2026almost}, which may be viewed as an extension of the cumulant method of~\cite{SW-22} by
conditioning on an auxiliary observation $\bR$. We describe it for the spiked Wigner model with a general estimand
$x=x(\btheta)$. Recall $\MMSE_{\leq D}(x) := \inf_{f \in \R[\bY]_{\leq D}} \E[(f(\bY)-x)^2]$.  Write $S_{ij}:=\theta_i\theta_j$, so that
$\bY=\bS/\sqrt n+\bZ$.
Let $\bR$ be any auxiliary observation such that $\bZ$ is
independent of $(\btheta,\bR)$, and define
\[
 \kappa_\alpha(\bR)
 :=\kappa\bigl(x,(S_e)_{e\in E(\alpha)}\mid\bR\bigr),
 \qquad
 \kappa_\emptyset(\bR):=\E[x\mid\bR],
\]
where each edge appears with its multiplicity in $\alpha$, and $\kappa(\cdot)$ in the right-hand side denotes the joint conditional cumulant (see \eqref{eq:def:cumulant}). Then, \cite[Lemma 2.2]{li2026almost} gives
\[
\MMSE_{\leq D}(x) \geq \E[x^2]-\sum_{\alpha:|\alpha|\leq D} \frac{\E[c_{\alpha}^{\circ}(\btheta)^2]}{\alpha!n^{|\alpha|} }\,,\quad\textnormal{where}\quad c_{\alpha}^{\circ}(\btheta)= \E [\kappa_{\alpha}(\bR)\mid \btheta].
\]
Taking $\bR$ to be constant recovers the cumulant bound
of~\cite[Theorem~2.2]{SW-22}, in which case
$c_\alpha^\circ(\btheta)
=\kappa(x,(S_e)_{e\in E(\alpha)})$ is deterministic. We remark that Even, Giraud, and Verzelen~\cite[Appendix~E.1]{even2025computational}
previously used conditional cumulants in a different way: their bound involves
$\E[\kappa_\alpha(\bR)^2]$, rather than
$\E[\E[\kappa_\alpha(\bR)\mid\btheta]^2]$ as in the bound above.

In our notation, $c_\alpha(\btheta)=c_\alpha^\circ(\btheta)$ is a valid dual certificate for use in \cref{cor:MMSE:lower:bound}. In fact, the construction by \cite{li2026almost} satisfies a stronger system of constraints that hold pointwise over $\bR$: the $\kappa_\alpha(\bR)$ are the unique solution to
\[ \sum_{\beta\leq\alpha}\binom{\alpha}{\beta}
\kappa_\beta(\bR)\,
\E[\bS^{\alpha-\beta}\mid\bR]
=
\E[x\bS^\alpha\mid\bR]
\qquad\text{almost surely}, \]
which recovers our dual constraints~\eqref{eq:dual:constraints:hermite:functional} for $c^\circ$ by taking expectation on both sides. Indeed, since $\bS$ is measurable with respect to $\btheta$, the tower
property gives
\[
\begin{aligned}
\E\left[\kappa_\beta(\bR)
          \E[\bS^{\alpha-\beta}\mid\bR]\right]=\E[\kappa_\beta(\bR)\bS^{\alpha-\beta}]=\E\left[\E[\kappa_\beta(\bR)\mid\btheta]
            \bS^{\alpha-\beta}\right]=\E[c_\beta^\circ(\btheta)\bS^{\alpha-\beta}].
\end{aligned}
\]
We do not know whether the approach by \cite{li2026almost} can be adapted to show degree-$o(n)$ hardness as we have done under the stability condition $\Psi''(q_*;\pi)>0$. The history certificate seems to make the degree-$o(n)$
analysis more tractable, where each coefficient is a sum over histories
of products of the scalar functions $A_r(\theta_v;q)$ and powers
of~$\theta_v$. This representation allows us to use the stability condition to prove the analytic estimates on $A_r$
in \cref{lem:analytic-1}, which we combine with graph reductions and counting in \cref{subsec:lower:positive-excess} to bound the certificate's norm.

For the spiked Wigner model, \cite{li2026almost} chooses
$\bR=\sqrt q\,\btheta+\bG$, where $\bG\sim\calN(0,I_n)$
is independent of $(\btheta,\bZ)$ and $q=q_*$ is the AMP
fixed point. The resulting $c_{\alpha}^\circ(\btheta)$ agrees with our
history certificate on rooted trees (see \cref{subsec:lower:excess-zero}), but their coefficients differ on cyclic graphs in general. We illustrate this distinction
for the target $x=\theta_1$ and the triangle
$\triangle=\{12,23,13\}$ rooted at $1$, also comparing them
with the AMP prediction from \cref{sec:amp-triangle}.

Write $X_i=\sqrt qR_i=q\theta_i+\sqrt q G_i$, and let
$m_i,v_i,k_i$ denote the posterior mean, variance, and third
cumulant of $\theta_i$ given $X_i$. For $q>0$, these satisfy
\begin{equation}\label{eq:def:posterior:cumulants}
m_i=\eta(X_i;q),\qquad
 v_i=\eta'(X_i;q),\qquad
 k_i=\eta''(X_i;q).
\end{equation}
Let $(c^{\mathrm{hist}}_{\alpha}(\btheta))_{|\alpha|\leq D}$ be the history coefficients from \cref{subsec:outline:history:certificate} with given $q$, and let
$c_\triangle^{\AMP}(\btheta)$ denote the triangle coefficient
predicted by the
AMP expansion in \cref{sec:amp-triangle}. The two certificate
constructions and this AMP prediction give
\begin{align}\label{eq:three:coefficients}
 c^{\mathrm{hist}}_\triangle(\btheta)
   &=\E_G\big[ 2\theta_1v_1v_2v_3
     +k_1(\theta_2m_2v_3+v_2\theta_3m_3+v_2v_3)\big], & (\text{History certificate, Def. \ref{def:history-cert}})\\
 c^{\circ}_\triangle(\btheta)
   &=\E_G\big[2m_1v_1v_2v_3
     +k_1(m_2^2v_3+v_2m_3^2+v_2v_3)\big], & (\text{Li's certificate~\cite{li2026almost}})\\
 c^{\AMP}_\triangle(\btheta)
   &=\E_G\big[2m_1v_1v_2v_3
     +k_1(m_2^2v_3+v_2m_3^2)\big], & (\text{Asymptotic AMP coefficient})
\end{align}
where $\E_G$ averages over $\bG$ with $\btheta$ fixed.

One can observe that the difference between $c_\triangle^{\text{hist}}(\btheta)$ and
$c_\triangle^{\circ}(\btheta)$ is that $c_\triangle^{\circ}(\btheta)$ employs an extra layer of averaging, replacing $\theta_i$ by $m_i$. The history coefficient was computed in
\cref{subsec:outline:history:certificate}, while the AMP
coefficient is obtained from the heuristic calculation in \cref{sec:amp-triangle}.\footnote{
$c_\triangle^\AMP(\btheta)$ is also the leading-order coefficient on the triangle for the unapproximated message passing algorithm, of which AMP is an approximation \cite[Section 16.4]{kz-notes}.
This is also shown by the cavity method in \cref{sec:amp-triangle}.}
We derive the formula
for $c_\triangle^\circ(\btheta)$ in the proof of the following
lemma, which shows that the two certificates give different triangle coefficients in general.

\begin{lemma}
\label{lem:distinct-triangle-coefficients}
Let $n\ge3$, $p\in (0,1)$, and
$\theta_1,\ldots,\theta_n\stackrel{\mathrm{iid}}{\sim}
\operatorname{Ber}(p)$. Let
$R_i=\sqrt q\,\theta_i+G_i$, where
$G_i\stackrel{\mathrm{iid}}{\sim}\cN(0,1)$ are independent
of $\btheta$ and fix $q \geq 0$.
For the estimand $x=\theta_1$, let
$c_\alpha^{\circ}(\btheta;q)$ denote the
conditional-cumulant coefficient defined above with this
choice of $\bR$, 
let $\triangle=\{12,23,31\}$ be rooted at $1$, and let
$c_\triangle^{\mathrm{hist}}(\btheta;q)$ be its history
coefficient from~\eqref{eq:history:coefficient}. Then, we have the following.
\begin{enumerate}[label=(\alph*)]
    \item With probability at least $\min\{p,1-p\}^3$, $c_\triangle^{\mathrm{hist}}(\btheta;q)\neq c_\triangle^{\circ}(\btheta;q)$.
    \item \label{item:dual-fail} Suppose $p\neq1/2$. Consider the coefficients
$(c_\beta^{\AMP}(\btheta))_{\beta\le\triangle}$, where
$c_\triangle^{\AMP}$ is the AMP
coefficient from \cref{sec:amp-triangle} and
$c_\beta^{\AMP}(\btheta)=c_\beta^{\mathrm{hist}}(\btheta)$
whenever $\beta$ is a tree. Then this family does not satisfy the dual
constraint~\eqref{eq:dual:constraints:hermite:functional}
for $\alpha=\triangle$.
\end{enumerate}
\end{lemma}

\begin{proof}
We first compute $\kappa_\triangle(\bR)$ for a general product
prior, using the first three posterior cumulants $m_i,v_i,k_i$
defined in~\eqref{eq:def:posterior:cumulants}. Conditional on $\bR$, the coordinates $\theta_i$ are independent.
Substituting
$\E[\theta_i^2\mid\bR]=m_i^2+v_i$ and
$\E[\theta_i^3\mid\bR]=m_i^3+3m_iv_i+k_i$
into the definition of the joint cumulant~\eqref{eq:def:cumulant} gives
\[
\begin{aligned}
\kappa_\triangle(\bR)
&=\kappa(\theta_1,\theta_1\theta_2,\theta_2\theta_3,
                   \theta_3\theta_1\mid\bR)\\
&=k_1\big[(m_2^2+v_2)(m_3^2+v_3)-m_2^2m_3^2\big]
  +2m_1v_1v_2v_3\\
&=k_1(v_2v_3+v_2m_3^2+m_2^2v_3)+2m_1v_1v_2v_3.
\end{aligned}
\]
Averaging over $\bG$ proves the formula
for $c_\triangle^\circ$ displayed above.

We now specialize to $\pi=\operatorname{Ber}(p)$. Write $\ell_i= \log (\frac{p}{1-p})+X_i-q/2$, where $X_i=q\theta_i+\sqrt{q}G_i$. By a direct calculation,
\[
m_i=\frac{e^{\ell_i}}{1+e^{\ell_i}},\qquad
 v_i=m_i(1-m_i),\qquad k_i=v_i(1-2m_i).
\]
For part~(a), write $A_{r,i}:=A_r(\theta_i;q)$.
Averaging the identities $m_i^2=m_i-v_i$ and $2m_iv_i=v_i-k_i$ gives
\[
\E_G[m_i^2]=A_{0,i}-A_{1,i},\qquad
\E_G[m_iv_i]=\frac{A_{1,i}-A_{2,i}}2.
\]
Independence of the $G_i$ therefore yields
\[
c_\triangle^\circ
=A_{1,1}A_{1,2}A_{1,3}
+A_{2,1}(A_{0,2}A_{1,3}+A_{1,2}A_{0,3}-2A_{1,2}A_{1,3}).
\]
Set
\[
a_0:=A_1(0;q),\qquad a_1:=A_1(1;q),\qquad
 b:=A_0(0;q),\qquad d:=A_2(0;q).
\]
Since $0<m_i<1$, we have $a_0,a_1,b>0$.
Subtracting from the history coefficient computed in
\cref{subsec:outline:history:certificate}, the difference
$\Delta:=c_\triangle^{\mathrm{hist}}-c_\triangle^\circ$ satisfies
\[
\Delta(0,1,1)=a_1^2(3d-a_0),\qquad
\Delta(0,0,1)=a_1\big[a_0(3d-a_0)-bd\big].
\]
If the first value vanishes, then $d=a_0/3$, making the second
value $-a_0a_1b/3<0$. Thus at least one is nonzero.
These configurations have probabilities $(1-p)p^2$ and $(1-p)^2p$,
respectively, each at least $\min\{p,1-p\}^3$. This proves~(a).

For part~(b), suppose for a contradiction that $(c_{\beta}^{\AMP}(\btheta))_{\beta\leq \triangle}$ satisfies the dual constraint~\eqref{eq:dual:constraints:hermite:functional} for $\alpha=\triangle$.
Since $(c_{\beta}^{\circ}(\btheta))_{\beta\leq \triangle}$ also satisfy the dual constraints,
and $c_{\beta}^\circ(\btheta)=c_{\beta}^{\AMP}(\btheta)$ for all $\beta \lneq \triangle$, we have from~\eqref{eq:three:coefficients} that
\[
0=\E[c_\triangle^\circ(\btheta)-c_\triangle^{\AMP}(\btheta)]
=\E[k_1v_2v_3]
=\E[k_1](\E[v_1])^2.
\]
Since $v_1>0$ almost surely, this implies $\E[k_1]=0$. We show that this is impossible when $p\neq1/2$.
At $q=0$, we have $\E[k_1]=p(1-p)(1-2p)\neq0$.
For $q>0$, let $f_q(x)$ be the density of $\ell_1\equiv \log(\frac{p}{1-p})+q\theta_1+\sqrt{q}G_1-q/2$. Then, by a direct calculation of this Gaussian mixture law,
\[
f_q(x)\propto
\cosh(x/2)\exp\left\{-\frac{(x-\log(p/(1-p)))^2}{2q}\right\}.
\]
In particular, we have $\frac{f_q(x)}{f_q(-x)}
=\left(\frac{p}{1-p}\right)^{2x/q}$. Write $k_1=K(\ell_1)$, where
$K(x)=\frac{e^x(1-e^x)}{(1+e^x)^3}$ is odd and strictly negative
for $x>0$. Consequently,
\[
\E[k_1]=\int_0^\infty K(x)\bigl[f_q(x)-f_q(-x)\bigr]dx.
\]
If $p>1/2$, then $f_q(x)>f_q(-x)$ for every $x>0$,
so $\E[k_1]<0$. Similarly, if $p<1/2$, then $\E[k_1]>0$.
In either case $\E[k_1]\neq0$, which is a contradiction, concluding the proof of~(b).
\end{proof}

\section{Omitted proofs from Section~\ref{sec:dual-analysis}}
\label{subsec:omitted:lower:proofs}

\begin{proof}[Proof of Lemma~\ref{lem:Gaussian:field}]
We prove the following, which will be useful later: for all integers $s, r \geq 0$,
\begin{equation}\label{eq:cumulant:Hermtie:coefficient}
    \E[\kappa_{s+1}(\Theta \mid \sF_q)h_r(G) \mid \Theta]=\frac{q^{r/2}}{\sqrt{r!}}A_{r+s}(\Theta; q)\,,
\end{equation}
where the RHS is understood as $A_{s}(\Theta; q)$ when $r=0$. Thus, using that $h_r \equiv \He_r/\sqrt{r!}$, ~\eqref{eq:cumulant:Hermtie:coefficient} with $r=0$ implies the conclusion of Lemma~\ref{lem:Gaussian:field}. Note that, when $q=0$, the conditioning field is trivial and both sides of~\eqref{eq:cumulant:Hermtie:coefficient} equal $\kappa_{s+1}(\Theta)\,\In_{\{r=0\}}$. Henceforth, we assume that $q>0$. By Bayes' rule, we have for all $y, t \in \R$, 
\[
    \log \E[e^{t\Theta} \mid q\Theta+\sqrt{q}G=y]=\log \E[e^{(y+t)\Theta-q\Theta^2/2}]-\log \E[e^{y\Theta-q\Theta^2/2}]\,.
\]
On the other hand, 
\[
    \eta(y; q)=\frac{\dx}{\dx{y}} \log \E[e^{y\Theta-q\Theta^2/2}]\,.
\]
Combining the preceding two displays and differentiating $s+1$ times in $t$ gives
\begin{equation}\label{eq:log:generating:Gaussian:field}
    \left. \frac{\dx^{s+1}}{\dx{t^{s+1}}} \Big\{ \log \E[e^{t\Theta} \mid q\Theta+\sqrt{q}G=y] \Big\} \right|_{t=0}=\partial_y^s \eta(y; q)\,.
\end{equation}
Using the standard log moment-generating function representation of cumulants (see e.g.~\cite[Section~2.3.4, Eq.~(2.9)]{mccullagh2018tensor}), the LHS of~\eqref{eq:log:generating:Gaussian:field} is precisely $\kappa_{s+1}(\Theta \mid q\Theta+\sqrt{q}G=y)$. Recalling $\sF_q=\sigma(q\Theta+\sqrt{q}G)$, the preceding discussion yields
\[
    \kappa_{s+1}(\Theta \mid \sF_q)=\partial_y^s \eta(q\Theta+\sqrt{q}G; q)\,.
\]
Taking expectation over the Gaussian randomness only, we obtain
\[
    \E[\kappa_{s+1}(\Theta \mid \sF_q)h_r(G) \mid \Theta]=\E[\partial_y^s \eta(q\Theta+\sqrt{q}G; q)h_r(G) \mid \Theta]=\frac{q^{r/2}}{\sqrt{r!}}A_{r+s}(\Theta; q)\,,
\]
where the second equality used the Gaussian integration by parts formula for Hermite polynomials and the definition of $A_{r+s}(\cdot; q)$. This concludes the proof.
\end{proof}

For $q \geq 0$, recall the definition of the power series
\begin{equation}\label{eq:Gamma:q}
    \Gamma_q(y) := \sum_{r \geq 0} \E[A_r(\Theta; q)^2]\, \frac{y^r}{r!}\,.
\end{equation}
\begin{lemma}
\label{lem:properties:of:Psi:Phi}
For all $q \geq 0$, it holds that
\[
    \Gamma_q(q)=\Phi(q)\, \quad \text{and}\, \quad \Gamma_q'(q)=\Phi'(q)\,.
\]
At $q=0$, we interpret $\Phi'(0)$ as the right derivative and use the convention $\Gamma_0'(0)=\E[A_1(\Theta;0)^2]$.
\end{lemma}

\begin{proof}
For $q=0$, we can compute directly $\Gamma_0(0)=\E[\Theta]^2=\Phi(0)$ and $\Gamma_0'(0)=\Var(\Theta)^2=\Phi'(0)$. So assume that $q>0$. Using~\eqref{eq:cumulant:Hermtie:coefficient}, 
\[
    \sum_{r=0}^\infty \E[A_{r+s}(\Theta;q)^2]\frac{q^r}{r!}=\E\left[ \sum_{r=0}^\infty \E[\kappa_{s+1}(\Theta \mid \sF_q)h_r(G) \mid \Theta]^2 \right] \stackrel{\text{Parseval}}{=} \E[\kappa_{s+1}(\Theta \mid \sF_q)^2]\,.
\]
Taking $s=0$ and equating both sides of the above display proves $\Gamma_q(q)=\Phi(q)$. Taking $s=1$ and recalling that the second cumulant is the variance, the preceding gives $\Gamma_q'(q)=\E[\Var(\Theta \mid \sF_q)^2]$, which, by~\cite[Proposition~9]{GWSV11}, is precisely $\Phi'(q)$.
\end{proof}

\begin{corollary}
\label{cor:perturbation:of:Psi:Phi}
Let $\rho=\frac{1}{2}(1+\Phi'(q_*))$. If $\Phi'(q_*)<1$, there exists $\eps_0 \equiv \eps_0(\pi)>0$ such that
\[
    \Phi(q_*+\eps) \leq q_*+\rho\eps\, \quad \text{and}\, \quad \sup_{y \in [0, q_*+\eps]} \Gamma_{q_*+\eps}'(y) \leq \rho\,, \quad \forall \eps \in [0, \eps_0)\,.
\]
\end{corollary}

\begin{proof}
The stability assumption $\Phi'(q_*)<1$ implies that $\Phi'(q_*)<\rho$. By continuity, there exists an $\eps_0 \equiv \eps_0(\pi)$ such that
\begin{equation}\label{eq:sup:Phi':eps}
    \sup_{s \in [0, \eps_0)} \Phi'(q_*+s) \leq \rho\,.
\end{equation}
Choose $\eps \in [0, \eps_0)$. Using $\Phi(q_*)=q_*$ and the fundamental theorem of calculus, 
\[
    \Phi(q_*+\eps)-q_*=\Phi(q_*+\eps)-\Phi(q_*)=\int_0^\eps \Phi'(q_*+s) \, \dx{s} \leq \rho\eps\,,
\]
which proves the first item. For the second item, observe that $\Gamma_q'$ is non-decreasing on $[0, q]$, which gives
\[
    \sup_{y \in [0, q_*+\eps]} \Gamma_{q_*+\eps}'(y)=\Gamma_{q_*+\eps}'(q_*+\eps) \stackrel{\text{Lemma~\ref{lem:properties:of:Psi:Phi}}}{=} \Phi'(q_*+\eps) \stackrel{\text{\eqref{eq:sup:Phi':eps}}}{\leq} \rho\,,
\]
as claimed. This concludes the proof.
\end{proof}

\begin{proof}[Proof of Lemma~\ref{lem:analytic-1}]
First assume the stability condition $\Phi'(q_*)<1$. Let $\Theta \sim \pi$ and $q=q_*+\eps$. We show that there exists a constant $C=C(\eps, \pi)>0$ such that, for any integers $p, s \geq 0$, 
\begin{equation}\label{eq:local:factor:estimate:A}
    M(p, s)=\sum_{r \geq 0} \frac{\Phi(q)^r}{r!} \E[|\Theta|^{2p}A_{r+s}(\Theta; q)^2] \leq C^{1+p+s}(p+s)!\,.
\end{equation}
Let $\delta=q-\Phi(q)$. By Corollary~\ref{cor:perturbation:of:Psi:Phi}, we have $\delta \geq (1-\rho)\eps$. The binomial theorem gives
\[
    \frac{(r+s)!}{s!r!} \Phi(q)^r \delta^s \leq (\Phi(q)+\delta)^{r+s}=q^{r+s}\,,
\]
and rearranging yields
\begin{equation}
    \frac{\Phi(q)^r}{r!} \leq \frac{s!q^{r+s}}{\delta^s(r+s)!}\,.
\end{equation}
Let $G \sim \cN(0, 1)$ be independent of $\Theta$. Inserting the above into~\eqref{eq:local:factor:estimate:A}, we obtain
\[
\begin{aligned}
    \sum_{r \geq 0} \frac{\Phi(q)^r}{r!} \E[|\Theta|^{2p}A_{r+s}(\Theta; q)^2]
    & \leq \frac{s!}{\delta^s} \E\left[ |\Theta|^{2p} \sum_{r \geq 0} \frac{q^{r+s}}{(r+s)!} A_{r+s}(\Theta; q)^2 \right]   \\
    & \leq \frac{s!}{\delta^s} \E\left[ |\Theta|^{2p} \sum_{r \geq 0} \frac{q^{r}}{r!} A_{r}(\Theta; q)^2 \right]           \\
    &= \frac{s!}{\delta^s} \E\left[ |\Theta|^{2p} \sum_{r \geq 0} \E[\eta(q\Theta+\sqrt{q}G; q)h_r(G) \mid \Theta]^2 \right ] \\
    &= \frac{s!}{\delta^s} \E\left[ |\Theta|^{2p} \E[\eta(q\Theta+\sqrt{q}G; q)^2 \mid \Theta] \right]\,,
\end{aligned}
\]
where the last equality used Parseval's identity in the Hermite basis. By the tower property and definition of the Bayes denoiser, 
\[
    \E\left[ |\Theta|^{2p} \E[\eta(q\Theta+\sqrt{q}G; q)^2 \mid \Theta] \right]=\E\left[ |\Theta|^{2p} \E[\Theta \mid q\Theta+\sqrt{q}G]^2 \right] \stackrel{\text{C-S}}{\leq} \sqrt{\E[|\Theta|^{4p}] \E[|\Theta|^4]}\,.
\]
The subgaussian assumption on $\pi$ implies that there is a constant $C=C(\pi) \geq 1$ such that, for any integer $m \geq 0$, $\E[|\Theta|^{4m}] \leq C^{2m}(2m)! \leq (2C)^{2m}m!^2$. Using $\delta \geq (1-\rho)\eps$, the preceding gives
\[
\begin{aligned}
    \sum_{r \geq 0} \frac{\Phi(q)^r}{r!} \E[|\Theta|^{2p}A_{r+s}(\Theta; q)^2] &\leq ((1-\rho)\eps)^{-s}(2C)^{2p+2}s!p! \\ &\leq 4C^2 \max\{((1-\rho)\eps)^{-1}, 4C^2\}^{s+p}(s+p)!\,.
\end{aligned}
\]
After possibly enlarging the constant $C=C(\eps, \pi)$, this proves~\eqref{eq:local:factor:estimate:A}.
Also, Corollary~\ref{cor:perturbation:of:Psi:Phi} gives $M(0,1)=\Gamma_q'(\Phi(q))\leq\rho$. This completes the proof of item (ii).

Below, we drop the stability condition, and allow $\Phi'(q_*)=1$. Take $q=q_* \geq 0$ and let $\rho=\frac{1}{2}(1+\Phi'(q_*))$. By direct calculation, 
\[
\begin{aligned}
    M(0,1)
    &=\sum_{r\geq 0}\frac{\Phi(q)^r}{r!}\E[A_{r+1}(\Theta;q)^2]
      =\Gamma_q'(\Phi(q))=\Gamma_{q_*}'(q_*)
      \stackrel{\text{Lemma~\ref{lem:properties:of:Psi:Phi}}}{=}
      \Phi'(q_*)\leq\rho\,.
\end{aligned}
\]
Finally, we establish the following general bound: there is a constant $C=C(\pi)>0$ such that, for all integers $p, s \geq 0$, 
\begin{equation}\label{eq:local:factor:estimate:C}
    M(p, s)=\sum_{r \geq 0} \frac{\Phi(q_*)^r}{r!} \E[|\Theta|^{2p}A_{r+s}(\Theta; q)^2] \leq C^{1+p+s}((p+s)!)^3\,.
\end{equation}
Recall $\sF_q=\sigma(q\Theta+\sqrt{q}G)$.
To prove~\eqref{eq:local:factor:estimate:C}, we use the cumulant representation~\eqref{eq:cumulant:k:copies} and the identity (see the proof of Lemma~\ref{lem:Gaussian:field})
\begin{equation}
    \E[\kappa_{s+1}(\Theta \mid \sF_q)h_r(G) \mid \Theta]=\frac{q^{r/2}}{\sqrt{r!}}A_{r+s}(\Theta; q)\,.
\end{equation}
This identity also holds at $q=0$, with the convention $0^0=1$: the field $\sF_0$ is trivial, and both sides equal $\kappa_{s+1}(\Theta)\,\In_{\{r=0\}}$, so the argument below applies for all $q\geq 0$.
Using $\Phi(q)=q$ and Parseval's identity, the LHS of~\cref{eq:local:factor:estimate:C} is bounded by
\[
\begin{aligned}
    \sum_{r \geq 0} \frac{q^r}{r!} \E[|\Theta|^{2p}A_{r+s}(\Theta; q)^2] &=\E\left[|\Theta|^{2p} \sum_{r \geq 0} \E[\kappa_{s+1}(\Theta \mid \sF_q)h_r(G) \mid \Theta]^2 \right] \stackrel{\text{Parseval}}{=} \E[ |\Theta|^{2p} \kappa_{s+1}(\Theta \mid \sF_q)^2]\,.
\end{aligned}
\]
Hölder's inequality implies the pointwise bound
\[
\begin{aligned}
    |\kappa_{s+1}(\Theta \mid \sF_q)| \leq \E[|\Theta|^{s+1} \mid \sF_q] \sum_{\sigma \in \cP([s+1])} (|\sigma|-1)!\,.
\end{aligned}
\]
The sum over partitions on the RHS is bounded by
\[
    \sum_{\sigma \in \cP([s+1])} |\sigma|! \leq \sum_{\sigma \in \cP([s+1])} |\sigma|! \prod_{I \in \sigma} |I|!=(s+1)!2^s \leq 4^ss!\,.
\]
Putting everything together, the LHS of~\cref{eq:local:factor:estimate:C} is at most
\[
    C^s(s!)^2 \E[|\Theta|^{2p} \E[|\Theta|^{s+1} \mid \sF_q]^2] \stackrel{\text{Hölder}}{\leq} C^s(s!)^2 \E[|\Theta|^{2(1+p+s)}] \leq C^{1+p+s}(s!)^2(1+p+s)!\,.
\]
The last step used the subgaussianity of $\pi$. Since this estimate clearly implies~\eqref{eq:local:factor:estimate:C}, the proof is complete.
\end{proof}

\section{Standard facts about AMP}
\label{app:scalar-facts}

\subsection{Convergence of Bayes AMP state evolution}
\label{app:scalar-limit}
The following lemma proves the state evolution description in \cref{sec:se} by showing that, under the initialization~\eqref{eq:q-init},
$q_t$ converges monotonically to $q_*$. 
\begin{lemma}
\label{lem:zero-fixed-point}
Assume that $\pi$ has finite moments of all orders with first moment $\mu$ and second moment $\lambda$.
Let $(q_t)_{t \geq 1}$ have the initialization~\eqref{eq:q-init} and recursion~\eqref{eq:phi},
and let $q_*$ be defined by~\eqref{eq:amp:fixed:point}.
Then the sequence $(q_t)_{t \geq 1}$ is monotonically increasing and converges to
$q_*\le\lambda$. More precisely:
\begin{enumerate}[label=(\roman*)]
 \item If $\mu\ne0$, then $q_*$ is the smallest nonnegative fixed point
 of $\Phi$, and $q_*\in[\mu^2,\lambda]$.
 \item If $\mu=0$ and $\lambda>1$, then $q_*$ is the smallest positive
 fixed point of $\Phi$, and $q_*\in[\lambda-\lambda^{-1},\lambda]$.
 \item If $\mu=0$ and $\lambda\le1$, then $q_t=q_*=0$ for all $t$.
\end{enumerate}
Furthermore, $(\Psi(q_t))_{t \geq 1}$ is monotonically decreasing and converges to $\Psi(q_*)$.
Finally, $\Phi(0)=\mu^2$ and $\Phi'(0)=(\lambda-\mu^2)^2$, with a
right derivative at zero.
\end{lemma}

\begin{proof}
The convergence argument in case~(ii) appears in
\cite[Appendix~F]{MV21}. Here, we include a proof covering
all three cases and the monotonicity of $\Psi(q_t)$.
For the definition of $\Phi(\cdot)$ in~\eqref{eq:scalar-channel} and $X_a=a\Theta+\sqrt{a}G$, $(\Theta,G)\sim \pi\otimes \cN(0,1)$, we have
\[
 \lambda-\Phi(a)=\E[(\Theta-\eta(X_a;a))^2]=\E\Big[\Big(\Theta-\E[\Theta\mid X_a]\Big)^2\Big].
 \]
 That is, $\lambda-\Phi(a)$ is a (scalar) minimum MSE. Thus, by standard properties of minimum MSE~\cite[Propositions~7 and~9]{GWSV11}, we have that $\Phi$ is continuous on $[0,\infty)$, smooth on $(0,\infty)$, and for $a>0$,
\[
 \Phi'(a)=\E[\operatorname{Var}(\Theta\mid X_a)^2]\ge0,\qquad
 \Phi'(0)=(\lambda-\mu^2)^2,
\]
where the derivative at zero is taken from the right.
Jensen's inequality and the tower property give
\[
\mu^2=\E[\E[\Theta\mid X_a]]^2\le\Phi(a)=\E[\E[\Theta\mid X_a]^2]
\le \lambda.
\]
Since $X_0=0$, $\Phi(0)=\mu^2$. Comparing the posterior mean with
estimators of the form $cX_a, c\in \R$, gives
\begin{equation}
 \Phi(a)\ge\frac{(\E[\Theta X_a])^2}{\E[X_a^2]}
 =\frac{a\lambda^2}{1+a\lambda},\qquad a>0.
 \label{eq:zero-linear}
\end{equation}
Since $\Psi'(q)=(q-\Phi(q))/2$~\eqref{eq:Psi-prime}, fixed points of $\Phi$ are exactly the stationary points
of $\Psi$, with $\Psi''=(1-\Phi')/2$.
Hence case~(iii) follows immediately from $q_1=0$ and
$\Phi'(0)=\lambda^2\le1$.

In case~(i), $q_1=\mu^2>0$ and $\Phi(q_1)\ge\mu^2=q_1$.
In case~(ii), $q_1=\lambda-\lambda^{-1}>0$, so
\eqref{eq:zero-linear} gives
\[
\Phi(q_1)\ge \frac{q_1\lambda^2}{1+\lambda q_1}=q_1.
\]
Thus, monotonicity and the bound $\Phi(a)\le\lambda$ imply
that $(q_t)$ is monotonically increasing. Moreover, by continuity, $(q_t)$ converges to $\bar{q}\leq \lambda$, which must then be the smallest fixed point of $\Phi(\cdot)$ above $q_1$. The preceding bounds also show that every positive fixed
point is at least $q_1$. Hence, $\bar q$ is the
smallest positive fixed point. Also, since $\Phi'(\bar q)\le1$ and in cases~(i) and~(ii), $0$ is either not a fixed point or has derivative
greater than one, we have $\bar q=q_*$. Finally, for $u\in[q_t,q_{t+1}]$, monotonicity gives
$\Phi(u)\ge\Phi(q_t)=q_{t+1}\ge u$, thus
\[
 \Psi(q_{t+1})-\Psi(q_t)
 =\frac12\int_{q_t}^{q_{t+1}}(u-\Phi(u))\,du
 \le0.
\]
Continuity of $\Psi$ completes the proof.
\end{proof}

\subsection{Limiting matrix MSE of Bayes AMP}
\label{app:amp-mse}

The following is essentially a consequence of~\cite{bayati2011dynamics,MV21} but this precise version does not seem to have appeared before.

\begin{proposition}[Matrix MSE of Bayes AMP]
\label{prop:amp-mse}
Assume the spiked Wigner model of \cref{def:wig} with a prior $\pi$ having finite moments of all orders. Let $(q_t)_{t\ge1}$ be defined by \eqref{eq:phi} and \eqref{eq:q-init}, and assume that $\eta_t=\eta(\cdot;q_t)$ is globally Lipschitz for every $t\ge1$.
Run Bayes AMP as described in \cref{sec:amp-intro}.
For every fixed integer $T\ge1$,
\begin{equation}\label{eq:amp-mse-fixed-T}
    \lim_{n\to\infty}\frac1{n^2}
    \E\bigl\|\hat\btheta^{(T)}(\hat\btheta^{(T)})^\T
            -\btheta\btheta^\T\bigr\|_F^2
    =\lambda^2-q_{T+1}^2.
\end{equation}
Consequently, with $q_*$ as defined in \eqref{eq:amp:fixed:point},
\begin{equation}\label{eq:amp-mse-double-limit}
    \lim_{T\to\infty}\lim_{n\to\infty}
    \MSE\bigl(\hat\btheta^{(T)}(\hat\btheta^{(T)})^\T;
        \btheta\btheta^\T\bigr)
    =\lambda^2-q_*^2=\MSE_\AMP.
\end{equation}
\end{proposition}

\begin{remark}[Lipschitz assumption]
The Lipschitz assumption is required by the state evolution theorem \cite[Theorem~1]{MV21}; the Lipschitz constant may depend on $t$. It also holds for bounded priors and for priors with log-concave densities; see \cite[Remark~2.3 and Lemma~F.1]{MV21}.
\end{remark}

\begin{proof}
If $\mu=0$ and $\lambda\le1$, all AMP estimates are zero, and the expected loss is
$(1-1/n)\lambda^2+\E[\Theta^4]/n\to\lambda^2$.
Otherwise $q_t>0$ for all $t\ge1$. State evolution follows from
\cite{bayati2011dynamics} when $\mu\ne0$ (see also
\cite[Remark~2.4]{MV21}), and from \cite[Theorem~1 and Section~2.4]{MV21}
for spectral initialization. In MV's notation, the signal is
$\btheta/\sqrt\lambda$, the iterate indexed by $t-1$ is
$\sqrt\lambda\,\bx^{(t)}$, and $\gamma_{t-1}=\lambda q_t$.

MV uses diagonal noise variance $2/n$, whereas $\bZ/\sqrt n$ has
diagonal variance $1/n$. Couple the models by adding an independent
diagonal matrix with entries $\mathcal N(0,1/n)$ to $\bY/\sqrt n$;
this perturbation has operator norm $o_p(1)$. Eigenvector perturbation,
using the positive limiting spectral gap when $\lambda>1$, and induction
in \eqref{eq:amp} imply that corresponding iterates differ by
$o_p(\sqrt n)$ in Euclidean norm.
Here we use that the denoisers are Lipschitz and their derivatives are
bounded and continuous; smoothness follows from \cite[Lemma~F.1]{MV21}.

Fix $T\ge1$ and write $\bh=\hat\btheta^{(T)}$, with the first denoising
step understood as $\hat\btheta^{(1)}=\eta_1(\bx^{(1)})$.
Applying state evolution to $\theta\eta_T(x)$ and $\eta_T(x)^2$, using
\eqref{eq:phi} and the sign convention in \cref{sec:amp-alg}, gives
\begin{equation}\label{eq:amp-mse-overlaps}
    \frac{|\langle\btheta,\bh\rangle|}{n}\stackrel p\longrightarrow q_{T+1},
    \qquad
    \frac{\|\bh\|_2^2}{n}\stackrel p\longrightarrow q_{T+1}.
\end{equation}
Also, $\|\btheta\|_2^2/n\to\lambda$ in probability by the law of large numbers. Thus, the identity
\begin{equation}\label{eq:amp-mse-loss-expansion}
    R_n:=\frac1{n^2}\|\bh\bh^\T-\btheta\btheta^\T\|_F^2
    =\left(\frac{\|\btheta\|_2^2}{n}\right)^2
     -2\left(\frac{\langle\btheta,\bh\rangle}{n}\right)^2
     +\left(\frac{\|\bh\|_2^2}{n}\right)^2
\end{equation}
gives $R_n\to\lambda^2-q_{T+1}^2$ in probability.

To pass to expectations, set $K_n=\|\bY/\sqrt n\|_{\mathrm{op}}$.
Since each $\eta_t$ grows at most linearly and $\eta_t'$ is bounded,
induction in \eqref{eq:amp}, starting from \eqref{eq:amp-init}, gives
\[
    \frac{\|\hat\btheta^{(T)}\|_2}{\sqrt n}
    \le C_T(1+K_n)^T
\]
for a constant $C_T$ independent of $n$.
The moment assumption and standard Gaussian operator-norm bounds imply
that $\|\btheta\|_2^2/n$ and
$K_n\le\|\btheta\|_2^2/n+\|\bZ/\sqrt n\|_{\mathrm{op}}$
have moments of every fixed order bounded uniformly in $n$.
Thus \eqref{eq:amp-mse-loss-expansion} gives
$\sup_n\E[R_n^2]<\infty$, so $(R_n)_n$ is uniformly integrable.
Its convergence in probability therefore proves \eqref{eq:amp-mse-fixed-T}.
Finally, \cref{lem:zero-fixed-point} gives $q_{T+1}\to q_*$,
proving~\eqref{eq:amp-mse-double-limit}.
\end{proof}

\subsection{Polynomial approximation of the state evolution}
\label{app:scalar-approximation}

\begin{proof}[Proof of \Cref{lem:approx-scalar}]
Recall $X_a$ and $\hat\eta(\cdot;a,d)$ from
\eqref{eq:scalar-channel} and~\eqref{eq:least-squares}, respectively.
If $\lambda=0$, then $\Theta=0$ almost surely and
$\eta=\hat\eta=0$, thus $a_t=q_t=0$ for all $t\geq 2$. Henceforth assume $\lambda>0$. By \Cref{lem:calibration},
\[
 \Phi_d(a):=\E[\hat\eta(X_a;a,d)^2]
 =\lambda-\min_{\deg p\le d}\E[(\Theta-p(X_a))^2],
 \qquad a_{t+1}=\Phi_d(a_t).
\]
Since $d\ge1$, the same linear estimator comparison
as in~\eqref{eq:zero-linear} gives
\[
 0<\frac{a\lambda^2}{1+a\lambda}
 \le\Phi_d(a)\le\Phi(a)\le\lambda,
 \qquad a>0,\quad d\ge1.
\]
In particular, $a_t,q_t>0$ for all $t\geq 1$.

We claim that $\Phi_d(b_d)\to\Phi(b)$ whenever $b_d\to b>0$ as $d\to\infty$.
The upper bound follows from $\Phi_d\le\Phi$ and continuity of $\Phi$.
For the lower bound, fix any polynomial $p$. For all $d\ge\deg p$,
\[
 \Phi_d(b_d)\ge\lambda-\E[(\Theta-p(X_{b_d}))^2]
 \longrightarrow\lambda-\E[(\Theta-p(X_b))^2].
\]
The law of $X_b$ has a finite exponential moment, since
$\E e^{\delta|X_b|}\le
\E e^{\delta b|\Theta|}\E e^{\delta\sqrt b|G|}<\infty$
for sufficiently small $\delta>0$.
Therefore polynomials are dense in $L^2$ of this law by
\cite[Theorem~1.1]{deJeu03}.
Taking the supremum over $p$ in the lower bound yields
$\liminf_d\Phi_d(b_d)\ge\Phi(b)$, proving the claim.
Starting from $a_1=q_1$, induction now gives
$a_{t+1}=\Phi_d(a_t)\to\Phi(q_t)=q_{t+1}$ for every fixed $t$.
\end{proof}

\subsection{Exceptional signal strengths}
\label{app:scalar-exceptional}

\begin{lemma}[Exceptional signal strengths]
\label{lem:critical-parameters}
Fix a subgaussian random variable $\widetilde\Theta$ with
$\E\widetilde\Theta^2=1$. For $\lambda>0$, let $\Phi_\lambda$
be the map~\eqref{eq:phi} for the prior
$\Theta_\lambda=\sqrt\lambda\widetilde\Theta$.
The set
\[
 \mathcal E:=\{\lambda>0:\text{there exists }q\ge0\text{ with }
 \Phi_\lambda(q)=q,\ \Phi_\lambda'(q)=1\}
\]
has Lebesgue measure zero. In particular,
$\Psi''(q_*;\pi_\lambda)>0$ for almost every $\lambda>0$,
where $\pi_\lambda$ is the law of $\Theta_\lambda$.
\end{lemma}

\begin{proof}
Let $f(s):=\Phi_1(s)$ for $s>0$, where $\Phi_1$ is the
state-evolution map for the prior $\widetilde\Theta$. As established in the proof of \cref{lem:zero-fixed-point}, $f$ is smooth and $f(s)\ge s/(1+s)>0$. Note that since $\Theta_{\lambda}=\sqrt{\lambda}\widetilde\Theta$,
\[
\Phi_\lambda(q)
=\lambda\E\left[\E\big[\widetilde\Theta\big|
\sqrt{\lambda q}\widetilde\Theta+G\big]^2\right]
=\lambda f(\lambda q).
\]
Thus, for $q>0$, the two critical fixed point conditions, with $s=\lambda q$,
are $ s=\lambda^2f(s)$ and $1=\lambda^2f'(s)$. Consequently, 
\[
\mathcal{E} \setminus \{1\} \subseteq \{L(s):s>0,\ L'(s)=0\},
\quad\textnormal{where}\quad L(s):=\sqrt{s/f(s)}.
\]
By Sard's theorem the set $\{L(s):s>0,\ L'(s)=0\}$ has Lebesgue measure zero. Moreover, at $q=0$, a critical fixed point occurs only when
$\E\widetilde\Theta=0$ and $\lambda=1$, so this contributes
at most one additional exceptional value. Therefore $\mathcal{E}$ has Lebesgue measure zero. The last assertion then follows from~\eqref{eq:Psi-prime}.
\end{proof}

\section{Achieving Bayes error at degree \texorpdfstring{$O(n)$}{O(n)}}
\label{app:bayes}

In this appendix we show that, for almost every signal
strength, polynomials of degree $O_\eps(n)$ achieve the Bayes
error up to any fixed additive error $\eps>0$. We will assume the prior is bounded to keep the proof simple, although we expect this can be relaxed to, say, subgaussian priors, with a more involved proof.

\subsection{Statement of result}

Fix a random variable $\widetilde\Theta\sim\widetilde\pi$ satisfying
\begin{equation}\label{eq:bayes-bounded-prior}
    \E[\widetilde\Theta^2]=1,
    \qquad |\widetilde\Theta|\le K\quad\text{almost surely},
\end{equation}
for some $K<\infty$. For $\lambda>0$, let $\pi_\lambda$ be the law of
$\sqrt\lambda\,\widetilde\Theta$, and consider the spiked Wigner model (\cref{def:wig}) with prior $\pi_\lambda$.
We will use the matrix MSE from \cref{eq:mat-mse}. To display the dependence on $n$ and
$\lambda$, write
\begin{align}
    \MMSE_n(\lambda)
    &:=\inf_{\bmf\text{ measurable}}
        \MSE(\bmf;\btheta\btheta^\T),
        \label{eq:bayes-mmse}\\
    \MMSE_{n,\le D}(\lambda)
    &:=\MMSE_{\le D}(\btheta\btheta^\T).
        \label{eq:bayes-degree-mmse}
\end{align}

\begin{theorem}\label{thm:bayes-polynomial}
Under \eqref{eq:bayes-bounded-prior}, there is an at-most-countable set
$\calE\subset(0,\infty)$ such that, for every fixed
$\lambda\notin\calE$, the limit
\[
    \MSE_\Bayes(\lambda):=\lim_{n\to\infty}\MMSE_n(\lambda)
\]
exists, and for every $\eps>0$ there is a finite $C_\eps=C_\eps(\lambda,\widetilde\pi)$ such that
\begin{equation}\label{eq:bayes-main-bound}
    \limsup_{n\to\infty}
        \MMSE_{n,\le C_\eps n}(\lambda)
    \le \MSE_\Bayes(\lambda)+\eps.
\end{equation}
\end{theorem}

The exceptional set $\calE$ consists exactly of those $\lambda>0$ for
which $\Psi(\cdot;\pi_\lambda)$ has more than one global minimizer, as in~\cite[Section~3.2]{lelarge2019fundamental}.
Whenever $\lambda \notin \calE$, we have
$\MSE_\Bayes(\lambda)=\lambda^2-q_\IT^2$ as in the introduction.

\subsection{Proof overview}
\label{sec:bayes-proof}

The main step is to construct a polynomial estimator whose error at signal strength
$\lambda$ is at most the Bayes error at a slightly weaker signal strength $a\lambda$ for $a \in (0,1)$.

\begin{proposition}[Comparison with a weaker signal]
\label{prop:bayes-weaker-signal}
Under \eqref{eq:bayes-bounded-prior}, for every fixed $\lambda>0$ and
$a\in(0,1)$ there is a finite $C=C(a,\lambda,\widetilde\pi)$ such that
\begin{equation}\label{eq:bayes-weaker-bound}
    \MMSE_{n,\le Cn}(\lambda)
    \le a^{-2}\MMSE_n(a\lambda)+o(1).
\end{equation}
\end{proposition}

The factor $a^{-2}$ accounts for the change in MSE when the estimand is scaled by $a$. The proof of \Cref{prop:bayes-weaker-signal} is deferred to \Cref{sec:bayes-weaker-signal}. Once we have this, the final ingredient is continuity of the Bayes MSE, which will allow us to take $a\uparrow1$.

\begin{lemma}[Continuity at non-exceptional signal strengths]\label{lem:bayes-regularity}
Under \eqref{eq:bayes-bounded-prior}, there is an at-most-countable set
$\calE\subset(0,\infty)$ such that, for every $\lambda\notin\calE$,
the limit $\MSE_\Bayes(\lambda)=\lim_{n\to\infty}\MMSE_n(\lambda)$
exists. Moreover, if $a_k\uparrow1$ and $a_k\lambda\notin\calE$, then
\begin{equation}\label{eq:bayes-left-continuity}
    a_k^{-2}\MSE_\Bayes(a_k\lambda)\longrightarrow\MSE_\Bayes(\lambda).
\end{equation}
\end{lemma}

The proof is deferred to \cref{sec:bayes-regularity}. The two ingredients above imply the desired result.

\begin{proof}[Proof of \cref{thm:bayes-polynomial}]
Fix $\lambda\notin\calE$ and $\eps>0$. Since $\calE$ is countable,
\cref{lem:bayes-regularity} allows us to choose $a\in(0,1)$ with
$a\lambda\notin\calE$ such that $a^{-2}\MSE_\Bayes(a\lambda)\le\MSE_\Bayes(\lambda)+\eps$. Taking the limit superior in \eqref{eq:bayes-weaker-bound} proves
\eqref{eq:bayes-main-bound} with $C_\eps=C(a,\lambda,\widetilde\pi)$.
\end{proof}

\subsection{Comparing to a weaker signal}
\label{sec:bayes-weaker-signal}

The argument we will use here applies to a more general additive Gaussian setting, which we now describe. Fix $a\in(0,1)$ and $B > 0$, and consider a random variable $\bX \in \R^m$ representing the ``signal,'' with the bound $\|\bX\|_2^2 \le Bn$ almost surely. Let $\bG,\bG'$ be independent standard Gaussian vectors
in $\R^m$, independent of $\bX$, and define
\[
    \bY=\bX+\bG,
    \qquad \bY'=a\bY+\sqrt{1-a^2}\,\bG'.
\]
Let $\bmf: \R^m \to \R^m$ be an arbitrary measurable estimator obeying the bound $\|\bmf(\by)\|_2^2 \le Bn$ for every $\by\in\R^m$. Using only the assumptions above, we will conclude there is a vector of polynomials $\bp: \R^m \to \R^m$ of degree $\le Cn$ where $C = C(a,B) > 0$, such that
\begin{equation}\label{eq:weaker-signal}
\E\|\bmp(\bY)-\bX\|_2^2 \le \E\|\bmf(\bY')-\bX\|_2^2+o(1).
\end{equation}
Here, $o(1)$ pertains to the limit $n\to\infty$ with $a,B$ fixed, while $m, \bX, \bmf$ may depend on $n$.

Once we have this, \cref{prop:bayes-weaker-signal} follows by taking $\bX$ to be the upper triangular part of the signal matrix $\btheta\btheta^\T/\sqrt{n}$, taking $\bmf$ to be the Bayes estimator $\bmf(\by)=\E[\bX\mid\bY'=\by]$, and $B=\lambda^2 K^4$. Note that $\bY'$ is distributed as an observation with signal strength $a\lambda$.

The rest of this section is devoted to the proof of~\eqref{eq:weaker-signal}. We proceed in three steps.

\paragraph{Step 1: Smooth the estimator.}
Define a smoothed estimator $\bmh := T_a\bmf$ where the Gaussian noise operator is defined as
\[
    T_t\bmg(\by)
    :=\E_{\bG'}\bmg\bigl(t\by+\sqrt{1-t^2}\,\bG'\bigr) \qquad \text{for } t \in (0,1).
\]
Conditional Jensen's inequality, with $(\bX,\bY)$ fixed, yields
\begin{equation}\label{eq:bayes-smoothing-risk}
    \E\|\bmh(\bY)-\bX\|_2^2
    \le \E\|\bmf(\bY')-\bX\|_2^2.
\end{equation}
It remains to approximate $\bmh$ by a polynomial of degree $Cn$.

\paragraph{Step 2: Approximate under Gaussian noise.}
Let $\Qb=\calN(0,I_m)$ and let $\Pb$ be the law of $\bY$.
For vector-valued functions, write
\[
    \|\bmg\|_{L^p(P)}:=\bigl(\E_P\|\bmg\|_2^p\bigr)^{1/p}.
\]
Since $\|\bmf\|_{L^2(\Qb)}\le\sqrt{Bn}$,
we may expand $\bmf=\sum_{k\ge0}\bmf_k$ in $L^2(\Qb)$, where
$\bmf_k$ is the coordinatewise orthogonal projection onto the span
of Hermite polynomials of total degree $k$. The well-known identity
$T_t\bmf_k=t^k\bmf_k$ gives
\[
    \bmh=\sum_{k\ge0}a^k\bmf_k,
    \qquad
    \bmp_d:=\sum_{k=0}^d a^k\bmf_k.
\]
Here, $\bp_d$ is a polynomial approximation to $\bh$, obtained by truncating the Hermite expansion. Next, we will bound the truncation error in $L^p(\Qb)$ for some $p>2$. Step~3 below will transfer this to a bound in the desired norm, $L^2(\Pb)$.
Choose $2<p<1+a^{-2}$ depending only on $a$, and put
\[
    t=(p-1)^{-1/2},
    \qquad r=a\sqrt{p-1}<1,
\]
so $tr=a$. Gaussian hypercontractivity~\cite{gross1975logarithmic}
states that $\|T_tg\|_{L^p(\Qb)}\le\|g\|_{L^2(\Qb)}$ for
scalar $g$. The same bound holds for vector-valued $\bmg$, since
the triangle inequality gives
$\|T_t\bmg(\by)\|_2\le T_t(\|\bmg\|_2)(\by)$, and hence
\[
    \|T_t\bmg\|_{L^p(\Qb)}
    \le\|T_t(\|\bmg\|_2)\|_{L^p(\Qb)}
    \le\|\bmg\|_{L^2(\Qb)}.
\]
Applying this bound and orthogonality yields
\begin{align}
    \|\bmh-\bmp_d\|_{L^p(\Qb)}
    &=\left\|T_t\!\left(\sum_{k>d}r^k\bmf_k\right)
        \right\|_{L^p(\Qb)}\notag\\
    &\le\left(\sum_{k>d}r^{2k}\|\bmf_k\|_{L^2(\Qb)}^2\right)^{1/2}
        \notag\\
    &\le r^{d+1}\|\bmf\|_{L^2(\Qb)}
    \le\sqrt{Bn}\,r^{d+1}.
    \label{eq:bayes-gaussian-error}
\end{align}

\paragraph{Step 3: Transfer to the spiked law.}
Let $L=\dd\Pb/\dd\Qb$ and $q=p/(p-2)$.
For $\bme=\bmp_d-\bmh$, H\"older's inequality gives
\[
    \|\bme\|_{L^2(\Pb)}^2
    =\E_{\Qb}[L\|\bme\|_2^2]
    \le\|L\|_{L^q(\Qb)}\|\bme\|_{L^p(\Qb)}^2.
\]
To bound the likelihood moment, condition on the signal $\bX$.
For a fixed shift $\bx$, the likelihood of $\calN(\bx,I_m)$
relative to $\Qb$ is
\[
    L_{\bx}(\by)
    =\exp\bigl(\langle\by,\bx\rangle-\|\bx\|_2^2/2\bigr).
\]
Gaussian integration gives
\[
    \E_{\Qb}[L_{\bx}^q]
    =\exp\!\left(\frac{q(q-1)}2\|\bx\|_2^2\right),
    \qquad
    \|L_{\bx}\|_{L^q(\Qb)}
    =\exp\!\left(\frac{q-1}2\|\bx\|_2^2\right).
\]
Since $L=\E_{\bX}L_{\bX}$, Minkowski's inequality and
$\|\bX\|_2^2\le Bn$ imply
\[
    \|L\|_{L^q(\Qb)}
    \le\E_{\bX}\|L_{\bX}\|_{L^q(\Qb)}
    \le\exp\!\left(\frac{q-1}{2}Bn\right).
\]
Combining these bounds with \eqref{eq:bayes-gaussian-error} gives
\begin{equation}\label{eq:bayes-planted-error}
    \|\bmp_d-\bmh\|_{L^2(\Pb)}
    \le\sqrt{Bn}\,\exp\!\left(\frac{q-1}{4}Bn\right)r^{d+1}.
\end{equation}
Choose
\[
    C>\frac{(q-1)B/4+1}{|\log r|},
    \qquad d=\lfloor Cn\rfloor,
\]
and set $\bmp=\bmp_d$. Then $C$ depends only on $a,B$, and
\eqref{eq:bayes-planted-error} is at most $\sqrt{Bn}\,e^{-n}$.
The triangle inequality and \eqref{eq:bayes-smoothing-risk} give
\[
    \bigl(\E\|\bmp(\bY)-\bX\|_2^2\bigr)^{1/2}
    \le\bigl(\E\|\bmf(\bY')-\bX\|_2^2\bigr)^{1/2}
        +\sqrt{Bn}\,e^{-n}.
\]
Since $\|\bmf(\by)\|_2,\|\bX\|_2\le\sqrt{Bn}$, the first term on the
right is at most $2\sqrt{Bn}$. Squaring yields
\[
    \E\|\bmp(\bY)-\bX\|_2^2
    \le\E\|\bmf(\bY')-\bX\|_2^2
        +4Bn e^{-n}+Bn e^{-2n},
\]
which proves \eqref{eq:weaker-signal}.

\subsection{Proof of Lemma~\ref{lem:bayes-regularity}}
\label{sec:bayes-regularity}

This is essentially a consequence of~\cite{lelarge2019fundamental}. Let $J(s)$ denote the limit $\lim_{n\to\infty}n^{-1}I(X;Y)$ in
Theorem~1 of \cite{lelarge2019fundamental} (arXiv v3),
with $P_0=\widetilde\pi$ and their signal parameter $\lambda$ replaced
by $s$, so $s=\lambda^2$ in our notation.
That theorem and Section~3.2 of the same reference imply that $J$ is
finite and concave, and for every $\lambda$ outside the at-most-countable set
\begin{equation}\label{eq:bayes-exceptional-set}
    \calE:=\{\lambda>0:J\text{ is not differentiable at }\lambda^2\},
\end{equation}
we have, in our normalization,
\begin{equation}\label{eq:bayes-limit-risk}
    \MMSE_n(\lambda)\longrightarrow
    \MSE_\Bayes(\lambda)=4\lambda^2J'(\lambda^2).
\end{equation}
Their discussion following equation~(2) also
allows diagonal observations. The equivalence between differentiability
and uniqueness of the optimizer, stated immediately before
Proposition~15 in Section~3.2, identifies $\calE$ with the set where
$\Psi(\cdot;\pi_\lambda)$ has multiple global minimizers.

Fix $\lambda\notin\calE$ and $a_k\uparrow1$ with
$a_k\lambda\notin\calE$. A standard fact for concave functions is that
$J'(s_k)\to J'(s)$ whenever $s_k\to s$ and $J$ is differentiable at $s$
and at every $s_k$. Applying this with $s_k=a_k^2\lambda^2$ and
$s=\lambda^2$ gives
\[
    a_k^{-2}\MSE_\Bayes(a_k\lambda)
    =4\lambda^2J'(a_k^2\lambda^2)
    \longrightarrow4\lambda^2J'(\lambda^2)
    =\MSE_\Bayes(\lambda).
\]

\phantomsection
\addcontentsline{toc}{section}{Statement on AI use}

\subsection*{Statement on AI use}
We used ChatGPT 5.5 and 5.6 Pro in developing the proof of the lower bound in \cref{thm:lowdeg}-\ref{item:lower}.
The human authors first derived the tree contributions
described in \cref{sec:stein-via-amp}.
Subsequent discussions with ChatGPT about extending these
calculations to a dual certificate containing cyclic terms
led us to identify the relevance of the Stein--cumulant
identity~\cite{barbour86asymptotic} and to develop the
history certificate in
\cref{subsec:outline:history:certificate}.
Further discussions helped us analyze the norm of this
certificate and identify the role of the stability condition
$\Psi''(q_*;\pi)>0$ in controlling the cyclic contributions
for degrees $D=o(n)$, as established in
\cref{prop:cyclic:norm:bound}.  We also used ChatGPT 6 to assist with editing, figure preparation, and also proving the upper bounds in \cref{thm:lowdeg}-\ref{item:upper} and \cref{thm:bayes-polynomial}. All AI-assisted calculations and suggestions were carefully checked by the authors, who take full responsibility for the contents of the paper.

\phantomsection
\addcontentsline{toc}{section}{Acknowledgments}

\subsection*{Acknowledgments}
YS was partially supported by NSF Award DMS-2601870. AW was partially supported by a Sloan Fellowship and NSF CAREER Award CCF-2338091. This work was supported in part by a grant of access to OpenAI models through the ChatGPT for Academic Researchers program.

YS and AW thank participants at a 2024 AIM workshop (American Institute of Mathematics) where some preliminary discussions related to this question occurred.

\phantomsection
\addcontentsline{toc}{section}{References}

{\small
\bibliographystyle{alpha}
\bibliography{main}

@article{even2025computational,
  author = {Even, Bertrand and Giraud, Christophe and Verzelen, Nicolas},
  title = {Computational lower bounds in latent models: clustering, sparse-clustering, biclustering},
  journal = {arXiv preprint arXiv:2506.13647},
  year = {2025}
}

@article{barbour86asymptotic,
  author = {Barbour, A. D.},
  title = {Asymptotic expansions based on smooth functions in the central limit theorem},
  journal = {Probability Theory and Related Fields},
  volume = {72},
  number = {2},
  pages = {289--303},
  year = {1986},
  doi = {10.1007/BF00699108}
}

@inproceedings{even2024computation,
  author = {Even, Bertrand and Giraud, Christophe and Verzelen, Nicolas},
  title = {Computation-information gap in high-dimensional clustering},
  booktitle = {Proceedings of the 37th Conference on Learning Theory},
  series = {Proceedings of Machine Learning Research},
  volume = {247},
  pages = {1646--1712},
  publisher = {PMLR},
  year = {2024},
  url = {https://proceedings.mlr.press/v247/even24a.html}
}

@article{fu2026low,
  author = {Fu, Daniel and Sohn, Youngtak},
  title = {Low-degree estimation thresholds in planted hypergraphs and tensor {PCA}},
  journal = {arXiv preprint arXiv:2605.30113},
  year = {2026}
}

@inproceedings{el2018estimation,
  author = {El Alaoui, Ahmed and Krzakala, Florent},
  title = {Estimation in the spiked {Wigner} model: A short proof of the replica formula},
  booktitle = {IEEE International Symposium on Information Theory (ISIT)},
  pages = {1874--1878},
  publisher = {IEEE},
  year = {2018}
}

@inproceedings{lesieur2016phase,
  author = {Lesieur, Thibault and De Bacco, Caterina and Banks, Jess and Krzakala, Florent and Moore, Cris and Zdeborov{\'a}, Lenka},
  title = {Phase transitions and optimal algorithms in high-dimensional {Gaussian} mixture clustering},
  booktitle = {54th Annual Allerton Conference on Communication, Control, and Computing (Allerton)},
  pages = {601--608},
  publisher = {IEEE},
  year = {2016}
}

@article{bayati2011dynamics,
  author = {Bayati, Mohsen and Montanari, Andrea},
  title = {The dynamics of message passing on dense graphs, with applications to compressed sensing},
  journal = {IEEE Transactions on Information Theory},
  volume = {57},
  number = {2},
  pages = {764--785},
  year = {2011}
}

@article{bolthausen2014iterative,
  author = {Bolthausen, Erwin},
  title = {An iterative construction of solutions of the {TAP} equations for the {Sherrington--Kirkpatrick} model},
  journal = {Communications in Mathematical Physics},
  volume = {325},
  number = {1},
  pages = {333--366},
  year = {2014}
}

@article{thouless1977solution,
  author = {D. J. Thouless and P. W. Anderson and R. G. Palmer},
  title = {Solution of `solvable model of a spin glass'},
  journal = {Philosophical Magazine},
  volume = {35},
  number = {3},
  pages = {593--601},
  year = {1977},
  doi = {10.1080/14786437708235992}
}

@article{donoho2009message,
  author = {David L. Donoho and Arian Maleki and Andrea Montanari},
  title = {Message-passing algorithms for compressed sensing},
  journal = {Proceedings of the National Academy of Sciences},
  volume = {106},
  number = {45},
  pages = {18914--18919},
  year = {2009},
  doi = {10.1073/pnas.0909892106},
  url = {https://www.pnas.org/doi/abs/10.1073/pnas.0909892106}
}

@article{krauthgamer2015semidefinite,
  author = {Robert Krauthgamer and Boaz Nadler and Dan Vilenchik},
  title = {Do semidefinite relaxations solve sparse {PCA} up to the information limit?},
  journal = {The Annals of Statistics},
  volume = {43},
  number = {3},
  pages = {1300--1322},
  year = {2015},
  doi = {10.1214/15-AOS1310}
}

@article{lesieur2017constrained,
  author = {Lesieur, Thibault and Krzakala, Florent and Zdeborov{\'a}, Lenka},
  title = {Constrained low-rank matrix estimation: Phase transitions, approximate message passing and applications},
  journal = {Journal of Statistical Mechanics: Theory and Experiment},
  volume = {2017},
  number = {7},
  pages = {073403},
  year = {2017}
}

@inproceedings{brennan19b,
  author = {Brennan, Matthew and Bresler, Guy},
  title = {Optimal Average-Case Reductions to Sparse {PCA}: From Weak Assumptions to Strong Hardness},
  booktitle = {Proceedings of the Thirty-Second Conference on Learning Theory},
  series = {Proceedings of Machine Learning Research},
  volume = {99},
  pages = {469--470},
  publisher = {PMLR},
  year = {2019},
  url = {https://proceedings.mlr.press/v99/brennan19b.html}
}

@article{BMVVX,
  author = {Banks, Jess and Moore, Cristopher and Vershynin, Roman and Verzelen, Nicolas and Xu, Jiaming},
  title = {Information-theoretic bounds and phase transitions in clustering, sparse {PCA}, and submatrix localization},
  journal = {IEEE Transactions on Information Theory},
  volume = {64},
  number = {7},
  pages = {4872--4894},
  year = {2018}
}

@article{gamarnik2022sparse,
  author = {David Gamarnik and Ilias Zadik},
  title = {Sparse high-dimensional linear regression. Estimating squared error and a phase transition},
  journal = {The Annals of Statistics},
  volume = {50},
  number = {2},
  pages = {880--903},
  year = {2022},
  doi = {10.1214/21-AOS2130}
}

@article{celentano2022fundamental,
  author = {Michael Celentano and Andrea Montanari},
  title = {Fundamental barriers to high-dimensional regression with convex penalties},
  journal = {The Annals of Statistics},
  volume = {50},
  number = {1},
  pages = {170--196},
  year = {2022},
  doi = {10.1214/21-AOS2100}
}

@book{billingsley1995probability,
  author = {Billingsley, Patrick},
  title = {Probability and Measure},
  edition = {third},
  publisher = {John Wiley \& Sons},
  year = {1995}
}

@article{SW-22,
  author = {Tselil Schramm and Alexander S. Wein},
  title = {Computational barriers to estimation from low-degree polynomials},
  journal = {The Annals of Statistics},
  volume = {50},
  number = {3},
  pages = {1833--1858},
  year = {2022},
  doi = {10.1214/22-AOS2179}
}

@article{nourdin2010cumulants,
  author = {Nourdin, Ivan and Peccati, Giovanni},
  title = {Cumulants on the {Wiener} space},
  journal = {Journal of Functional Analysis},
  volume = {258},
  number = {11},
  pages = {3775--3791},
  year = {2010}
}

@inproceedings{dia2016mutual,
  author = {Barbier, Jean and Dia, Mohamad and Macris, Nicolas and Krzakala, Florent and Lesieur, Thibault and Zdeborov{\'a}, Lenka},
  title = {Mutual information for symmetric rank-one matrix estimation: A proof of the replica formula},
  booktitle = {Advances in Neural Information Processing Systems},
  volume = {29},
  pages = {424--432},
  year = {2016}
}

@article{wein2019kikuchi,
  author = {Wein, Alexander S. and {El Alaoui}, Ahmed and Moore, Cristopher},
  title = {The {Kikuchi} hierarchy and tensor {PCA}},
  journal = {Journal of the ACM},
  volume = {72},
  number = {5},
  pages = {35:1--35:40},
  year = {2025},
  doi = {10.1145/3762806}
}

@article{wein2025computational,
  author = {Wein, Alexander S.},
  title = {Computational complexity of statistics: New insights from low-degree polynomials},
  journal = {arXiv preprint arXiv:2506.10748},
  year = {2025}
}

@article{li2026almost,
  author = {Li, Zhangsong},
  title = {Almost Sharp Equivalence between Approximate Message Passing and Low-Degree Polynomials},
  journal = {arXiv preprint arXiv:2609.06988},
  year = {2026}
}

@inproceedings{sohn2025sharp,
  author = {Sohn, Youngtak and Wein, Alexander S.},
  title = {Sharp phase transitions in estimation with low-degree polynomials},
  booktitle = {Proceedings of the 57th Annual ACM Symposium on Theory of Computing},
  pages = {891--902},
  publisher = {ACM},
  year = {2025},
  doi = {10.1145/3717823.3718294}
}

@article{bayati2015universality,
  author = {Bayati, Mohsen and Lelarge, Marc and Montanari, Andrea},
  title = {Universality in polytope phase transitions and message passing algorithms},
  journal = {The Annals of Applied Probability},
  volume = {25},
  number = {2},
  pages = {753--822},
  year = {2015},
  doi = {10.1214/14-AAP1010}
}

@article{montanari2022equivalence,
  author = {Montanari, Andrea and Wein, Alexander S.},
  title = {Equivalence of approximate message passing and low-degree polynomials in rank-one matrix estimation},
  journal = {Probability Theory and Related Fields},
  volume = {191},
  pages = {181--233},
  year = {2025},
  doi = {10.1007/s00440-024-01322-z}
}

@inproceedings{jones2025fourier,
  author = {Jones, Chris and Pesenti, Lucas},
  title = {{Fourier} analysis of iterative algorithms},
  booktitle = {52nd International Colloquium on Automata, Languages, and Programming (ICALP)},
  series = {Leibniz International Proceedings in Informatics (LIPIcs)},
  volume = {334},
  pages = {102:1--102:21},
  publisher = {Schloss Dagstuhl--Leibniz-Zentrum f{\"u}r Informatik},
  year = {2025},
  doi = {10.4230/LIPIcs.ICALP.2025.102}
}

@article{GWSV11,
  author = {Guo, D. and Wu, Y. and Shamai, S. and Verd{\'u}, S.},
  title = {Estimation in {Gaussian} Noise: Properties of the Minimum Mean-Square Error},
  journal = {IEEE Transactions on Information Theory},
  volume = {57},
  number = {4},
  pages = {2371--2385},
  year = {2011}
}

@article{deJeu03,
  author = {de Jeu, M.},
  title = {Determinate multidimensional measures, the extended {Carleman} theorem and quasi-analytic weights},
  journal = {The Annals of Probability},
  volume = {31},
  number = {3},
  pages = {1205--1227},
  year = {2003}
}

@article{MV21,
  author = {Montanari, Andrea and Venkataramanan, Ramji},
  title = {Estimation of low-rank matrices via approximate message passing},
  journal = {The Annals of Statistics},
  volume = {49},
  number = {1},
  pages = {321--345},
  year = {2021},
  doi = {10.1214/20-AOS1958}
}

@phdthesis{hopkins-thesis,
  author = {Hopkins, Samuel B.},
  title = {Statistical Inference and the Sum of Squares Method},
  school = {Cornell University},
  year = {2018},
  doi = {10.7298/X4416V82},
  url = {https://www.samuelbhopkins.com/thesis.pdf}
}

@incollection{wu-xu-survey,
  author = {Wu, Yihong and Xu, Jiaming},
  title = {Statistical Problems with Planted Structures: Information-Theoretical and Computational Limits},
  booktitle = {Information-Theoretic Methods in Data Science},
  pages = {383--424},
  publisher = {Cambridge University Press},
  year = {2021},
  doi = {10.1017/9781108616799.014}
}

@incollection{ld-notes,
  author = {Kunisky, Dmitriy and Wein, Alexander S. and Bandeira, Afonso S.},
  title = {Notes on Computational Hardness of Hypothesis Testing: Predictions Using the Low-Degree Likelihood Ratio},
  booktitle = {Mathematical Analysis, Its Applications and Computation},
  series = {Springer Proceedings in Mathematics \& Statistics},
  volume = {385},
  pages = {1--50},
  publisher = {Springer},
  year = {2022},
  doi = {10.1007/978-3-030-97127-4_1}
}

@article{disordered-insights,
  author = {Gamarnik, David and Moore, Cristopher and Zdeborov{\'a}, Lenka},
  title = {Disordered Systems Insights on Computational Hardness},
  journal = {Journal of Statistical Mechanics: Theory and Experiment},
  volume = {2022},
  number = {11},
  pages = {114015},
  year = {2022},
  doi = {10.1088/1742-5468/ac9cc8}
}

@misc{kz-notes,
  author = {Krzakala, Florent and Zdeborov{\'a}, Lenka},
  title = {Statistical Physics Methods in Optimization and Machine Learning: An Introduction to Replica, Cavity and Message-Passing Techniques},
  year = {2024},
  howpublished = {\url{https://sphinxteam.github.io/EPFLDoctoralLecture2022/Notes.pdf}},
  note = {Lecture notes, March 9, 2024}
}

@misc{zadik-notes,
  author = {Zadik, Ilias},
  title = {An Introduction to Computational-Statistical Tradeoffs in High-Dimensional Statistics},
  year = {2026},
  howpublished = {\url{https://iliaszadik.github.io/files/Book_v1.pdf}},
  note = {Draft book/lecture notes, June 2026}
}

@article{amp-tutorial,
  author = {Feng, Oliver Y. and Venkataramanan, Ramji and Rush, Cynthia and Samworth, Richard J.},
  title = {A unifying tutorial on approximate message passing},
  journal = {Foundations and Trends in Machine Learning},
  volume = {15},
  number = {4},
  pages = {335--536},
  year = {2022},
  doi = {10.1561/2200000092}
}

@article{FR18,
  author = {Fletcher, Alyson K. and Rangan, Sundeep},
  title = {Iterative reconstruction of rank-one matrices in noise},
  journal = {Information and Inference: A Journal of the IMA},
  volume = {7},
  number = {3},
  pages = {531--562},
  year = {2018},
  doi = {10.1093/imaiai/iax014}
}

@inproceedings{DM14,
  author = {Deshpande, Yash and Montanari, Andrea},
  title = {Information-theoretically optimal sparse {PCA}},
  booktitle = {IEEE International Symposium on Information Theory (ISIT)},
  pages = {2197--2201},
  publisher = {IEEE},
  year = {2014},
  doi = {10.1109/ISIT.2014.6875223}
}

@article{LFW23,
  author = {Li, Gen and Fan, Wei and Wei, Yuting},
  title = {Approximate message passing from random initialization with applications to {$\mathbb{Z}_2$} synchronization},
  journal = {Proceedings of the National Academy of Sciences},
  volume = {120},
  number = {31},
  pages = {e2302930120},
  year = {2023},
  doi = {10.1073/pnas.2302930120}
}

@article{BBP,
  author = {Baik, Jinho and {Ben Arous}, G{\'e}rard and P{\'e}ch{\'e}, Sandrine},
  title = {Phase transition of the largest eigenvalue for nonnull complex sample covariance matrices},
  journal = {The Annals of Probability},
  volume = {33},
  number = {5},
  pages = {1643--1697},
  year = {2005},
  doi = {10.1214/009117905000000233}
}

@article{FP-wigner,
  author = {F{\'e}ral, Delphine and P{\'e}ch{\'e}, Sandrine},
  title = {The largest eigenvalue of rank one deformation of large {Wigner} matrices},
  journal = {Communications in Mathematical Physics},
  volume = {272},
  number = {1},
  pages = {185--228},
  year = {2007},
  doi = {10.1007/s00220-007-0209-3}
}

@article{capitaine2009largest,
  author = {Capitaine, Mireille and Donati-Martin, Catherine and F{\'e}ral, Delphine},
  title = {The largest eigenvalues of finite rank deformation of large {Wigner} matrices: Convergence and nonuniversality of the fluctuations},
  journal = {The Annals of Probability},
  volume = {37},
  number = {1},
  pages = {1--47},
  year = {2009},
  doi = {10.1214/08-AOP394}
}

@article{benaychgeorges2011eigenvalues,
  author = {Benaych-Georges, Florent and Nadakuditi, Raj Rao},
  title = {The eigenvalues and eigenvectors of finite, low rank perturbations of large random matrices},
  journal = {Advances in Mathematics},
  volume = {227},
  number = {1},
  pages = {494--521},
  year = {2011},
  doi = {10.1016/j.aim.2011.02.007}
}

@article{lelarge2019fundamental,
  author = {Lelarge, Marc and Miolane, L{\'e}o},
  title = {Fundamental limits of symmetric low-rank matrix estimation},
  journal = {Probability Theory and Related Fields},
  volume = {173},
  pages = {859--929},
  year = {2019},
  doi = {10.1007/s00440-018-0845-x}
}

@inproceedings{ld-false-1,
  author = {Rares-Darius Buhai and Jun-Ting Hsieh and Aayush Jain and Pravesh K. Kothari},
  title = {The Quasi-Polynomial Low-Degree Conjecture is False},
  booktitle = {66th Annual Symposium on Foundations of Computer Science (FOCS)},
  pages = {2577--2590},
  publisher = {IEEE},
  year = {2025},
  doi = {10.1109/FOCS63196.2025.00134}
}

@article{ld-false-2,
  author = {Songtao Mao},
  title = {The Polynomial-Time Low-Degree Conjecture is False},
  journal = {arXiv preprint arXiv:2607.20318},
  year = {2026}
}

@book{mccullagh2018tensor,
  author = {McCullagh, Peter},
  title = {Tensor Methods in Statistics},
  publisher = {Dover Publications},
  year = {2018},
  isbn = {9780486823782},
  url = {https://www.stat.uchicago.edu/~pmcc/tensorbook/DoverEdition.pdf}
}

@article{gross1975logarithmic,
  author = {Gross, Leonard},
  title = {Logarithmic {Sobolev} Inequalities},
  journal = {American Journal of Mathematics},
  volume = {97},
  number = {4},
  pages = {1061--1083},
  year = {1975},
  doi = {10.2307/2373688}
}

@inproceedings{BR-reduction,
  author = {Berthet, Quentin and Rigollet, Philippe},
  title = {Complexity Theoretic Lower Bounds for Sparse Principal Component Detection},
  booktitle = {Proceedings of the 26th Annual Conference on Learning Theory},
  series = {Proceedings of Machine Learning Research},
  volume = {30},
  pages = {1046--1066},
  publisher = {PMLR},
  year = {2013},
  url = {https://proceedings.mlr.press/v30/Berthet13.html}
}

@inproceedings{HS17,
  author = {Hopkins, Samuel B. and Steurer, David},
  title = {Efficient {Bayesian} Estimation from Few Samples: Community Detection and Related Problems},
  booktitle = {58th Annual IEEE Symposium on Foundations of Computer Science (FOCS)},
  pages = {379--390},
  publisher = {IEEE},
  year = {2017},
  doi = {10.1109/FOCS.2017.42}
}

@inproceedings{HKP+17,
  author = {Hopkins, Samuel B. and Kothari, Pravesh K. and Potechin, Aaron and Raghavendra, Prasad and Schramm, Tselil and Steurer, David},
  title = {The Power of Sum-of-Squares for Detecting Hidden Structures},
  booktitle = {58th Annual IEEE Symposium on Foundations of Computer Science (FOCS)},
  pages = {720--731},
  publisher = {IEEE},
  year = {2017},
  doi = {10.1109/FOCS.2017.72}
}

@article{sos-clique,
  author = {Barak, Boaz and Hopkins, Samuel B. and Kelner, Jonathan and Kothari, Pravesh K. and Moitra, Ankur and Potechin, Aaron},
  title = {A Nearly Tight Sum-of-Squares Lower Bound for the Planted Clique Problem},
  journal = {SIAM Journal on Computing},
  volume = {48},
  number = {2},
  pages = {687--735},
  year = {2019},
  doi = {10.1137/17M1138236}
}

@inproceedings{CMSW25,
  author = {Chin, Byron and Mossel, Elchanan and Sohn, Youngtak and Wein, Alexander S.},
  title = {Stochastic block models with many communities and the {Kesten--Stigum} bound},
  booktitle = {Proceedings of the 38th Conference on Learning Theory},
  series = {Proceedings of Machine Learning Research},
  volume = {291},
  pages = {1253--1258},
  publisher = {PMLR},
  year = {2025},
  url = {https://proceedings.mlr.press/v291/chin25a.html}
}

@article{SGSW26,
  author = {Skeja, Anda and {Guti{\'e}rrez Espinoza}, Daniel and Skerman, Fiona and Wein, Alexander S.},
  title = {Sharp Low-Degree Thresholds for Planted-vs-Planted Testing},
  journal = {arXiv preprint arXiv:2606.05266},
  year = {2026}
}

@article{alg-contig,
  author = {Li, Zhangsong},
  title = {Computational Lower Bounds for Correlated Random Graphs via Algorithmic Contiguity},
  journal = {arXiv preprint arXiv:2502.09832},
  year = {2025}
}

@article{functional-ineq,
  author = {Li, Shuchen and Zadik, Ilias},
  title = {Functional Inequalities from the Fundamental Limits of {Bayesian} Estimation},
  journal = {arXiv preprint arXiv:2611.05831},
  year = {2026}
}
}

\end{document}